\PassOptionsToPackage{dvipsnames}{xcolor}

\documentclass[11pt]{article}
\usepackage{geometry}
\usepackage{tikz}
\usetikzlibrary{arrows.meta,positioning,calc,fit,decorations.pathreplacing}
\usepackage[normalem]{ulem}
\usepackage{a4wide}
\usepackage{hyperref}
\usepackage[english]{babel}
\usepackage[T1]{fontenc}
\usepackage[latin1]{inputenc}
\usepackage{dsfont}
\usepackage{mathrsfs}
\usepackage{amsmath}
\usepackage{bbm}
\usepackage{cleveref}
\usepackage{amssymb}
\usepackage{verbatim}
\usepackage{amsthm}
\usepackage{graphicx}
\usepackage{xcolor}
\usepackage{fancybox}
\usepackage{enumitem}
\usepackage{eurosym}
\usepackage{amsthm}

\hypersetup{
    colorlinks=true,
    linkcolor=Violet, 
    citecolor=OliveGreen, 
    urlcolor=brown
}
\usepackage[expansion=false]{microtype}
\allowdisplaybreaks
\newtheorem{theorem}{Theorem}[section]
\newtheorem{example}[theorem]{Example}
\newtheorem{prop}[theorem]{Proposition}
\newtheorem{lemma}[theorem]{Lemma}
\newtheorem{corollary}[theorem]{Corollary}
\newtheorem{definition}{Definition}
\newtheorem{remark}{Remark}

\newtheoremstyle{named}{}{}{\bfshape}{}{\bfseries}#1{\thmnote{#3 }}
\theoremstyle{named}

\newtheoremstyle{named1}{}{}{\bfshape}{}{\bfseries}#1 {\thmnote{#3 }}
\theoremstyle{named1}

\def\Q{\mathbb{Q}}

\def\B{\mathbf{B}}

\def\Tr{\mathrm{Tr}}

\newcommand{\R}{\mathbb{R}}

\newcommand{\1}{\mathbf{1}}

\newcommand{\Law}{\operatorname{Law}}

\newcommand{\Imm}{\operatorname{Im}}
\newcommand{\C}{\mathbb{C}}

\newcommand{\E}{\mathbb{E}}
\renewcommand{\i}{\sqrt{-1}}

\newcommand{\weak}{\Rightarrow}
\newcommand{\CCR}{\operatorname{CCR}}
\newcommand{\cS}{\mathcal{S}}
\newcommand{\cH}{\mathcal{H}}

\newcommand{\cD}{\mathcal{D}}

\newcommand{\cL}{\mathcal{L}}     
\newcommand{\cN}{\mathcal{N}}

\newcommand{\Id}{I}
\newcommand{\parallelsum}{\mathbin{\#}} 

\usepackage{stmaryrd}
\SetSymbolFont{stmry}{bold}{U}{stmry}{m}{n}
\usepackage{ytableau}
\providecommand{\Young}{\mathbb Y_{n,d}}
\providecommand{\ot}{\otimes}
\providecommand{\WSS}{\operatorname{WSS}}
\providecommand{\spec}{\operatorname{spec}}
\providecommand{\cU}{\mathcal U}
\providecommand{\cV}{\mathcal V}
\providecommand{\cF}{\mathcal F}
\providecommand{\ind}{\mathbf1}
\usepackage{cases}

\makeatother

\theoremstyle{plain}

\theoremstyle{remark}

\providecommand{\notename}{Note}
\providecommand{\theoremname}{Theorem}

\global\long\def\R{\mathbb{R}}%
\global\long\def\C{\mathbb{C}}%
\global\long\def\B{\mathcal{B}}%
\global\long\def\Tr{{\rm Tr}\,}%
\global\long\def\re{{\rm Re}\,}%
\global\long\def\conv#1{\stackrel{#1}{\rightsquigarrow}}%

\newlist{lamport}{enumerate}{4}
\setlist[lamport,1]{
  label=\(\langle 1\rangle\arabic*.\),
  ref=\(\langle 1\rangle\arabic*\),
  leftmargin=3.4em,
  labelsep=0.55em,
  itemsep=0.8em,
  topsep=0.6em
}
\setlist[lamport,2]{
  label=\(\langle 2\rangle\arabic*.\),
  ref=\(\langle 2\rangle\arabic*\),
  leftmargin=3.6em,
  labelsep=0.55em,
  itemsep=0.65em,
  topsep=0.45em
}
\setlist[lamport,3]{
  label=\(\langle 3\rangle\arabic*.\),
  ref=\(\langle 3\rangle\arabic*\),
  leftmargin=3.8em,
  labelsep=0.55em,
  itemsep=0.5em,
  topsep=0.35em
}
\setlist[lamport,4]{
  label=\(\langle 4\rangle\arabic*.\),
  ref=\(\langle 4\rangle\arabic*\),
  leftmargin=4em,
  labelsep=0.55em,
  itemsep=0.4em,
  topsep=0.3em
}

\title{Asymptotic efficiency of tests for quantum statistical models}
\author{
Rathindra Nath Karmakar\textsuperscript{*}\\
Graduate School of Mathematics, Kyushu University\\
\href{mailto:karmakar.rathindra.735@s.kyushu-u.ac.jp}{\texttt{karmakar.rathindra.735@s.kyushu-u.ac.jp}}
\and
Samriddha Lahiry\textsuperscript{*}\\
Department of Statistics and Data Science,\\ National University of Singapore\\
\href{mailto:slahiry@nus.edu.sg}{\texttt{slahiry@nus.edu.sg}}
\and 
Koichi Yamagata\textsuperscript{*}\\
Faculty of Electrical, Information and Communication Engineering,\\
Institute of Science and Engineering,\\
Kanazawa University\\
\href{mailto:yamagata@se.kanazawa-u.ac.jp}{\texttt{yamagata@se.kanazawa-u.ac.jp}}}
\date{}

\begin{document}

\maketitle

\begingroup
\renewcommand{\thefootnote}{\fnsymbol{footnote}}
\setcounter{footnote}{0}
\footnotetext{\textsuperscript{*}Authors are listed in alphabetical order of their surnames.}
\endgroup

\begin{abstract}
  In classical statistics, local asymptotic normality provides a general framework for characterizing optimal testing power and constructing asymptotically efficient procedures. Extending this framework to quantum statistical models introduces additional challenges arising from the noncommutativity of observables and the need to optimize over quantum measurements. In this paper, we develop a theory of asymptotic efficiency for quantum hypothesis testing using quantum local asymptotic normality and the asymptotic representation theorem. We establish sharp nuisance-maximin power bounds for one-sided tests and two-sided tests within the unbiased class by reducing the corresponding classical--quantum Gaussian shift experiments to classical testing problems. This reduction is justified through a quantum Hunt--Stein argument that accommodates both commuting and noncommuting Gaussian coordinates. Under suitable asymptotic invariance conditions, the resulting bounds transfer to the original quantum experiments and are attained by scalar symmetric logarithmic derivative tests. We further develop a two-stage adaptive construction that achieves the same local asymptotic power without knowledge of the reference parameter.

As an application, we study testing of smooth spectral
functionals of quantum states. We compare scalar symmetric
logarithmic derivative tests with the weak Schur
sampling plug-in test and with optimally weighted
linear statistics of its fluctuations. This comparison
distinguishes loss from the plug-in weighting from a
gap between the best such spectral statistic and the
scalar score procedure. We characterize equality of
the corresponding variances and illustrate the two
gaps through explicit examples.
\end{abstract}

\newpage

\setcounter{tocdepth}{2}
\renewcommand{\contentsname}{}
\tableofcontents

\section{Introduction}

Quantum statistics studies inferential questions arising in quantum-mechanical experiments and has become central to uncertainty quantification in emerging quantum technologies. Its objectives, including estimation and hypothesis testing, mirror those of classical statistics, although the underlying models differ. Quantum states, which play the role of probability distributions in classical statistics, are represented by operators that need not commute. Further, statistical procedures involve both the design of quantum measurements and the subsequent analysis of the resulting classical data. Consequently, characterizing optimal procedures generally requires joint optimization over measurements and decision rules.

The decision-theoretic formulation of quantum inference dates back to the foundational work of Helstrom and Holevo \cite{Helstrom,Hol-2nd-ed-11}. The rapid development of quantum technologies has renewed interest in this framework. Over the last two decades, quantum statistical methods have become increasingly important in applications such as quantum state tomography, quantum imaging, quantum metrology, and quantum sensing. These developments have raised fundamental questions concerning the limits of statistical inference in quantum systems, including the determination of efficiency bounds and minimax risks and the construction of procedures that attain them.

\subsection{Inference in quantum statistical models}

A quantum statistical model is a family
\(
\mathcal Q=\{\rho_\theta:\theta\in\Theta\},
\)
where the parameter space $\Theta$ is a subset of a Euclidean space and each $\rho_\theta$ is a density operator, that is, a positive semidefinite Hermitian operator of unit trace. We briefly outline the principal inferential problems arising in such models.

Given measurement outcomes from repeated preparations of a
quantum system, one seeks to test statistical hypotheses concerning
the underlying state or its unknown parameters.
The hypotheses
may also concern structural properties of the state, such as purity,
separability, entanglement, low-rank structure, or membership in a
specified family. Testing such properties without reconstructing the
full state has attracted considerable attention in quantum information
theory \cite{montanaro2016survey}.

A standard example of a quantum hypothesis-testing problem is quantum
state discrimination, in which the objective is to distinguish between
specified quantum states through a suitable measurement and decision
rule. In the asymptotic setting of many independent, identically
prepared copies, a central question concerns the optimal exponential
rates at which testing errors decay. For discrimination between two
fixed states, the quantum Stein lemma and the quantum Chernoff bound
characterize these rates in asymmetric and symmetric testing,
respectively
\cite{MR1806811,Audenaert,Audenaert&Nussbaum,Li,Nussbaum&Szkola}.
More broadly, quantum decision theory provides a framework for studying
the choice of measurements and decision rules, extending the classical
theory of statistical experiments and decision procedures to
noncommutative statistical models \cite{Helstrom,Hol-2nd-ed-11}.

Within this broader framework, the present paper considers hypotheses
about unknown parameters, with an emphasis on the asymptotic efficiency
of testing procedures. At a prescribed significance level, a central
criterion is local power, the ability to detect alternatives that
approach the null hypothesis as the number of repeated preparations
increases. This perspective assesses the sensitivity of tests to small
departures from the null and provides a basis for comparing their
efficiency.

The study of local power is naturally related to quantum parameter
estimation, since the precision of an estimator can govern the
sensitivity of a test constructed from it. Quantum analogues of the
classical Cram\'er--Rao inequality provide lower bounds on the covariance
matrices of locally unbiased estimators and serve as benchmarks for
attainable precision \cite{Braunstein&Caves,Helstrom,Hol-2nd-ed-11}.
These bounds underpin the study of asymptotic efficiency in quantum
estimation and play a central role in quantum metrology, imaging,
and sensing. As in classical statistics, the connection between
estimation and testing operates in both directions. Testing bounds
can yield lower bounds on estimation risk, while asymptotically
normal estimators can be used, under suitable regularity conditions
and with consistent covariance estimation, to construct Wald-type
tests. We study optimal testing through local asymptotic normality, building on its use in establishing the efficiency of quantum estimators \cite{Fujiwara2023}.


\subsection{Local asymptotic normality in quantum statistical models}

In classical statistics, one approach to establishing optimality is to approximate a complicated statistical experiment by a simpler limit experiment in the sense of Le Cam \cite{LeCam}. An important instance is \emph{local asymptotic normality} (LAN). Under LAN, a regular i.i.d.\ model localized around a fixed parameter value has a Gaussian shift limit. Such approximations allow optimality bounds derived in the Gaussian limit experiment to be carried over to the original model. A procedure whose asymptotic risk or power attains the corresponding bound is then asymptotically optimal in the original problem.

Beyond local approximations, asymptotic equivalence has also been established for a variety of nonparametric models. Under suitable smoothness conditions, both nonparametric regression \cite{BrownLow} and density estimation \cite{Nussbaum96} are asymptotically equivalent to Gaussian white noise experiments. These equivalence results provide a unified framework for deriving minimax rates, efficiency bounds, and optimal procedures across seemingly different statistical models. LAN techniques have also proved useful in the study of functional estimation, where suitably regular functionals can often be analyzed through associated local parametric submodels \cite{khasminskiibook}.

The quantum analogue of this theory, known as \emph{quantum local asymptotic normality} (q-LAN), was developed in \cite{GutaJencova2007,Kahn&Guta,YFG13,Fujiwara2020}. These results show that, under suitable regularity conditions, a collection of \(n\) independent, identically prepared finite-dimensional quantum systems can be approximated locally around a fixed state by a quantum Gaussian shift experiment. A quantum analogue of nonparametric asymptotic equivalence was established in \cite{BGN-QAE}, where local ensembles of pure states in an infinite-dimensional Hilbert space were shown to be asymptotically equivalent to coherent-state models, providing a quantum counterpart of the classical Gaussian white noise experiment. Related asymptotic approximations have also been established for quantum Markov models \cite{guta2015equivalence,guta2017information,girotti2026asymptotic} and quantum time series models \cite{nussbaum2025asymptotic}.
Quantum LAN has become a fundamental tool for deriving asymptotic lower bounds and constructing procedures that attain them \cite{Fujiwara2023,lowLAN,minimaxpure,GutaKahnproc}. Much of these works concern parameter estimation. Here we use q-LAN to establish analogous optimality results for hypothesis testing.

Three related approaches to Gaussian approximation are relevant here. One uses a local expansion of the log-likelihood ratio in terms of a central sequence $\Delta_n$. The others study strong convergence in Le Cam distance and weak convergence of statistical experiments. Quantum counterparts of these approaches were developed in \cite{YFG13}, \cite{Kahn&Guta}, and \cite{GutaJencova2007}, respectively. We use the likelihood-ratio formulation of \cite{YFG13}, together with the quantum asymptotic representation theorem and the quantum analogue of Le Cam's third lemma. The representation theorem provides the comparison with the Gaussian experiment, and the third lemma gives the limiting laws of our test statistics under local alternatives. Their classical counterparts play the same roles in the analysis of testing efficiency in Chapter~15 of \cite{vaart}.

\subsection{Efficiency of tests in classical statistics}

In the classical setup, Le Cam's asymptotic representation theorem \cite[Theorem~7.10]{vaart} provides a link between the original statistical experiment and its Gaussian limit. We give the precise statement in Appendix~\ref{sec:class-lan-testing}. Roughly speaking, if a sequence of experiments converges weakly to a limit experiment and a sequence of statistics has limiting distributions under local alternatives, then a single, possibly randomized procedure in the limit experiment realizes these distributions.

This observation underlies several results in asymptotic optimality theory. In estimation, it leads to the local asymptotic minimax theorem, which gives a lower bound on local asymptotic risk in terms of the minimax risk in the Gaussian shift experiment. In testing, the same principle yields an upper bound on local power. For a one-sided hypothesis about a differentiable scalar functional, Theorem~15.4 of \cite{vaart} bounds the asymptotic power of any sequence of level-$\alpha$ tests at each fixed local alternative by the power of the corresponding Gaussian test. Addendum~15.5 gives a sufficient condition for attainment. 

Fujiwara and Yamagata \cite{Fujiwara2023} proved a quantum version of the asymptotic representation theorem and used it to establish efficiency bounds for estimators in quantum models. Our goal is to develop an analogous framework for quantum hypothesis testing. In particular, we use q-LAN to derive asymptotically optimal tests by analyzing the corresponding quantum Gaussian limit experiment.

\subsection{Our contributions}

Our main contributions are as follows.

\begin{enumerate}
\item We first solve the problem of testing a linear functional of the mean in a quantum Gaussian model with known covariance. A shift of the mean that leaves the linear functional
unchanged preserves the hypotheses. Using this symmetry, we show
that restricting attention to tests invariant under these nuisance
shifts preserves the greatest attainable worst-case power. The
argument builds on the quantum Hunt--Stein method used by Kumagai
and Hayashi for Gaussian tests with nuisance parameters
\cite{KumagaiHayashi2013}. In the present setting, every invariant test can be
implemented by measuring the linear combination of Gaussian coordinates
corresponding to the functional and applying a classical rule to its
outcome. The outcome has a classical Gaussian distribution, so the
upper-tail test of its mean gives the sharp one-sided power bound, even when competing tests can
use arbitrary quantum measurements. The two-tailed normal test gives
the corresponding bound within the class of unbiased tests.

\item We next transfer these bounds to q-LAN models whose central sequences are asymptotically $\mathcal D$-invariant. The asymptotic representation theorem relates the power of tests under local alternatives to the power of tests in the corresponding quantum Gaussian experiment. Applying the Gaussian bound to the local linear approximation of a differentiable functional gives an upper bound on worst-case local power. For two-sided hypotheses, this comparison is within the unbiased class. In the independent-copy model, a test based on a weighted combination of the symmetric logarithmic derivatives, the quantum analogues of the usual likelihood scores, attains this bound in local asymptotic power.

\item The weighted combination of observables on which the optimal test is based, typically depends on the unknown state, even when
the null value of the functional is known. For smooth independent-copy
models of constant rank, we give a two-stage adaptive procedure. A preliminary
measurement on a vanishing fraction of the copies provides an estimate
of the parameter. We use that estimate to choose the observable
measured separately on each remaining copy. The average of the
resulting observations supplies a one-step correction to the preliminary
estimate of the functional. The resulting Wald tests have the same
limiting power as the test designed at the true reference point,
uniformly over compact sets of local perturbations. This construction
is related to Fujiwara's work on adaptive quantum estimation
\cite{adaptive}.

\item We apply our results to testing smooth spectral
functionals at states with distinct positive eigenvalues
and a nonvanishing functional derivative. A natural
measurement for this problem is \emph{weak Schur sampling}
(WSS), which originates in the study of symmetries of
product states
\cite{AlickiRudnickiSadowski1988,KeylWerner2001}
and has been used in spectrum testing and entropy
estimation
\cite{ChildsHarrowWocjan2007,ODonnellWright2021,AcharyaIssaShendeWagner2020}.
The sufficiency result of \cite{Lahiry2026QCS} shows
that WSS retains all information available through
unitarily invariant measurements, although this need
not imply optimality in restricted parametric models.
We identify two distinct sources of efficiency loss:
suboptimal processing of the WSS outcome and information
discarded by the measurement itself. In particular, we
characterize the best limiting power among linear
statistics of the WSS fluctuations having the correct
first-order mean in every local direction, and show how
this improves upon the usual plug-in test. We then
compare this power with that of the scalar score test,
which performs at least as well and attains the local
maximin bound under our q-LAN assumptions. We derive
necessary and sufficient conditions for equality in both
comparisons. These conditions have a geometric
interpretation in terms of the structure of the tangent space,
closely related to the projection-based approach to
quantum semiparametric estimation developed in
\cite{tsang_semipara}.

\end{enumerate}

\subsection{Organization of the paper}

Section~\ref{sec:prelim} introduces the notation and background on quantum statistical models. We define q-LAN, present the version of the quantum asymptotic representation theorem used in our optimality proofs, and compare classical and quantum hypothesis testing. Section~\ref{sec:gauss} gives the optimal power bounds in quantum Gaussian models, and Section~\ref{sec:iid_opt} derives the corresponding bounds for q-LAN models, including the i.i.d.\ case. Section~\ref{sec:wald} develops an adaptive version of the classical Wald test that attains the same limiting power while Section~\ref{sec:spec_test} applies these results to tests of spectral functionals. Section~\ref{sec:disc} discusses extensions of our framework. The proofs and auxiliary results are collected in the appendices.

\subsection{Notation}
For a vector $v$ in a separable Hilbert space $\mathcal H$, we write $|v\rangle$ for the vector (the ket) and $\langle v|$ for the corresponding element $v^*$ of the dual space $\mathcal H^*$ (the bra). The inner product $\langle u|v\rangle\in\mathbb C$ is linear in its second argument and conjugate-linear in its first. The rank-one operator $M:=|u\rangle\langle v|$ acts by $M|w\rangle=\langle v|w\rangle\,|u\rangle$. For an operator $A$, we also write $\langle u|A|v\rangle$ for $\langle u|Av\rangle$. The space of bounded linear operators on $\mathcal H$ is denoted by $\mathcal L(\mathcal H)$. We write $\mathcal T_1(\mathcal H)\subset\mathcal L(\mathcal H)$ for the trace-class operators, equipped with the trace norm
\[
\|A\|_1:=\Tr\!\left[(A^\dagger A)^{1/2}\right].
\]

We write $\|\cdot\|$ for the Hilbert space norm and $\mathbf{1}$ for the identity operator, with the underlying space understood from context. For probability measures $\mu$ and $\nu$, we use $\|\mu-\nu\|_{\mathrm{TV}}:=\sup_B|\mu(B)-\nu(B)|$, where the supremum is over measurable sets $B$. For probability mass functions $p$ and $q$ on a countable space, $\|p-q\|_1:=\sum_x|p(x)-q(x)|=2\|p-q\|_{\mathrm{TV}}$. We write $a\vee b$ and $a\wedge b$ for $\max(a,b)$ and $\min(a,b)$, respectively, and $a_+:=a\vee 0$. The symbols $\lfloor a\rfloor$ and $\lceil a\rceil$ denote the largest integer less than or equal to $a$ and the smallest integer greater than or equal to $a$, respectively. We write $a_n\asymp b_n$ if $c<\liminf_n(a_n/b_n)\leq\limsup_n(a_n/b_n)<C$ for some constants $c,C>0$. Throughout, $c$ and $C$ denote generic constants. Finally, $\Phi$ denotes the standard Gaussian distribution function, and $\Phi(a,b)$ denotes the probability that a standard Gaussian random variable lies in $[a,b]$.






\section{Preliminaries}\label{sec:prelim}

\subsection{States, Measurements and Observables}
A \textit{state of a quantum system} is described by a self-adjoint operator
$\rho$ on a complex Hilbert space {$\mathcal{H}$, which is positive }%
$(\rho\geq0$) and normalized to $\mathrm{Tr}\left(  \rho\right)  =1$ (a
density operator). A state is called \textit{pure} if it is of the form $\rho=|\psi\rangle\langle\psi|$, otherwise it is called a \textit{mixed state}. We denote the set of states by $\mathcal{S}(\mathcal{H})$. It can be shown that $\mathcal{S}(\mathcal{H})\subset \mathcal{T}_1(\mathcal{H})$.

Data on a quantum system are obtained from
\textit{observables} which are self-adjoint operators in the Hilbert space $\mathcal{H}$. If  $S$ is a self-adjoint operator in {$\mathcal{H}$
with spectral decomposition }$S=\sum_{j}\lambda_{j}\Pi_{j}$ where $\Pi_{j}$
are projectors, then a measurement generates a \textit{discrete random variable} $X_{S}$ taking
values in the set of eigenvalues $\left\{  \lambda_{1},\lambda_{2}%
,\ldots\right\}  $ with probabilities $p_{j}=\mathrm{Tr}(\rho\cdot\Pi_{j})$. Whenever the first absolute moment is finite, the expectation of
\(X_S\) under the state \(\rho\) is given by the
\textit{Born-von Neumann postulate:}%
\begin{equation}
E_{\rho}X_{S}=\sum_{j}\lambda_{j}\mathrm{Tr}\left[  \rho\Pi_{j}\right]
=\mathrm{Tr}\left[  \rho S\right]  . \label{trace-rule}%
\end{equation}

More generally, a measurement with outcomes in a measurable space
$(\Omega,\mathfrak B)$ is described by a positive operator-valued
measure (POVM).
\begin{definition}\label{def.POVM}
A positive operator valued measure (POVM) is a map $M:\mathfrak{B}\to \mathcal{L}(%
\mathcal{H})$ having the following properties
\begin{itemize}
\item[1)] positivity: $M(B) \geq 0$ for all events $B\in\mathfrak{B}$ (hence M(B) is self-adjoint)
\item[2)] $\sigma$-additivity: $M(\cup_i B_i) = \sum_i {M}(B_i)$ for
any countable set of mutually disjoint events $B_i$ (here the convergence is in the weak operator topology of $\mathcal{L}(\mathcal{H})$)
\item[3)] normalization: $M(\Omega) = \mathbf{1}$.
\end{itemize}
\end{definition}

The simplest example of a POVM is given by an \emph{effect} $E$, which is a bounded operator satisfying $0\leq E\leq \mathbf{1}$. Indeed $\{E,\mathbf{1}-E\}$ forms a two-outcome POVM. If the operators $M(B)$ are also orthogonal projections, i.e. $M(A)^2=M(A)$ and $M(B)M(A)=0$ when $A\cap B=\emptyset$, then it is called a \textit{simple measurement}. The collection of projectors $\{\Pi_j\}$ in the spectral decomposition $S=\sum_j \lambda_j\Pi_j$ is an example of a simple measurement. The outcome of the measurement has probability
distribution
\begin{equation}
P_\rho(B) = \mathrm{Tr}(\rho M(B)), \qquad B\in \mathfrak{B}.
\label{prob_measurement}
\end{equation}

The spectral theorem associates with every self-adjoint operator $S$,
with dense domain $D(S)\subset\mathcal H$, a projection-valued spectral
measure $M$ such that
\[
    S=\int_{\sigma(S)}x\,M(dx),
\]
in the sense of the spectral calculus. Here $\sigma(S)$ is the spectrum of $S$ and $M$ is a POVM, also called spectral measure associated with the operator $S$. When $S$ is an observable with a continuous spectrum, it generates a \textit{continuous random variable} $X_S$ with probabilities given by (\ref{prob_measurement}). 
Whenever the first absolute moment is finite,
\[
    E_\rho[X_S]
    =
    \int_{\sigma(S)} x\,\Tr\{\rho M(dx)\},
\]
which we also denote by $\Tr(\rho S)$.

\subsection{Quantum Gaussian states and convergence in distribution}
\label{subsec:gaussian-states}

Here we introduce the Gaussian models which arise as limiting
experiments under q-LAN. We begin with the usual Hilbert-space
representation of quantum Gaussian states and subsequently explain
how classical Gaussian coordinates are incorporated when the
commutation matrix is degenerate.

To define a one-mode Gaussian state, consider the Hilbert space
$\mathcal H=L^2(\mathbb R)$ and the position and momentum operators
\[
    (Qf)(x)=xf(x),
    \qquad
    (Pf)(x)=-i\frac{d}{dx}f(x),
\]
defined on appropriate dense domains. They satisfy the canonical
commutation relation $[Q,P]=iI$. For $u=(u_1,u_2)^T\in\mathbb R^2$, the linear combination
$Z_u=u_1Q+u_2P$ is a self-adjoint observable. Given a density operator $\rho$ on
$\mathcal H$, its quantum characteristic function is
\[
    \widetilde W_\rho(u)
    =
    \operatorname{Tr}\!\left(
        \rho e^{iZ_u}
    \right).
\]

More generally, for an $m$-mode system we consider
$\mathcal H_m=\bigotimes_{j=1}^mL^2(\mathbb R)$
and the canonical observables
\[
    R=(Q_1,P_1,\ldots,Q_m,P_m)^T.
\]
Here $Q_j$ and $P_j$ act on the $j$th tensor factor and satisfy
\begin{equation}
    [Q_j,P_k]=i\delta_{jk}I,
    \qquad
    [Q_j,Q_k]=[P_j,P_k]=0.
    \label{joint_commutation}
\end{equation}
Writing $R=(R_1,\ldots,R_{2m})^T$, these relations may be expressed as
\begin{equation}
    [R_j,R_k]=2iS_{jk}I,
    \label{eq:gaussian-CCR}
\end{equation}
where
\[
    S
    =
    \frac12
    \bigoplus_{j=1}^m
    \Omega,
    \qquad
    \Omega
    =
    \begin{pmatrix}
        0&1\\
        -1&0
    \end{pmatrix}.
\]
The matrix $S$ is real, skew-symmetric, and nonsingular. For
$\xi\in\mathbb R^{2m}$, define the quantum characteristic function by
\[
    \widetilde W_\rho(\xi)
    :=
    \operatorname{Tr}\left(
        \rho e^{i\xi^TR}
    \right).
\]

\begin{definition}[Quantum Gaussian state]
\label{def:quantum-gaussian-state}
Let $\mu\in\mathbb R^{2m}$ and let $\Sigma=V+iS$, where $V$ is a real
symmetric matrix and $\Sigma\succeq0$. A density operator $\rho$ on
$\mathcal H_m$ is called a \emph{quantum Gaussian state} with mean
$\mu$ and covariance matrix $\Sigma$ if
\begin{equation}
    \operatorname{Tr}\left(
        \rho e^{i\xi^TR}
    \right)
    =
    \exp\left(
        i\xi^T\mu-\frac12\xi^TV\xi
    \right),
    \qquad
    \xi\in\mathbb R^{2m}.
    \label{eq:quantum-gaussian-cf}
\end{equation}
We use the notation $\rho=\mathbb N_{2m}(\mu,\Sigma)$.
\end{definition}

The real part $V=\operatorname{Re}\Sigma$ is the symmetrized covariance
matrix, while $S=\operatorname{Im}\Sigma$ encodes the commutation
relations. With the convention
\[
[R_j,R_k]=2iS_{jk}I,
\qquad
\Sigma=V+iS,
\]
the matrix \(\Sigma\) is also the ordered second-moment matrix of the
centered Gaussian observables. More precisely,
\begin{equation}
\Tr\!\left[
\rho\,
(R_j-\mu_j I)(R_k-\mu_k I)
\right]
=
V_{jk}+iS_{jk}
=
\Sigma_{jk}.
\label{eq:ordered-gaussian-covariance}
\end{equation}
Indeed, the symmetric part of the left-hand side is \(V_{jk}\), while
\[
\frac12
\Tr\!\left[
\rho\,[R_j,R_k]
\right]
=
iS_{jk}.
\]
Further since $S$ is skew-symmetric, $a^T\Sigma a=a^TVa$ for every
$a\in\mathbb R^{2m}$. It follows immediately from the characteristic
function that measuring the single observable $a^TR$ under
$\mathbb N_{2m}(\mu,\Sigma)$ produces the ordinary Gaussian distribution
\begin{equation}
    a^TR
    \sim
    N(a^T\mu,a^TVa).
    \label{eq:scalar-quantum-gaussian}
\end{equation}

\paragraph{Gaussian models with classical coordinates.}

The Gaussian limits arising in q-LAN need not have a nonsingular
commutation matrix. We therefore allow a general real skew-symmetric
matrix $S\in\mathbb R^{d\times d}$. For a $k$-mode quantum system,
consider the canonical observables $\{Q_i,P_i\}_{i=1}^k$ satisfying
\eqref{joint_commutation}. In addition, allow $\ell$ classical Gaussian
random variables $C_1,\ldots,C_\ell$, which commute with each other and
with all the quantum coordinates. We write the resulting
$d=2k+\ell$ Gaussian coordinates as
\begin{equation}
    \label{hybrid_notation}
    (X_1,\ldots,X_d)
    \equiv
    (C_1,\ldots,C_\ell,Q_1,P_1,\ldots,Q_k,P_k).
\end{equation}
Their commutation structure is represented compactly by
$[X_i,X_j]=2iS_{ij}I$, where
\[
    S
    =
    0_{\ell\times\ell}
    \oplus
    \frac12\bigoplus_{j=1}^k\Omega.
\]

We represent a product classical--quantum experiment by
$\varrho=f\otimes\rho$, where $f$ is a probability density on
$\mathbb R^\ell$ and $\rho$ is a density operator on
$\mathcal H_k=L^2(\mathbb R)^{\otimes k}$. Its characteristic function
is defined by
\begin{equation}
    \label{hybrid_cf}
    E_{\varrho}
    \left(
        e^{i\sum_{j=1}^d u_jX_j}
    \right)
    :=
    \int_{\mathbb R^\ell}
    \operatorname{Tr}
    \left(
        \rho
        e^{\,i\sum_{j=\ell+1}^{d}u_jX_j}
    \right)
    f(y)
    e^{\,i\sum_{j=1}^{\ell}u_jy_j}
    \,dy.
\end{equation}

\begin{definition}[Classical--quantum Gaussian experiment]
A product classical--quantum experiment $\varrho=f\otimes\rho$ is
called Gaussian if its characteristic function has the form
\[
    E_{\varrho}
    \left(
        e^{i\sum_{j=1}^d u_jX_j}
    \right)
    =
    \exp\left(
        iu^T\mu-\frac12u^TVu
    \right),
    \qquad
    u\in\mathbb R^d,
\]
for some $\mu\in\mathbb R^d$ and real symmetric covariance matrix
$V$, and we associate to it the complex covariance
$\Sigma=V+iS$.
\end{definition}

A Gaussian experiment with a possibly degenerate commutation
matrix and $\operatorname{Re}\Sigma\succ0$ admits, after an
invertible real linear transformation of its coordinates,
a product classical--quantum representation. We use
$N(\mu,\Sigma)$ for both purely quantum and classical--quantum
Gaussian experiments. When $S=\operatorname{Im}\Sigma$ is
nonsingular, this experiment is represented by the Gaussian
density operator $\mathbb N_d(\mu,\Sigma)$ introduced above;
when $S$ is degenerate, it denotes the corresponding
classical--quantum Gaussian experiment. In either case, we
write $\phi_{\mu,\Sigma}$ for the associated normal state
on the von Neumann algebra generated by the Gaussian Weyl
operators. The precise von Neumann algebra formulation,
including the degenerate case, is given in
Appendix~\ref{sec:vn_alg}.

It is convenient to denote the Gaussian coordinates collectively by
$X=(X_1,\ldots,X_d)^T$. In the classical--quantum case, some of the
$X_j$ are ordinary classical Gaussian random variables and the
remaining coordinates are canonical quantum observables. For every
$a\in\mathbb R^d$, the scalar $a^TX$ has the ordinary Gaussian
distribution
\begin{equation}
    a^TX
    \sim
    N\left(
        a^T\mu,
        a^T\operatorname{Re}(\Sigma)a
    \right).
    \label{eq:scalar-general-gaussian}
\end{equation}

\paragraph{Convergence in distribution.}

Let $\mathcal H^{(n)}$ be a sequence of Hilbert spaces,
$\rho^{(n)}$ states on $\mathcal H^{(n)}$, and
\[
    R^{(n)}
    =
    (R_1^{(n)},\ldots,R_d^{(n)})
\]
lists of observables. We say that
\[
    R^{(n)}
    \stackrel{\rho^{(n)}}{\leadsto}
    N(\mu,\Sigma)
\]
if their ordered quantum characteristic functions converge to those
of the limiting Gaussian experiment. More precisely, for every
$r\ge1$ and every
$\xi_1,\ldots,\xi_r\in\mathbb R^d$,
\begin{equation}
    \lim_{n\to\infty}
    \operatorname{Tr}\left[
        \rho^{(n)}
        \prod_{t=1}^r
        e^{i\xi_t^TR^{(n)}}
    \right]
    =
    \mathbb E_{N(\mu,\Sigma)}
    \left[
        \prod_{t=1}^r
        e^{i\xi_t^TX}
    \right].
    \label{eq:quantum-convergence-distribution}
\end{equation}
When the limiting model is purely quantum, the expectation on the
right-hand side is the corresponding trace against the Gaussian
density operator. When classical coordinates are present, the
expectation additionally averages over the classical Gaussian
component.

\subsection{q-LAN and D-Extendibility}
\label{subsec:qLAN-Dext}

Next we formalize this using the q-LAN framework \cite{Fujiwara2020, Fujiwara2023}, written here in the functional notation of Kumagai and Hayashi \cite{KumagaiHayashi2013}. In classical statistics, an experiment is defined to be a family of probability measures on a sample space and denoted by $\mathcal{E}=\{P_{\theta},\theta\in \Theta\}$ where $\Theta$ is the parameter space. Analogously, we define a quantum statistical model as follows:
\begin{definition}
A quantum statistical model over a parameter space $\Theta$ consists of a
family of quantum states $\mathcal{Q} = \{\rho_\theta :\, \theta\in \Theta \}
$ on a Hilbert space $\mathcal{H}$, indexed by an unknown parameter $%
\theta\in \Theta$.
\end{definition} 

For a differentiable parametric family of states
$
\mathcal S
=
\{\rho_\theta:\theta\in\Theta\subset\mathbb R^d\},
$
the \emph{symmetric logarithmic derivative} (SLD) in the \(j\)th
parameter direction at \(\theta\) is a self-adjoint operator
\(L_{\theta,j}\) satisfying
\begin{equation}
\partial_j\rho_\theta
=
\frac12
\left(
\rho_\theta L_{\theta,j}
+
L_{\theta,j}\rho_\theta
\right),
\qquad j=1,\ldots,d.
\label{eq:general-sld}
\end{equation}
Taking the trace in
\eqref{eq:general-sld} also gives
\begin{equation}
\Tr(\rho_\theta L_{\theta,j})=0.\label{eq:trace_SLD}\end{equation}
When \(\rho_\theta\) is strictly positive, the SLD is unique. Indeed,
in an eigenbasis of \(\rho_\theta\), with eigenvalues
\(\lambda_1,\ldots,\lambda_r>0\),
\[
(L_{\theta,j})_{ab}
=
\frac{2(\partial_j\rho_\theta)_{ab}}
{\lambda_a+\lambda_b}.
\]

Whenever the SLDs exist, we define their ordered second-moment matrix
and the real SLD Fisher information matrix by
\begin{equation}
J_{\theta,jk}
:=
\Tr\!\left(
\rho_\theta
L_{\theta,j}L_{\theta,k}
\right),
\qquad
G_{\theta,jk}
:=
\operatorname{Re}J_{\theta,jk}
=
\frac12
\Tr\!\left[
\rho_\theta
\left(
L_{\theta,j}L_{\theta,k}
+
L_{\theta,k}L_{\theta,j}
\right)
\right].
\label{eq:general-sld-information}
\end{equation}
Thus \(G_\theta=\operatorname{Re}J_\theta\) is real symmetric and
positive semidefinite. 

Next for finite-dimensional states $\rho$ and $\sigma$, we consider the
square-root likelihood ratio $\mathcal R(\sigma\mid\rho)$ and the
quantum Lebesgue decomposition of \cite{Fujiwara2020}. A chosen
representative\\$R=\mathcal R(\sigma\mid\rho)$ satisfies
\[
 \sigma=R\rho R+\sigma^\perp,
 \qquad R\succeq0,\quad \sigma^\perp\succeq0,
 \quad \operatorname{Tr}(\rho\sigma^\perp)=0.
\]
When $\rho$ and $\sigma$ are faithful, the positive solution is unique
and is
\[
 R=\rho^{-1/2}
       (\rho^{1/2}\sigma\rho^{1/2})^{1/2}\rho^{-1/2}.
\]
At nonfaithful states the displayed decomposition alone does not
assert uniqueness on the ambient Hilbert space.
In the next definition we also use the notations $o_{L^2}(\rho^{(n)})$ and
$o_D(Z^{(n)},\rho^{(n)})$ to denote remainders
infinitesimal in $L^2$ and in distribution,
respectively, following
\cite[Definition~1.2]{Fujiwara2023};
see also \cite[Definition~5.5]{Fujiwara2020}
for the precise definition of the latter.

\begin{definition}[q-LAN]
A sequence of quantum statistical models $\cS^{(n)}=\{\rho_{\theta}^{(n)} \mid \theta\in\Theta\subseteq \R^d\}$ is called \textbf{q-LAN} at $\theta_0$ if for every fixed local parameter $h\in\R^d$, the square-root likelihood ratio $\mathcal{R}_h^{(n)} = \mathcal{R}(\rho^{(n)}_{\theta_0+h/\sqrt{n}}\mid \rho^{(n)}_{\theta_0})$ can be expanded  as:
\begin{equation}
    \log\left(\mathcal{R}^{(n)}_h + o_{L^2}(\rho^{(n)}_{\theta_0})\right)^2 = h^\top \Delta^{(n)} - \frac{1}{2}(h^\top J h)\Id^{(n)} + o_D(h^\top \Delta^{(n)}, \rho^{(n)}_{\theta_0})
\end{equation}
where
\[
\Delta^{(n)}
=
\left(
\Delta_1^{(n)},\ldots,\Delta_d^{(n)}
\right)^T
\]
is a \(d\)-tuple of self-adjoint observables, called the
\emph{quantum central sequence}, satisfying
\[
\Delta^{(n)}
\stackrel{\rho^{(n)}_{\theta_0}}{\leadsto}
\cN(0,J),
\]
where $J\in\mathbb C^{d\times d}$ is Hermitian positive semidefinite
and $\operatorname{Re}J\succ0$.
\end{definition}

\begin{definition}[Holevo's commutation operator]
\label{def:commutation-operator}
Let $\rho$ be a density matrix on a finite-dimensional
Hilbert space, and define the left- and right-multiplication
superoperators by
\[
\mathsf L_\rho(X):=\rho X,
\qquad
\mathsf R_\rho(X):=X\rho.
\]
When $\rho$ is faithful, Holevo's commutation operator is
the linear superoperator
\[
\cD_\rho
:=
\i(\mathsf L_\rho+\mathsf R_\rho)^{-1}
(\mathsf L_\rho-\mathsf R_\rho).
\]
Equivalently, it is characterized by
\begin{equation}
\rho\cD_\rho(X)+\cD_\rho(X)\rho
=
\i(\rho X-X\rho).
\label{eq:commutation-operator}
\end{equation}
For a possibly nonfaithful state, we use the canonical
extension obtained by replacing the inverse with the
Moore--Penrose inverse. In an eigenbasis of $\rho$, with
eigenvalues $\lambda_1,\ldots,\lambda_r$, this gives
\[
\bigl(\cD_\rho(X)\bigr)_{ab}
=
\begin{cases}
\displaystyle
\i\frac{\lambda_a-\lambda_b}
        {\lambda_a+\lambda_b}X_{ab},
& \lambda_a+\lambda_b>0,\\[2mm]
0,&\lambda_a+\lambda_b=0.
\end{cases}
\]
For nonfaithful $\rho$, this is the canonical representative
of the operator on the quotient by
\[
\mathcal K_\rho
:=
\{X:X\rho=\rho X=0\};
\]
see \cite[Appendix~C]{Fujiwara2023}.
\end{definition}

\begin{remark}[$\cD$-invariance]
For a finite-dimensional quantum statistical model, let
\[
\mathcal T_\rho
=
\operatorname{span}_{\R}\{L_1,\ldots,L_d\}
\]
denote the real SLD tangent space. When $\rho$ is faithful,
this space is called $\cD_\rho$-invariant if
\[
\cD_\rho(\mathcal T_\rho)
\subseteq\mathcal T_\rho.
\]
When this condition fails, one may enlarge the SLD tangent
space by introducing additional observables so that the
extended space is $\cD_\rho$-invariant.
\end{remark}

For the quantum asymptotic representation theorem,
we use the following asymptotic formulation of
$\cD$-extendibility, following
\cite[Definition~2.2]{Fujiwara2023}.

\begin{definition}[$\cD$-extendibility]
\label{def:D-extendibility}
Let
$R^{(n)}=(R_1^{(n)},\ldots,R_m^{(n)})^T$
be a list of self-adjoint observables. The sequence
$R^{(n)}$ is called \emph{asymptotically
$\cD$-invariant} at $\theta_0$ if
\[
R^{(n)}
\stackrel{\rho^{(n)}_{\theta_0}}{\leadsto}
\cN(0,\Sigma),
\]
where $\Sigma\in\mathbb C^{m\times m}$ is Hermitian,
$\Sigma\succeq0$, and $\re\Sigma\succ0$, and if,
for every $\xi,\eta\in\R^m$,
\begin{equation}
\begin{aligned}
&\lim_{n\to\infty}
\Tr\left(
\sqrt{\rho^{(n)}_{\theta_0}}\,
e^{\i\xi^\top R^{(n)}}
\sqrt{\rho^{(n)}_{\theta_0}}\,
e^{\i\eta^\top R^{(n)}}
\right)\\
&\quad=
\exp\left[
-\frac12
\begin{pmatrix}\xi\\\eta\end{pmatrix}^{\top}
\begin{pmatrix}
\Sigma & \Sigma\parallelsum\Sigma^\top\\
\Sigma\parallelsum\Sigma^\top & \Sigma^\top
\end{pmatrix}
\begin{pmatrix}\xi\\\eta\end{pmatrix}
\right].
\label{eq:D-invariant-sandwich}
\end{aligned}
\end{equation}
Here $\parallelsum$ denotes the operator geometric mean,
extended to positive semidefinite matrices by continuity.

A q-LAN model $\cS^{(n)}$ with $d$-dimensional quantum
central sequence $\Delta^{(n)}$ is called
\emph{$\cD$-extendible} at $\theta_0$ if there exist an
asymptotically $\cD$-invariant sequence $R^{(n)}$
and a fixed real matrix $F\in\R^{m\times d}$ such that
\[
\Delta^{(n)}=F^\top R^{(n)}.
\]
In this case, the covariance matrix of the central
sequence is $J=F^\top\Sigma F$.
We call the model \emph{asymptotically
$\cD$-invariant} when the central sequence itself
is asymptotically $\cD$-invariant, that is, when
the preceding conditions hold with $m=d$,
$R^{(n)}=\Delta^{(n)}$, and $\Sigma=J$.
\end{definition}

We are now ready to state the asymptotic representation theorem for the models which are q-LAN and $\cD$-extendible. As shown in \cite{Fujiwara2023}, this theorem is crucial in establishing optimality of estimators in quantum i.i.d. models. We will show that this theorem  also serves as one of the key ingredients in deriving upper bounds on local testing power in q-LAN models.

To describe the limiting Gaussian states 
let $\tau=\Sigma F$
and define
\[
\nu_h:=(\operatorname{Re}\tau)h,
\qquad
\varrho_h:=\phi_{\nu_h,\Sigma},
\qquad h\in\R^d.
\]  

\begin{theorem}[Quantum asymptotic representation theorem
\cite{Fujiwara2023}]
\label{thm:representation}
Suppose that $\cS^{(n)}$ is q-LAN and $\cD$-extendible
at $\theta_0$. Let $M^{(n)}$ be POVMs with a common
outcome space $(\R^s,\mathcal B(\R^s))$, and define
\[
\cL_h^{(n)}(B)
:=\Tr\!\left(
\rho_{\theta_0+h/\sqrt n}^{(n)}M^{(n)}(B)
\right),\qquad B\in\mathcal B(\R^s).
\]
Assume that $\cL_h^{(n)}\weak\cL_h$ for every fixed
$h\in\R^d$. Then there exists a POVM $M^{(\infty)}$
on the limiting Gaussian experiment
$\{\varrho_h:h\in\R^d\}$ such that
\[
\cL_h(B)=\varrho_h\!\left(M^{(\infty)}(B)\right),
\qquad h\in\R^d,\quad B\in\mathcal B(\R^s).
\]
\end{theorem}
The weak-convergence assumption can be replaced by
tightness under each fixed local parameter if one
passes to a common subsequence, as follows.

\begin{theorem}[A subsequential version of Quantum Asymptotic Representation Theorem]
\label{thm:common_sub}
Let
\[
   \mathcal S^{(n)}
   =
   \{\,\rho_\theta^{(n)}:\theta\in\Theta\subset\R^d\,\}
\]
be q-LAN and $\mathcal D$-extendible at $\theta_0$,
with $\mathcal D$-extension $R^{(n)}$, limiting covariance
$\Sigma$, and $\tau=\Sigma F$.
Let \(M^{(n)}\) be arbitrary locally tight POVMs on \(\R^s\), i.e. assume that
for every fixed \(h\in\R^d\) the laws
\[
   \mathcal L_h^{(n)}(B)
   =
   \Tr\rho_{\theta_0+h/\sqrt n}^{(n)}M^{(n)}(B)
\]
are tight.
Then there is a subsequence \(n_k\), chosen independently of \(h\), and a POVM
\[
   M^{(\infty)}:\B(\R^s)\to\CCR(\operatorname{Im}\Sigma)_+
\]
such that, for every fixed \(h\in\R^d\),
\[
   \mathcal L_h^{(n_k)}
   =
   \Tr\rho_{\theta_0+h/\sqrt{n_k}}^{(n_k)}M^{(n_k)}
   \weak
   \mathcal L_h,
\]
where the limiting law is represented on the limiting quantum Gaussian shift by
\[
   \cL_h(B)
=
\varrho_h\bigl(M^{(\infty)}(B)\bigr),
\]
In particular, every locally tight sequence of POVMs has a common subsequence
along which its outcome laws converge under all fixed local alternatives.
\end{theorem}

\subsection{Hypothesis testing in quantum statistical models}

We begin by recalling the classical formulation of hypothesis testing before introducing its quantum analogue. The objective is to provide a common decision-theoretic framework that will later be used to formulate optimality criteria for quantum tests.

\subsubsection{Classical hypothesis testing}

Let $\mathcal P=\{P_\theta:\theta\in\Theta\}$
be a parametric statistical model and suppose that
$
X_1,\ldots,X_n
\stackrel{\mathrm{i.i.d.}}{\sim}
P_\theta
$. The basic hypothesis testing problem consists of deciding between the hypotheses
\[
H_0:\theta\in\Theta_0,
\qquad
H_1:\theta\in\Theta_1,
\]
where
\(
\Theta_0\cap\Theta_1=\varnothing.
\)

More generally, many statistical testing problems concern a functional of the underlying distribution rather than the parameter itself. Let
$\vartheta:\mathcal P\rightarrow\mathbb R
$
be a statistical functional and write
$
\vartheta(P_\theta)=\psi(\theta).
$
Testing a property of the functional may then be formulated as
\begin{equation}
H_0:\psi(\theta)\in\mathcal C,
\qquad
H_1:\psi(\theta)\notin\mathcal C,
\label{functional_testing}
\end{equation}
for a prescribed subset
\(
\mathcal C\subseteq\mathbb R.
\)

A (possibly randomized) test is a measurable function
$
\phi_n:\mathcal X^n\rightarrow[0,1],
$
where
\(
\phi_n(x_1,\ldots,x_n)
\)
denotes the probability of rejecting the null hypothesis after observing
\(
(x_1,\ldots,x_n).
\)
The performance of a test is summarized by its power function
\[
\pi_{\phi_n}(\theta)
:=
\mathbb E_\theta
\bigl[
\phi_n(X_1,\ldots,X_n)
\bigr].
\]

A test is said to have level $\alpha$ if
\[
\sup_{\theta:\,\psi(\theta)\in\mathcal C}
\pi_{\phi_n}(\theta)
\le
\alpha.
\]

Let
\(
\mathbb T_\alpha^{(n)}
\)
denote the collection of all level-$\alpha$ tests. A level-$\alpha$ test
\(
\phi_n^\star\in\mathbb T_\alpha^{(n)}
\)
is called \emph{uniformly most powerful (UMP)} if
\[
\pi_{\phi_n^\star}(\theta)
\ge
\pi_{\phi_n}(\theta),
\qquad
\forall\,
\phi_n\in\mathbb T_\alpha^{(n)},
\quad
\forall\,
\theta:\psi(\theta)\notin\mathcal C.
\]

Uniformly most powerful tests need not exist, particularly for two-sided testing problems. In such situations one often restricts attention to the class of unbiased tests. A level-$\alpha$ test is called \emph{unbiased} if
\[
\sup_{\theta:\,\psi(\theta)\in\mathcal C}
\pi_{\phi_n}(\theta)
\le
\alpha,
\qquad
\inf_{\theta:\,\psi(\theta)\notin\mathcal C}
\pi_{\phi_n}(\theta)
\ge
\alpha.
\]

Let
\(
\mathbb T_{\alpha,u}^{(n)}
\)
denote the collection of all level-$\alpha$ unbiased tests. A test
\(
\phi_n^\star
\in
\mathbb T_{\alpha,u}^{(n)}
\)
is called
\emph{uniformly most powerful unbiased (UMPU)}
if
\[
\pi_{\phi_n^\star}(\theta)
\ge
\pi_{\phi_n}(\theta),
\qquad
\forall\,
\phi_n
\in
\mathbb T_{\alpha,u}^{(n)},
\quad
\forall\,
\theta:\psi(\theta)\notin\mathcal C.
\]

Although introduced here for an i.i.d.\ sample, these
definitions apply equally to a general sequence of
classical experiments
$\{P_\theta^{(n)}:\theta\in\Theta\}$:
a test is a measurable $[0,1]$-valued function of
the observation in the $n$th experiment, with power
$\mathbb E_{P_\theta^{(n)}}[\phi_n]$.

\subsubsection{Hypothesis testing with nuisance parameters}

Many hypothesis testing problems involve nuisance parameters that are not of direct inferential interest but nevertheless affect the distribution of the observations. Assume in this subsection that the parameter space has
the product form $\Theta=\Theta_\eta\times\Lambda$,
and write $\theta=(\eta,\lambda)$, where
$\eta\in\Theta_\eta$ is the parameter of interest
and $\lambda\in\Lambda$ is a nuisance parameter.
Let $\Theta_0,\Theta_1\subseteq\Theta_\eta$ be
disjoint sets. The testing problem takes the form
\begin{equation}
H_0:\eta\in\Theta_0,
\qquad
H_1:\eta\in\Theta_1.
\label{nuisance_testing}
\end{equation}

The power function now depends on both the parameter of interest and the nuisance parameter,
\[
\pi_{\phi_n}(\eta,\lambda)
=
\mathbb E_{\eta,\lambda}
\bigl[
\phi_n(X_1,\ldots,X_n)
\bigr],
\]
and the level-$\alpha$ condition becomes
\[
\sup_{\eta\in\Theta_0,\;
\lambda\in\Lambda}
\pi_{\phi_n}(\eta,\lambda)
\le
\alpha.
\]

Uniformly most powerful tests may also fail to exist
in the presence of nuisance parameters. A natural decision-theoretic criterion is therefore to optimize the worst-case power over all possible nuisance parameters. Accordingly, for every
\(\eta\in\Theta_1\),
define the minimum (or worst-case) power function
\[
\underline{\pi}_{\phi_n}(\eta)
:=
\inf_{\lambda\in\Lambda}
\pi_{\phi_n}(\eta,\lambda).
\]

A level-$\alpha$ test
\(
\phi_n^\star
\)
is called
\emph{maximin}
if
\[
\underline{\pi}_{\phi_n^\star}(\eta)
\ge
\underline{\pi}_{\phi_n}(\eta),
\qquad
\forall\,
\phi_n
\in
\mathbb T_\alpha^{(n)},
\quad
\forall\,
\eta\in\Theta_1.
\]

Similarly,
a level-$\alpha$ unbiased test
\(
\phi_n^\star
\)
is called
\emph{unbiased maximin}
if
\[
\underline{\pi}_{\phi_n^\star}(\eta)
\ge
\underline{\pi}_{\phi_n}(\eta),
\qquad
\forall\,
\phi_n
\in
\mathbb T_{\alpha,u}^{(n)},
\quad
\forall\,
\eta\in\Theta_1.
\]

Here maximin is used in the uniform nuisance
sense: the same test maximizes the minimum power
over $\lambda$ for every $\eta\in\Theta_1$.
This is the UMP min-max criterion in the terminology
of \cite[Section~3]{KumagaiHayashi2013}; it is stronger
than optimizing a single minimum over both
$\eta$ and $\lambda$.
We are now ready to discuss the quantum analogues of the notions introduced above.


\subsubsection{Quantum hypothesis testing}

The quantum analogue of the classical testing problem replaces
probability distributions with quantum states and randomized
test functions with binary positive operator-valued measures.

Let
\[
\mathcal S^{(n)}
=
\{\rho_\theta^{(n)}:\theta\in\Theta\}
\]
be a sequence of quantum statistical models on Hilbert spaces
$\mathcal H^{(n)}$. A quantum test is an effect
$0\le T_n\le I$ on $\mathcal H^{(n)}$, or equivalently
the binary POVM $\{T_n,I-T_n\}$, where the outcome
associated with $T_n$ means rejection of the null.
Its power function is
\[
\pi_{T_n}(\theta)
:=
\Tr\!\left(\rho_\theta^{(n)}T_n\right).
\]
The i.i.d.\ experiment is the special case
\[
\rho_\theta^{(n)}=\rho_\theta^{\otimes n},
\qquad
\mathcal H^{(n)}=\mathcal H^{\otimes n}.
\]
Unless stated otherwise, no tensor-product structure is
assumed. This formulation includes the general q-LAN model
sequences considered in Section~\ref{sec:iid_opt}.

The same definitions apply to an experiment represented
by normal states on a von Neumann algebra (see Appendix \ref{sec:vn_alg} for details): a test is
an effect in that algebra, and its rejection probability
is obtained by evaluating the state at the effect.
In particular, for the Gaussian state
$\phi_{\mu,\Sigma}$, tests belong to
$\mathcal M=\operatorname{CCR}(\operatorname{Im}\Sigma)$
and have rejection probability $\phi_{\mu,\Sigma}(T)$.
For the q-LAN limiting state indexed by the local
parameter $h$, we use the notation $\varrho_h(T)$.

The notions of level-\(\alpha\), uniformly most powerful (UMP), and uniformly most powerful unbiased (UMPU) tests are defined exactly as in the classical setting by replacing expectations with quantum expectations. Thus a quantum test \(T_n\) has level \(\alpha\) if
\[
\sup_{\theta:\,\psi(\theta)\in\mathcal C}
\pi_{T_n}(\theta)
\le
\alpha,
\]
and the corresponding classes of level-\(\alpha\) and level-\(\alpha\) unbiased quantum tests are denoted by
$
\mathcal T_\alpha^{(n)}
\text{ and }
\mathcal T_{\alpha,u}^{(n)},
$
respectively. A test
\(
T_n^\star\in\mathcal T_\alpha^{(n)}
\)
is called uniformly most powerful if
\[
\pi_{T_n^\star}(\theta)
\ge
\pi_{T_n}(\theta),
\qquad
\forall\,
T_n\in\mathcal T_\alpha^{(n)},
\quad
\forall\,
\theta:\psi(\theta)\notin\mathcal C,
\]
while
\(
T_n^\star\in\mathcal T_{\alpha,u}^{(n)}
\)
is uniformly most powerful unbiased if the same inequality holds over the smaller class
\(
\mathcal T_{\alpha,u}^{(n)}.
\)

The formulation with nuisance parameters is likewise identical to the classical setting. Writing
$
\theta=(\eta,\lambda),
$
the testing problem is again given by
\eqref{nuisance_testing}, and the power function becomes
\[
\pi_{T_n}(\eta,\lambda)
=
\Tr\!\left(
\rho_{\eta,\lambda}^{(n)}T_n
\right).
\]
For each
\(
\eta\in\Theta_1,
\)
define the minimum power function
\[
\underline{\pi}_{T_n}(\eta)
:=
\inf_{\lambda\in\Lambda}
\pi_{T_n}(\eta,\lambda).
\]
A level-\(\alpha\) quantum test is called \emph{maximin} if its minimum power dominates that of every competing level-\(\alpha\) test, namely,
\[
\underline{\pi}_{T_n^\star}(\eta)
\ge
\underline{\pi}_{T_n}(\eta),
\qquad
\forall\,
T_n\in\mathcal T_\alpha^{(n)},
\quad
\forall\,
\eta\in\Theta_1.
\]
Similarly, a test
$T_n^\star\in\mathcal T_{\alpha,u}^{(n)}$
is called \emph{unbiased maximin} if the same
inequality holds for every
$T_n\in\mathcal T_{\alpha,u}^{(n)}$
and every $\eta\in\Theta_1$.

The decision-theoretic notions introduced above provide the framework for studying optimal hypothesis tests in both classical and quantum statistical models. In the absence of nuisance parameters, one seeks uniformly most powerful or uniformly most powerful unbiased tests whenever they exist. In contrast, the presence of nuisance parameters naturally leads to the maximin criterion, which optimizes the minimum power over all possible values of the nuisance parameter.

Unfortunately, the resulting maximin optimization problem is typically intractable. A powerful approach to overcoming this difficulty is to exploit symmetries of the statistical model. Suppose that the nuisance parameter space admits a group action under which the statistical model is covariant. One may then restrict attention to tests that are invariant under the corresponding group action. Since the class of invariant tests is generally much smaller than the class of all tests, constructing optimal invariant tests is often considerably simpler than solving the original maximin problem.

The Hunt--Stein principle justifies a reduction to
invariant tests under suitable covariance and
averaging assumptions. In the nuisance-maximin
form used here, a uniformly most powerful invariant
test attains the best minimum power over nuisance
parameters. 
Quantum versions of this reduction are discussed in
\cite[Section~3]{KumagaiHayashi2013}. In
Appendix~\ref{sec:quantum-hunt-stein}, we prove the
version needed here for Euclidean nuisance groups
acting on a von Neumann algebra. Its normal-state
formulation covers both quantum and
classical--quantum Gaussian experiments and supplies
the symmetry reduction used in Section~\ref{sec:gauss}.



\section{Hypothesis testing in a quantum Gaussian model with nuisance parameters}
\label{sec:gauss}
In this section we establish upper bounds for power functions in quantum Gaussian models and describe tests achieving them. Fix $0<\alpha<1$ and write
$\Phi(z_\alpha)=1-\alpha$, where $\Phi$ is the standard normal
c.d.f. Throughout this section, $\Sigma=V+iS\succeq0$, with
$V\succ0$ and $S$ real skew-symmetric. 
Throughout this section, the Gaussian shift experiment is
\[
\bigl\{\phi_{h,\Sigma}:h\in\R^d\bigr\},
\qquad
\phi_{h,\Sigma}\sim N(h,\Sigma),
\]
where $\phi_{\mu,\Sigma}$ is the normal Gaussian state
defined in Appendix~\ref{sec:vn_alg}.
All tests are effects in
$\mathcal M=\operatorname{CCR}(\operatorname{Im}\Sigma)$.
When writing $h=(y,z)$, we use the shorthand
$\phi_{y,z}:=\phi_{(y,z),\Sigma}$.
For every effect $T\in\mathcal M$, its rejection
probability under mean $h$ is
$\phi_{h,\Sigma}(T)$.

\subsection{One-sided tests}
Write $h=(y,z)$, with $y\in\mathbb R$ and
$z\in\mathbb R^{d-1}$, and write the Gaussian coordinates as
$X=(Y,Z_1,\ldots,Z_{d-1})^T$, so $Y=X_1$.
The hypotheses are $H_0:y\le0$ versus $H_1:y>0$, with $z$
nuisance. A level-$\alpha$ test satisfies
$\phi_{y,z}(T)\le\alpha$ for every $y\le0$ and every $z$.
Unlike the classical Gaussian problem, pointwise uniformly most
powerful quantum tests need not exist. We therefore compare
$\inf_z\phi_{y,z}(T)$ at each $y>0$.

\begin{theorem}
\label{thm:qGausTest_HS1}
Every level-$\alpha$ effect $T\in\mathcal M$ satisfies
\[
 \inf_{z\in\mathbb R^{d-1}}\phi_{y,z}(T)
 \le 1-\Phi\left(z_\alpha-\frac{y}{\sqrt{V_{11}}}\right),
 \qquad y>0.
\]
Equality is attained by the spectral projection
$T_\alpha=\mathbf1_{(z_\alpha\sqrt{V_{11}},\infty)}(Y)$,
which does not depend on $y$ or $z$.
\end{theorem}

The following formulation involving testing of a linear functional of $h$ will be used in the q-LAN model.
\begin{theorem}
\label{thm:qGausTest_HS}
For $c\in\mathbb R^d\setminus\{0\}$, consider
$H_0:c^Th\le0$ versus $H_1:c^Th>0$ in the model
$\{\phi_{h,\Sigma}:h\in\R^d\}$.
Every level-$\alpha$ effect $T\in\mathcal M$ satisfies
\[
 \inf_{c^Th\ge s}\phi_{h,\Sigma}(T)
 \le 1-\Phi\left(z_\alpha-
                  \frac{s}{\sqrt{c^TVc}}\right),
 \qquad s>0.
\]
Equality is attained by
$T_\alpha=\mathbf1_{(z_\alpha\sqrt{c^TVc},\infty)}(c^TX)$,
simultaneously for all $s>0$.
\end{theorem}
\begin{proof}
Choose an invertible real matrix $L$ with first row $c^T$.
The coordinates $X'=LX$ have mean $Lh$ and covariance
$\Sigma'=L\Sigma L^T$, with $\operatorname{Re}\Sigma'\succ0$
and $V'_{11}=c^TVc$. Apply
Theorem~\ref{thm:qGausTest_HS1} at first mean coordinate $s$:
\[
 \inf_{c^Th\ge s}\phi_{h,\Sigma}(T)
 \le\inf_{c^Th=s}\phi_{h,\Sigma}(T)
 \le1-\Phi\left(z_\alpha-\frac{s}{\sqrt{c^TVc}}\right).
\]
The proposed test has power depending only on $c^Th$, increasing
in that argument, so its infimum is its value at $s$.
\end{proof}

\subsection{Two-sided tests}
For the same model and nonzero $c$, consider
$H_0:c^Th=0$ versus $H_1:c^Th\ne0$. As in the classical case UMP test does not exist in the two-sided problem, so we restrict ourselves to unbiased tests.
A level-$\alpha$ unbiased effect satisfies
$\phi_{h,\Sigma}(T)\le\alpha$ on the null and
$\phi_{h,\Sigma}(T)\ge\alpha$ on the alternative.
The Hunt-Stein reduction is now used within this unbiased class.
\begin{theorem}
\label{thm:qGausTest_two}
Every level-$\alpha$ unbiased effect $T\in\mathcal M$ satisfies
\[
 \inf_{|c^Th|\ge s}\phi_{h,\Sigma}(T)
 \le
 1-\Phi\left(z_{\alpha/2}-\frac{s}{\sqrt{c^TVc}}\right)
   +\Phi\left(-z_{\alpha/2}-\frac{s}{\sqrt{c^TVc}}\right),
 \qquad s>0.
\]
Equality is attained by
$T_\alpha=\mathbf1_{\{|x|>z_{\alpha/2}\sqrt{c^TVc}\}}(c^TX)$,
simultaneously for all $s>0$.
\end{theorem}

\section{Efficient testing in q-LAN models}\label{sec:iid_opt}

In this section, we study efficient testing in q-LAN
models. The upper bounds on power follow from
Theorem~\ref{thm:common_sub} and the Gaussian testing
results of Section~\ref{sec:gauss}. For motivation and
for completeness, the corresponding classical theory
is recalled in Appendix~\ref{sec:class-lan-testing}.

Let
$\mathcal S^{(n)}=\{\rho_\theta^{(n)}:\theta\in\Theta\}$
be a q-LAN model sequence whose central sequence is
asymptotically $\mathcal D$-invariant at
$\theta_0\in\Theta$, where $\Theta\subset\R^d$ is open.
Write
\[
\Delta^{(n)}
\stackrel{\rho_{\theta_0}^{(n)}}{\leadsto}N(0,J),
\qquad G:=\operatorname{Re}J\succ0,
\]
and define, for this general model sequence,
\[
X^{(n)}:=G^{-1}\Delta^{(n)},
\qquad \Sigma:=G^{-1}JG^{-1}.
\]

In the i.i.d.\ case $\rho_\theta^{(n)}=\rho_\theta^{\otimes n}$,
let $L_{\theta_0,j}$ be the single-copy SLDs from
\eqref{eq:general-sld}, and write
\[
L_{\theta_0,j}^{[a]}
:=I^{\otimes(a-1)}\otimes L_{\theta_0,j}
\otimes I^{\otimes(n-a)}.
\]
The SLD central sequence is
\begin{equation}
\Delta_j^{(n)}
:=\frac1{\sqrt n}\sum_{a=1}^nL_{\theta_0,j}^{[a]},
\qquad j=1,\ldots,d,
\label{eq:iid-collective-sld}
\end{equation}
with $J=J_{\theta_0}$ and $G=G_{\theta_0}$ as in
\eqref{eq:general-sld-information}.

Under the assumptions of this section, the sandwiched
Le Cam third lemma (Corollary 3.7 of \cite{Fujiwara2023}) gives, for every fixed $h\in\R^d$,
\begin{equation}
\left(X^{(n)},
\rho_{\theta_0+h/\sqrt n}^{(n)}\right)
\conv{}N(h,\Sigma),
\qquad\Sigma=G^{-1}JG^{-1}.
\label{eq:iid-efficient-gaussian-covariance}
\end{equation}
In the asymptotic representation theorem, this
normalized extension has $F=G$ and
$(\operatorname{Re}\Sigma)F=I$. Its limiting states
are therefore $\varrho_h=\phi_{h,\Sigma}$, so the
Gaussian testing results of Section~\ref{sec:gauss}
apply directly.

For finite $n$, write
\[
\Theta_n:=\{h\in\R^d:\theta_0+h/\sqrt n\in\Theta\}.
\]
All finite-sample local probabilities below are
understood on $\Theta_n$. Every fixed finite set of
local parameters is contained in $\Theta_n$ for all
sufficiently large $n$.
In the finite-submodel power criteria below,
we adopt the convention $\inf\varnothing=+\infty$.

\subsection{One-sided tests}

By the quantum asymptotic representation theorem, limiting
measurement outcome laws in the q-LAN model admit a
representation in the Gaussian shift experiment
$\mathcal S^{(\infty)}
=\{\varrho_h=\phi_{h,\Sigma}:h\in\R^d\}$.
Together with Theorem~\ref{thm:qGausTest_HS}, this
yields the following local asymptotic power bound.

\begin{theorem}
\label{thm:qLAN-test}
Let $\Theta\subset\R^d$ be open and let
$\psi:\Theta\to\R$ be differentiable at
$\theta_0$, with nonzero gradient
$\dot\psi_{\theta_0}$ and $\psi(\theta_0)=0$.
Consider the testing problem
\begin{equation}
\begin{cases}
H_0:\ \psi(\theta)\le0,\\
H_1:\ \psi(\theta)>0.
\end{cases}
\label{eq:qlan_test}
\end{equation}
Every sequence of level-$\alpha$ tests $T^{(n)}$
for the original null, i.e.,
\[
\sup_{\theta:\psi(\theta)\le0}
\Tr\bigl(\rho_\theta^{(n)}T^{(n)}\bigr)
\le\alpha
\qquad\text{for every }n,
\]
satisfies, for every $s>0$,
\begin{equation}
\begin{aligned}
&\inf_{\substack{D\subset\R^d\\0<|D|<\infty}}
\limsup_{n\to\infty}
\inf_{h\in D\cap H_{1,s}^{(n)}}
\Tr\bigl(
\rho_{\theta_0+h/\sqrt n}^{(n)}T^{(n)}
\bigr)\\
&\qquad\le
1-\Phi\left(
z_\alpha-
\frac{s}{
\sqrt{\dot\psi_{\theta_0}
\Sigma\dot\psi_{\theta_0}^T}}
\right),
\end{aligned}
\label{eq:qlan-one-sided-power-bound}
\end{equation}
where $D$ ranges over nonempty finite subsets
of $\R^d$ and
\[
H_{1,s}^{(n)}
:=
\left\{
h\in\R^d:
\sqrt n\,\psi\left(
\theta_0+\frac{h}{\sqrt n}
\right)\ge s
\right\}.
\]

Moreover, the scalar tests
\[
T_\alpha^{(n)}
=
\left\{
\dot\psi_{\theta_0}X^{(n)}
>
z_\alpha
\sqrt{\dot\psi_{\theta_0}
\Sigma\dot\psi_{\theta_0}^T}
\right\}
\]
attain equality in the local asymptotic power
criterion \eqref{eq:qlan-one-sided-power-bound}
for every $s>0$, and satisfy
\[
\lim_{n\to\infty}
\Tr\bigl(\rho_{\theta_0}^{(n)}
T_\alpha^{(n)}\bigr)=\alpha.
\]
\end{theorem}

\begin{remark}
The upper bounds in Theorems~\ref{thm:qLAN-test},
\ref{thm:qLAN-test2}, and
\ref{thm:qLAN-test3_exact_level} apply to tests satisfying
exact finite-sample level, with unbiasedness additionally
required in the two-sided cases. The proposed scalar tests
attain the corresponding Gaussian power bounds in the
stated local asymptotic criteria, but are not asserted to
satisfy these exact finite-sample constraints. All claims
of asymptotic optimality, including those in the
introduction, are understood in this sense.
\end{remark}

The scalar test constructed in this theorem is an
oracle procedure because its defining observable
depends on $\theta_0$. For the smooth finite-dimensional
i.i.d.\ models considered in Section~\ref{sec:wald},
a two-stage adaptive construction attains the same
local Gaussian power without knowledge of $\theta_0$.
It uses single-copy measurements and may be regarded as a quantum analogue
of a Wald procedure.

\subsection{Two-sided tests}


We consider two formulations of the two-sided
testing problem: the original null hypothesis
$\psi(\theta)=0$ and its linearized counterpart
$\dot\psi_{\theta_0}h=0$ in the local parameter.
The first requires an additional continuity
assumption to transfer the level condition to
the tangent null, whereas the second is treated
directly under the q-LAN assumptions. 
\begin{theorem}
\label{thm:qLAN-test2}
Let $\Theta\subset\R^d$ be open, and let
$\psi:\Theta\to\R$ be continuous in a neighborhood
of $\theta_0$ and differentiable at $\theta_0$,
with nonzero gradient $\dot\psi_{\theta_0}$ and
$\psi(\theta_0)=0$. Consider the testing problem
\begin{equation}
\begin{cases}
H_0:\ \psi(\theta)=0,\\
H_1:\ \psi(\theta)\ne0.
\end{cases}
\label{eq:qlan_test2}
\end{equation}
Suppose further that, for every sequence
$h_n\to h\in\R^d$,
\begin{equation}
\left\|
\rho_{\theta_0+h_n/\sqrt n}^{(n)}
-
\rho_{\theta_0+h/\sqrt n}^{(n)}
\right\|_1
\rightarrow0.
\label{trace_conv_local}
\end{equation}

Under the asymptotically $\mathcal D$-invariant
q-LAN assumptions of this section, every sequence
of level-$\alpha$ unbiased tests $T^{(n)}$
for the original hypotheses satisfies, for every
$s>0$,
\begin{equation}
\begin{aligned}
&\inf_{\substack{D\subset\R^d\\0<|D|<\infty}}
\limsup_{n\to\infty}
\inf_{h\in D\cap H_{1,s}^{(n)}}
\Tr\left(
\rho_{\theta_0+h/\sqrt n}^{(n)}
T^{(n)}
\right)\\
&\qquad\le
1-\Phi\left(
-z_{\alpha/2}
-\frac{s}{\sqrt{
\dot\psi_{\theta_0}\Sigma
\dot\psi_{\theta_0}^{T}}},
z_{\alpha/2}
-\frac{s}{\sqrt{
\dot\psi_{\theta_0}\Sigma
\dot\psi_{\theta_0}^{T}}}
\right),
\end{aligned}
\label{eq:qlan-two-sided-power-bound}
\end{equation}
where $D$ ranges over nonempty finite subsets
of $\R^d$ and
\[
H_{1,s}^{(n)}
:=
\left\{
h\in\R^d:
\left|
\sqrt n\,\psi\left(
\theta_0+\frac{h}{\sqrt n}
\right)
\right|\ge s
\right\}.
\]

Moreover, the scalar tests
\[
T_\alpha^{(n)}
=
\left\{
\left|\dot\psi_{\theta_0}X^{(n)}\right|
>
z_{\alpha/2}
\sqrt{
\dot\psi_{\theta_0}\Sigma
\dot\psi_{\theta_0}^{T}}
\right\}
\]
attain equality in the local asymptotic power
criterion \eqref{eq:qlan-two-sided-power-bound}
for every $s>0$, and satisfy
\[
\lim_{n\to\infty}
\Tr\left(
\rho_{\theta_0}^{(n)}T_\alpha^{(n)}
\right)=\alpha.
\]
\end{theorem}


When the hypotheses are formulated directly in the
local tangent parameter, the corresponding result
holds without assumption~\eqref{trace_conv_local},
as stated next.

\begin{theorem}
\label{thm:qLAN-test3_exact_level}
Let $\Theta\subset\R^d$ be open and let
$\psi:\Theta\to\R$ be differentiable at
$\theta_0$, with $\psi(\theta_0)=0$ and
nonzero gradient $\dot\psi_{\theta_0}$.
Under the asymptotically $\mathcal D$-invariant
q-LAN assumptions of this section, put
$c:=\dot\psi_{\theta_0}^T$ and consider
the tangent testing problem
\begin{equation}
\begin{cases}
H_0:\ c^Th=0,\\
H_1:\ c^Th\ne0.
\end{cases}
\label{eq:qlan_test3}
\end{equation}

For each $n$, let
$\Theta_n:=\{h\in\R^d:
\theta_0+h/\sqrt n\in\Theta\}$.
Every sequence of level-$\alpha$ unbiased
tests $T^{(n)}$ for these local hypotheses
satisfies, for each $s>0$,
\begin{equation}
\begin{aligned}
&\inf_{\substack{\varnothing\ne D\subset H_{1,s}\\
                 |D|<\infty}}
\limsup_{n\to\infty}
\inf_{h\in D}
\Tr\left(
\rho_{\theta_0+h/\sqrt n}^{(n)}
T^{(n)}
\right)\\
&\qquad\le
1-\Phi\left(
-z_{\alpha/2}
-\frac{s}{\sqrt{
\dot\psi_{\theta_0}\Sigma
\dot\psi_{\theta_0}^T}},
z_{\alpha/2}
-\frac{s}{\sqrt{
\dot\psi_{\theta_0}\Sigma
\dot\psi_{\theta_0}^T}}
\right),
\end{aligned}
\label{eq:qlan-tangent-two-sided-bound}
\end{equation}
where
\[
H_{1,s}:=
\{h\in\R^d:|c^Th|\ge s\}.
\]
Here level and unbiasedness are understood
in the finite-sample sense for the local
experiment indexed by $\Theta_n$.

Moreover, the scalar tests
\[
T_\alpha^{(n)}
=
\left\{
\left|\dot\psi_{\theta_0}X^{(n)}\right|
>
z_{\alpha/2}
\sqrt{\dot\psi_{\theta_0}
\Sigma\dot\psi_{\theta_0}^T}
\right\}
\]
attain equality in
\eqref{eq:qlan-tangent-two-sided-bound}
for every $s>0$. They have limiting
rejection probability $\alpha$ at every
fixed $h$ satisfying $c^Th=0$, and
limiting rejection probability at least
$\alpha$ at every fixed alternative.
\end{theorem}

\section{Quantum adaptive Wald tests}
\label{sec:wald}

The scalar SLD tests depend on the unknown parameter through the
measurement itself. We implement them using a preliminary estimate on a
vanishing fraction of the copies, followed by a scalar measurement on the
remaining copies. A one-step correction removes the first-order error of
the preliminary estimate. The resulting procedure has the same limiting
local power as the oracle scalar SLD test, while using no knowledge
of the true parameter.
Fujiwara's sequential adaptive method remains an alternative of
practical interest.

\subsection{Model and scalar influence observable}
\label{subsec:wald-model}

Let $\Theta\subset\mathbb R^d$ be open, and let
$\{\rho_\theta:\theta\in\Theta\}$ be a model on $\mathcal H=\mathbb C^r$,
where $r\ge2$ is fixed. Assume that $\theta\mapsto\rho_\theta$ is twice
continuously differentiable in trace norm and that
$\operatorname{rank}\rho_\theta=k$ for every
$\theta\in\Theta$, where $1\le k\le r$ is fixed.
Let \(L_{\theta,j}\), \(j=1,\ldots,d\), denote
any operator-norm continuous choice of SLDs satisfying
\eqref{eq:general-sld}. Such a choice exists under the
constant-rank assumption; see Appendix~\ref{app:wald-identities}. Let
$
J_\theta
$ and $\quad
G_\theta=\operatorname{Re}J_\theta
$
be as in \eqref{eq:general-sld-information}. Assume $G_\theta\succ0$ on $\Theta$. Let $\psi:\Theta\to\mathbb R$ be
$C^2$, and write $c_\theta=\nabla\psi(\theta)=\dot\psi_\theta^T$ as a
column vector. The hypotheses are
\begin{equation}
 \begin{aligned}
 H_0^+ &: \psi(\theta)\le0, & H_1^+ &: \psi(\theta)>0,\\
 H_0^{\pm} &: \psi(\theta)=0, & H_1^{\pm} &: \psi(\theta)\ne0.
 \end{aligned}
 \label{eq:qlan_test1}
\end{equation}
 Define
\begin{equation}
 v_\theta:=G_\theta^{-1}c_\theta,
 \qquad Y_\theta:=\sum_{j=1}^d(v_\theta)_jL_{\theta,j},
 \qquad V_{\mathrm{eff}}(\theta):=c_\theta^TG_\theta^{-1}c_\theta.
 \label{eq:wald-influence}
\end{equation}
The assumptions imply that $Y_\theta$ is operator-norm continuous and
locally bounded. Its spectral measure is
$M_\eta(B):=\mathbf1_B(Y_\eta)$ for Borel sets $B\subset\mathbb R$.

We distinguish the true-state parameter $\vartheta$ from the
measurement-design parameter $\eta$. For fixed $\eta$, the outcome law is
$P_{\vartheta;\eta}(B)=\operatorname{Tr}\{\rho_\vartheta M_\eta(B)\}$,
with mean and variance
\begin{equation}
 \mu(\vartheta,\eta):=\operatorname{Tr}(\rho_\vartheta Y_\eta),
 \qquad
 \nu(\vartheta,\eta):=\operatorname{Tr}(\rho_\vartheta Y_\eta^2)
                         -\mu(\vartheta,\eta)^2.
 \label{eq:wald-two-arguments}
\end{equation}
All derivatives of this outcome law used below are with respect to
$\vartheta$, holding $\eta$ fixed.  The SLD equations give
\begin{equation}
 \mu(\eta,\eta)=0,
 \qquad
 \left.\nabla_\vartheta\mu(\vartheta,\eta)\right|_{\vartheta=\eta}
       =c_\eta,
 \qquad
 \nu(\eta,\eta)=V_{\mathrm{eff}}(\eta).
 \label{eq:wald-influence-identities}
\end{equation}
Consequently, on a sufficiently small neighborhood of any fixed point,
\begin{equation}
 \big|\psi(\eta)+\mu(\vartheta,\eta)-\psi(\vartheta)\big|
 \le C\|\vartheta-\eta\|^2.
 \label{eq:wald-quadratic-cancellation}
\end{equation}
The constant is uniform for both parameters in that neighborhood.
The identities and the bound are proved in
Appendix~\ref{app:adaptive-wald}.

\subsection{Preliminary localization}
\label{subsec:wald-pilot}

Fix a compact set $K\subset\Theta$ with nonempty interior, specified
independently of the unknown parameter. 
Choose integers $m_n<n$ such that
\begin{equation}
 m_n\rightarrow\infty,
 \qquad \frac{m_n}{n}\rightarrow0,
 \qquad \frac{\sqrt n}{m_n}\rightarrow0,
 \qquad N_n:=n-m_n.
 \label{eq:wald-split}
\end{equation}
For example, $m_n=\lfloor n^{2/3}\rfloor$ satisfies these conditions.
Use the first $m_n$ copies to construct a measurable $K$-valued pilot
estimate $\widetilde\theta_n$. The pilot measurement and its processing
rule may be adaptive but must not depend on the unknown true parameter.
The required property is
\begin{equation}
 \sqrt{m_n}(\widetilde\theta_n-\theta)
       =O_{\mathbb P_\theta^{(n)}}(1)
 \quad\text{for each fixed }\theta\in\operatorname{int}K.
 \label{eq:wald-pilot-rate}
\end{equation}
Here $\mathbb P_\theta^{(n)}$ denotes the classical law of all recorded
outcomes of the procedure on $\rho_\theta^{\otimes n}$; a statement about
the pilot uses only its marginal law. In particular, no uniform
local-alternative rate is assumed in \eqref{eq:wald-pilot-rate}.

Fujiwara's sequential-design consistency and asymptotic-normality
results \cite[Theorems~8--9]{adaptive} provide such a pilot when their
positivity, identifiability, smoothness, and design-selection hypotheses
are satisfied; see also \cite{Fujiwara2011Corrigendum}.

The following fixed-measurement
construction supplies a concrete alternative and also fits the
fixed-design specialization of Fujiwara's framework under his smoothness
hypotheses.

\begin{prop}[An informationally complete pilot]
\label{prop:wald-ic-pilot}
In addition to the assumptions of Subsection~\ref{subsec:wald-model},
assume that $\theta\mapsto\rho_\theta$ is injective on $K$.
Put $s=r^2-1$, choose a real basis $H_1,\ldots,H_s$ of the traceless
Hermitian matrices with $\|H_j\|_{\mathrm{op}}\le1$, and fix
$\varepsilon\in(0,1)$. The effects
\begin{equation}
 M_{j,\pm}^{\mathrm{pil}}:=\frac{I\pm\varepsilon H_j}{2s},
 \qquad j=1,\ldots,s,
 \label{eq:wald-ic-povm}
\end{equation}
form a strictly positive informationally complete POVM. Let
$\mathcal X=\{1,\ldots,s\}\times\{+,-\}$ and
$p_\theta(x)=\operatorname{Tr}(\rho_\theta M_x^{\mathrm{pil}})$.
Measure this POVM separately on the pilot copies and, for outcomes
$W_1,\ldots,W_{m_n}$, choose
\begin{equation}
 \widetilde\theta_n\in
 \operatorname*{argmax}_{\eta\in K}
 \sum_{a=1}^{m_n}\log p_\eta(W_a),
 \label{eq:wald-pilot-mle}
\end{equation}
using a fixed choice for each possible outcome-count vector.
For every $\theta\in\operatorname{int}K$, this estimator is consistent
and satisfies
\begin{equation}
 \sqrt{m_n}(\widetilde\theta_n-\theta)
 \weak N_d(0,(I_\theta^{\mathrm{pil}})^{-1}),
 \qquad
 I_\theta^{\mathrm{pil}}
 :=\sum_{x\in\mathcal X}
   \frac{\nabla p_\theta(x)\nabla p_\theta(x)^T}{p_\theta(x)}\succ0.
 \label{eq:wald-pilot-mle-clt}
\end{equation}
In particular, \eqref{eq:wald-pilot-rate} holds.
\end{prop}

The proof of Proposition~\ref{prop:wald-ic-pilot}
is given in Appendix~\ref{app:adaptive-wald}.

\subsection{The two-stage estimator and tests}
\label{subsec:wald-procedure}

Let $\mathcal F_n^{\mathrm{pil}}$ be the sigma-field generated by the pilot
data. Recall that $M_\eta(B)=\mathbf1_B(Y_\eta)$, with
$Y_\eta$ defined in \eqref{eq:wald-influence}. Now
measure $M_{\widetilde\theta_n}$ separately on each of the remaining
$N_n$ copies, and record the numerical outcomes
$Z_{n,1},\ldots,Z_{n,N_n}$. Under the true parameter $\vartheta$, these
outcomes are conditionally independent with common law
$P_{\vartheta;\widetilde\theta_n}$. Define
\begin{equation}
 \widehat\psi_n:=\psi(\widetilde\theta_n)+\overline Z_n,
 \qquad
 \overline Z_n:=\frac1{N_n}\sum_{a=1}^{N_n}Z_{n,a}.
 \label{eq:wald-one-step}
\end{equation}
The correction has conditional mean
$\mu(\vartheta,\widetilde\theta_n)$, so
\eqref{eq:wald-quadratic-cancellation} makes the residual conditional
bias quadratic in the pilot error.
This explains the requirement $\sqrt n/m_n\to0$.
The condition $m_n/n\to0$ ensures that the preliminary stage has no
first-order cost in the final variance.

Now set
\begin{equation}
 \widehat V_n:=\max\{V_{\mathrm{eff}}(\widetilde\theta_n),n^{-1}\},
 \qquad
 W_n:=\frac{\sqrt n\,\widehat\psi_n}{\sqrt{\widehat V_n}}.
 \label{eq:wald-studentization}
\end{equation}
The deterministic lower bound makes the statistic well-defined even if
the pilot falls at a critical point of $\psi$; it is asymptotically
inactive at the regular reference points considered below.
For $\alpha\in(0,1)$ let $z_\alpha=\Phi^{-1}(1-\alpha)$, where $\Phi$
is the standard normal distribution function. The one-sided test rejects
when $W_n>z_\alpha$, and the two-sided test rejects when
$|W_n|>z_{\alpha/2}$. Denote the corresponding quantum effects on
$\mathcal H^{\otimes n}$ by $T_{n,\alpha}^+$ and $T_{n,\alpha}^{\pm}$.
Both effects incorporate the pilot measurement, the pilot-dependent
second-stage measurement, and the final classical decision.
We note that with the pilot in
Proposition~\ref{prop:wald-ic-pilot}, every quantum measurement acts on a
single copy.

\subsection{Local limit and asymptotic power}
\label{subsec:wald-local-theorem}

\begin{theorem}[Local limit of the adaptive scalar estimator]
\label{thm:wald-local-clt}
Under the assumptions of Subsection~\ref{subsec:wald-model}, suppose
$K$, the pilot, and the sample split satisfy
\eqref{eq:wald-split}--\eqref{eq:wald-pilot-rate}. Fix
$\theta_0\in\operatorname{int}K$ with $c_0:=c_{\theta_0}\ne0$, and set
$V_0:=V_{\mathrm{eff}}(\theta_0)>0$.
For every sequence $h_n\to h\in\mathbb R^d$, put
$\theta_n=\theta_0+h_n/\sqrt n$. Under $\mathbb P_{\theta_n}^{(n)}$,
\begin{equation}
 \begin{gathered}
 \sqrt{m_n}(\widetilde\theta_n-\theta_n)
       =O_{\mathbb P_{\theta_n}^{(n)}}(1),\qquad
 \widehat V_n\rightarrow V_0\quad\text{in probability},\\
 \sqrt n\{\widehat\psi_n-\psi(\theta_n)\}
       \weak N(0,V_0).
 \end{gathered}
 \label{eq:wald-centered-local-clt}
\end{equation}
If also $\psi(\theta_0)=0$, then
\begin{equation}
 \sqrt n\,\widehat\psi_n\weak N(c_0^Th,V_0),
 \qquad
 W_n\weak N(c_0^Th/\sqrt{V_0},1).
 \label{eq:wald-shifted-local-clt}
\end{equation}
\end{theorem}

\begin{corollary}[Local power and local asymptotic level]
\label{cor:wald-local-power}
Suppose the assumptions of Theorem~\ref{thm:wald-local-clt} hold and
$\psi(\theta_0)=0$. Define
\begin{equation}
 \begin{aligned}
 \beta_\alpha^+(t;V)
   &:=1-\Phi(z_\alpha-t/\sqrt V),\\
 \beta_\alpha^{\pm}(t;V)
   &:=1-\Phi(z_{\alpha/2}-t/\sqrt V)
       +\Phi(-z_{\alpha/2}-t/\sqrt V).
 \end{aligned}
 \label{eq:wald-power-functions}
\end{equation}
For $\bullet\in\{+,\pm\}$, set
$\pi_n^\bullet(h)=\operatorname{Tr}
 (\rho_{\theta_0+h/\sqrt n}^{\otimes n}T_{n,\alpha}^\bullet)$.
For every compact $H\subset\mathbb R^d$,
\begin{equation}
 \sup_{h\in H}
 \left|\pi_n^\bullet(h)-\beta_\alpha^\bullet(c_0^Th;V_0)\right|
 \rightarrow0.
 \label{eq:wald-uniform-local-power}
\end{equation}
In particular, the tests have local uniform asymptotic level at most
$\alpha$ for their respective actual nonlinear nulls:
\begin{equation}
 \begin{aligned}
 \limsup_n\sup_{\substack{h\in H:\psi(\theta_0+h/\sqrt n)\le0}}
       \pi_n^+(h)&\le\alpha,\\
 \limsup_n\sup_{\substack{h\in H:\psi(\theta_0+h/\sqrt n)=0}}
       \pi_n^{\pm}(h)&\le\alpha.
 \end{aligned}
 \label{eq:wald-local-level}
\end{equation}
An empty supremum in these level statements is taken to be zero.
At $h=0$, both rejection probabilities converge to $\alpha$.
\end{corollary}

The proofs are given in Appendix~\ref{app:adaptive-wald}.
The limiting powers in \eqref{eq:wald-uniform-local-power} coincide with
those of the scalar SLD tests in Section~\ref{sec:iid_opt}, with variance
$V_0=c_0^TG_{\theta_0}^{-1}c_0$. The guarantees established here
are local uniform asymptotic level and convergence to these power
functions. Since $\beta_\alpha^{\pm}(t;V_0)\ge\alpha$, with equality at
$t=0$, the two-sided tests are also locally asymptotically unbiased
for the tangent hypotheses $c_0^Th=0$ versus $c_0^Th\ne0$.


\subsection{Sequential adaptive estimation}
\label{wald-sequential}
Fujiwara's sequential adaptive maximum-likelihood method
\cite{adaptive,Fujiwara2011Corrigendum} provides another route to
adaptive Wald testing. Starting from an initial estimate
$\widehat\theta_0\in K$, at step $k$ one measures
$M_{\widehat\theta_{k-1}}$, records the outcome $x_k$, and updates
\[
\begin{aligned}
 \widehat\theta_k&\in\operatorname*{argmax}_{\vartheta\in K}
                    \ell_k(\vartheta),\\
 \ell_k(\vartheta)&:=\sum_{j=1}^k
 \log\operatorname{Tr}\!\left\{
 \rho_\vartheta M_{\widehat\theta_{j-1}}(\{x_j\})\right\}.
\end{aligned}
\]
The realized measurement settings are held fixed when maximizing the
likelihood. Both the estimate and the measurement are thus updated
throughout the experiment. Our numerical experience suggests that
this sequential scheme may approach its limiting performance at
smaller sample sizes than a two-stage procedure.

Let $A_\theta$ denote the classical Fisher information of
$P_{\vartheta;\theta}$ in $\vartheta$, evaluated at $\vartheta=\theta$.
For the scalar influence observable used here,
\[
 A_\theta v_\theta=c_\theta.
\]
Indeed, the SLD equations imply that the score in direction
$v_\theta$ equals the observed eigenvalue of $Y_\theta$ almost surely; multiplying
by the score vector and taking expectations gives the identity.
Consequently, whenever $A_\theta\succ0$,
\[
 c_\theta^TA_\theta^{-1}c_\theta
 =c_\theta^TG_\theta^{-1}c_\theta
 =V_{\mathrm{eff}}(\theta).
\]
Under Fujiwara's positivity, identifiability, smoothness, and
design-selection hypotheses, the adaptive MLE is strongly consistent
and asymptotically normal with covariance $A_\theta^{-1}$
\cite[Theorems~8--9]{adaptive}.
At a regular null point $\theta_0\in\operatorname{int}K$, with
$\psi(\theta_0)=0$ and $c_{\theta_0}\ne0$, suppose additionally that
\[
 \sqrt n(\widehat\theta_n-\theta_n)
 \ \weak\ N_d(0,A_{\theta_0}^{-1})
 \quad\text{under }\theta_n=\theta_0+h_n/\sqrt n
\]
for every $h_n\to h$. The delta method then gives the same local powers
as in \eqref{eq:wald-power-functions}
for the Wald statistic
\[
 W_n^{\mathrm{seq}}
 :=\frac{\sqrt n\,\psi(\widehat\theta_n)}
 {\sqrt{\max\{V_{\mathrm{eff}}(\widehat\theta_n),n^{-1}\}}}.
\]
This local regularity is an additional requirement beyond
fixed-parameter asymptotic normality.
The constant-rank assumption can be dropped for this sequential approach
whenever a scalar-optimal measurement can be chosen continuously and
the induced classical model satisfies the likelihood and local
regularity conditions above.

If $A_\theta$ is singular, one can use the fixed informationally
complete POVM of Proposition~\ref{prop:wald-ic-pilot} with a fixed
probability $\epsilon\in(0,1)$, and
$M_{\widehat\theta_{k-1}}$ otherwise, recording which measurement was
used and updating the likelihood accordingly. The information at the
true design is then
\[
 A_{\epsilon,\theta}
 =(1-\epsilon)A_\theta+\epsilon I_\theta^{\mathrm{pil}}\succ0.
\]
Under the corresponding likelihood and design regularity conditions,
the scalar asymptotic variance is
$V_\epsilon(\theta):=c_\theta^TA_{\epsilon,\theta}^{-1}c_\theta$.
The information inequalities
$A_\theta,I_\theta^{\mathrm{pil}}\preceq G_\theta$ and
$A_\theta v_\theta=c_\theta$ give
\[
 V_{\mathrm{eff}}(\theta)
 \le V_\epsilon(\theta)
 \le \frac{V_{\mathrm{eff}}(\theta)}{1-\epsilon}.
\]
Thus a sufficiently small fixed $\epsilon$ gives an arbitrarily
small asymptotic variance loss; the corresponding Wald statistic
uses $V_\epsilon$ for studentization. We use the two-stage construction
for its direct proof of the local limit, while sequential adaptation
remains of interest for its performance at moderate sample sizes.


\section{Application to spectral functional testing}\label{sec:spec_test}

In this section, we apply the preceding results
to testing smooth spectral functionals. We compare
the local asymptotic powers of scalar SLD tests
and tests based on weak Schur sampling (WSS),
distinguishing the effects of classical
post-processing and measurement choice.
The comparison identifies conditions under which
the WSS procedures attain the same limiting power
as the scalar SLD test, which achieves the local
maximin bound under the assumptions of Section~\ref{sec:iid_opt}.


Let
\(
\mathcal S
=
\left\{
\rho_\theta:
\theta\in\Theta\subset\R^d
\right\}
\)
be a $C^1$ parametric family of full-rank states
on $\C^r$, where $\Theta$ is open.
Fix $\theta_0\in\Theta$ and write
\begin{equation}
\rho_0:=\rho_{\theta_0}
=
U_0\operatorname{diag}(\boldsymbol\alpha_0)U_0^*,
\qquad
\alpha_{0,1}>\cdots>\alpha_{0,r}>0.
\label{eq:spec-reference}
\end{equation}

Denote
the decreasingly ordered eigenvalue vector of \(\rho_\theta\) by
\(
\boldsymbol{\alpha}_\theta
=
\left(
\alpha_{\theta,1},\ldots,\alpha_{\theta,r}
\right)^T
\in
\Delta_r^\circ,
\)
where
\[
\Delta_r
:=
\left\{
x\in[0,1]^r:
\sum_{i=1}^r x_i=1
\right\},
\qquad
\Delta_r^\circ
:=
\left\{
x\in(0,1)^r:
\sum_{i=1}^r x_i=1
\right\}.
\]
Assume that
\(
\theta\longmapsto\boldsymbol{\alpha}_\theta
\)
is differentiable at a fixed point \(\theta_0\in\Theta\).

Let \(\psi\) be a spectral functional.  Thus there exists a symmetric
function
\[
\varphi:\Delta_r\rightarrow\R
\]
such that
$
\psi(\theta)
=
\varphi(\boldsymbol{\alpha}_\theta).
$
Assume that \(\varphi\) is continuously differentiable in a
neighborhood of
\(
\boldsymbol{\alpha}_0
:=
\boldsymbol{\alpha}_{\theta_0},
\)
and that
\(
\psi(\theta_0)=0,
\dot\psi_{\theta_0}\neq0.
\)

We consider the one-sided hypotheses
\[
H_0:\psi(\theta)\le0,
\qquad
H_1:\psi(\theta)>0.
\]
Next define
\[
A_0
:=
D_\theta\boldsymbol{\alpha}_{\theta_0}
\in\R^{r\times d},
\qquad
b_0
:=
\nabla\varphi(\boldsymbol{\alpha}_0)
\in\R^r,
\]
where gradients are represented as column vectors.  By the chain rule,
\[
c_0
:=
\dot\psi_{\theta_0}^{T}
=
A_0^Tb_0
\in\R^d.
\]
Let
\[
C_0
:=
\operatorname{diag}(\boldsymbol{\alpha}_0)
-
\boldsymbol{\alpha}_0\boldsymbol{\alpha}_0^T.
\]

For $b_0$, choose any local $C^1$ extension of
$\varphi$ from the simplex to a neighborhood
in $\mathbb R^r$. Two such extensions have
gradients differing by a multiple of $\mathbf1$.
Since $\mathbf1^TA_0=0$ and $C_0\mathbf1=0$,
this choice does not affect $c_0$, $C_0b_0$, or the
variance formulas below. We next introduce two tests:
one based on the scalar SLD statistic from
Section~\ref{sec:iid_opt}, and the other based on
weak Schur sampling.

\subsection{Tests based on SLD statistics and weak Schur sampling}
\label{subsec:wss-efficiency}

\textbf{Test based on the SLD statistic}\\

Let $L_{0,j}:=L_{\theta_0,j}$, $j=1,\ldots,d$,
be the SLDs at $\theta_0$, and put
\[
J_0:=J_{\theta_0},
\qquad
G_0:=\re J_0.
\]
Assume that $G_0\succ0$.
Recall from Section \ref{sec:iid_opt}, that
\[
\Delta_{n,j}
:=
\frac1{\sqrt n}
\sum_{a=1}^nL_{0,j}^{[a]},
\qquad
X^{(n)}:=G_0^{-1}\Delta_n,
\]
and
\[
\Sigma_{\theta_0}
:=
G_0^{-1}J_0G_0^{-1}.
\]




Put $v_0:=G_0^{-1}c_0$ and define
\[
Y_0:=\sum_{j=1}^d(v_0)_jL_{0,j}.
\]
Then
\[
c_0^TX^{(n)}
=
\frac1{\sqrt n}\sum_{a=1}^nY_0^{[a]}.
\]
Thus the scalar statistic can be obtained by measuring
$Y_0$ separately on each copy and summing the outcomes.
This is the observable $Y_{\theta_0}$ defined in
Section~\ref{sec:wald}. The scalar SLD variance
for $\psi$ is
\begin{equation}
\label{eq:qlan-spectral-efficient-variance}
V_{\mathrm{eff}}
:=
\dot\psi_{\theta_0}
\Sigma_{\theta_0}
\dot\psi_{\theta_0}^{T}
=
c_0^T\Sigma_{\theta_0}c_0.
\end{equation}
For the canonical SLD covariance
$
\Sigma_{\theta_0}
=
G_0^{-1}J_0G_0^{-1},
$
we have
\[
V_{\mathrm{eff}}
=
c_0^TG_0^{-1}J_0G_0^{-1}c_0
=
c_0^TG_0^{-1}c_0,
\]
since $\operatorname{Im}J_0$ is real
antisymmetric.
Fix a significance level \(\gamma\in(0,1)\), and let
$
z_\gamma:=\Phi^{-1}(1-\gamma).
$ 
The procedures below are local oracle procedures at
$\theta_0$: their defining observables, weights, and
variance calibrations may depend on $\theta_0$.
Consider the local parameter sequence
\[
\theta_n
=
\theta_0+\frac{h}{\sqrt n},
\qquad
h\in\R^d.
\]

Define the q-LAN estimator of the functional by
\begin{equation}
\label{eq:qlan-functional-estimator}
\widehat\psi_{n,1}
:=
\psi(\theta_0)
+
\frac{1}{\sqrt n}\,
\dot\psi_{\theta_0}X^{(n)}
=
\frac{1}{\sqrt n}\,
c_0^TX^{(n)}.
\end{equation}
The corresponding test, calibrated to have local
asymptotic level $\gamma$, is
\begin{equation}
\label{eq:qlan-spectral-test}
T_{n,1}
:=
\left\{
c_0^TX^{(n)}
>
z_\gamma\sqrt{V_{\mathrm{eff}}}
\right\}.
\end{equation}
Here and below, \(\{Y>a\}\) denotes the spectral projection of the
self-adjoint observable \(Y\) corresponding to the interval
\((a,\infty)\).\\
The test $T_{n,1}$ in \eqref{eq:qlan-spectral-test} is
locally asymptotically optimal under the q-LAN and
$\mathcal D$-invariance assumptions of
Theorem~\ref{thm:qLAN-test}. Without $\mathcal D$-invariance,
we do not assert such optimality. Nevertheless, we
retain the notation $V_{\mathrm{eff}}$ in
\eqref{eq:qlan-spectral-efficient-variance} for the
scalar SLD variance, irrespective of whether the
optimality conditions hold. We will, however, show that
\[
V_{\mathrm{eff}}\leq V_{\mathrm{sp}}
\leq V_{\mathrm{WSS}},
\]
where $V_{\mathrm{sp}}$ is the minimum variance
among linear statistics of the Gaussian WSS
fluctuations having the correct local mean in
every direction. Thus $V_{\mathrm{eff}}$ is
also a lower bound for these variances, although
it need not be attained within this class.

\textbf{Test based on weak Schur Sampling}\\

Under the Schur-Weyl decomposition (see Appendix \ref{sec:SW} for details),
\[
(\C^r)^{\otimes n}
\cong
\bigoplus_{\substack{\lambda\vdash n\\\ell(\lambda)\leq r}}
\mathcal Q_\lambda\otimes\mathcal P_\lambda,
\]
let \(\Pi_\lambda\) denote the orthogonal projection onto the
\(\lambda\)-isotypic component.  Weak Schur sampling is the projective
measurement
\[
\left\{
\Pi_\lambda:
\lambda\vdash n,\ \ell(\lambda)\leq r
\right\}.
\]
If the input state is \(\rho_\theta^{\otimes n}\), its outcome
$
\lambda^{(n)}
=
\left(
\lambda_1^{(n)},\ldots,\lambda_r^{(n)}
\right)
$
has the Schur-Weyl distribution (see Appendix \ref{sec:SW} for details) denoted by:
\[
\lambda^{(n)}
\sim
\operatorname{SW}^n(\boldsymbol{\alpha}_\theta).
\]

Define the empirical Young-diagram spectrum and the WSS plug-in
estimator by
\[
\widehat{\boldsymbol{\alpha}}_n
:=
\frac{\lambda^{(n)}}{n},
\qquad
\widehat\psi_{n,2}
:=
\varphi\left(
\widehat{\boldsymbol{\alpha}}_n
\right).
\]
We define the weak-Schur-sampling variance by
\begin{equation}
\label{eq:wss-asymptotic-variance}
V_{\mathrm{WSS}}
:=
b_0^TC_0b_0
=
\sum_{i=1}^r
\alpha_{0,i}b_{0,i}^2
-
\left(
\sum_{i=1}^r\alpha_{0,i}b_{0,i}
\right)^2.
\end{equation}
The corresponding WSS test is
\begin{equation}
\label{eq:wss-spectral-test}
T_{n,2}
:=
\sum_{\substack{\lambda\vdash n,\ \ell(\lambda)\leq r:\\
\sqrt n\,[\varphi(\lambda/n)-\varphi(\boldsymbol{\alpha}_0)]
>
z_\gamma\sqrt{V_{\mathrm{WSS}}}}}
\Pi_\lambda.
\end{equation}
Equivalently, after observing \(\lambda^{(n)}\), the test rejects
\(H_0\) precisely when
\[
\frac{
\sqrt n\,
[\widehat\psi_{n,2}-\psi(\theta_0)]
}{
\sqrt{V_{\mathrm{WSS}}}
}
>
z_\gamma.
\]


\subsection{Weak Schur sampling under local alternatives}
\label{subsec:spec-wss-limit}

\begin{lemma}[Local WSS central limit theorem and delta method]
\label{lem:wss-local-delta-method}
Suppose \eqref{eq:spec-reference} holds. Let
$\boldsymbol\alpha_n\in\Delta_r^\circ$ satisfy
\[
 \sqrt n(\boldsymbol\alpha_n-\boldsymbol\alpha_0)
 \rightarrow a,
 \qquad
 \mathbf1^Ta=0,
\]
and let
$\lambda^{(n)}\sim\operatorname{SW}^n(\boldsymbol\alpha_n)$.
Then
\begin{align}
 \sqrt n(\widehat{\boldsymbol\alpha}_n-\boldsymbol\alpha_n)
 &\weak N_r(0,C_0),
 \label{eq:spec-wss-vector-clt}
 \\
 \sqrt n(\widehat{\boldsymbol\alpha}_n-\boldsymbol\alpha_0)
 &\weak N_r(a,C_0),
 \label{eq:spec-wss-shift-clt}
 \\
 \sqrt n\{
   \varphi(\widehat{\boldsymbol\alpha}_n)
   -
   \varphi(\boldsymbol\alpha_n)
 \}
 &\weak N(0,V_{\mathrm{WSS}}).
 \label{eq:spec-wss-delta-clt}
\end{align}
The vector Gaussian laws are supported on
$\{x\in\mathbb R^r:\mathbf1^Tx=0\}$.
\end{lemma}

\begin{theorem}[Local comparison of the q-LAN and WSS tests]
\label{thm:qlan-wss-comparison}
Assume that the model is q-LAN and \(D\)-extendible at \(\theta_0\). Then, under
$
\rho_{\theta_0+h/\sqrt n}^{\otimes n},
$
the following assertions hold.

\begin{enumerate}[label=\textup{(\roman*)}]
\item The q-LAN estimator satisfies
\begin{equation}
\label{eq:qlan-functional-local-clt}
\sqrt n
\left[
\widehat\psi_{n,1}
-
\psi\left(\theta_0+\frac{h}{\sqrt n}\right)
\right]
\weak
N(0,V_{\mathrm{eff}}).
\end{equation}

\item The WSS estimator satisfies
\begin{equation}
\label{eq:wss-functional-local-clt}
\sqrt n
\left[
\widehat\psi_{n,2}
-
\psi\left(\theta_0+\frac{h}{\sqrt n}\right)
\right]
\weak
N(0,V_{\mathrm{WSS}}).
\end{equation}

\item Let
$
\kappa(h):=\dot\psi_{\theta_0}h=c_0^Th,
$
then the limiting powers of the two tests $T_{n,1}$ and $T_{n,2}$ are
\begin{align}
\label{eq:qlan-local-power-spectral}
\lim_{n\to\infty}
\Tr\left[
\rho_{\theta_0+h/\sqrt n}^{\otimes n}T_{n,1}
\right]
&=
1-\Phi\left(
z_\gamma-\frac{\kappa(h)}{\sqrt{V_{\mathrm{eff}}}}
\right),\\
\label{eq:wss-local-power-spectral}
\lim_{n\to\infty}
\Tr\left[
\rho_{\theta_0+h/\sqrt n}^{\otimes n}T_{n,2}
\right]
&=
1-\Phi\left(
z_\gamma-\frac{\kappa(h)}{\sqrt{V_{\mathrm{WSS}}}}
\right).
\end{align}

\item Along every one-dimensional submodel
\[
t\longmapsto\theta_0+th
\]
for which \(\kappa(h)>0\), the Pitman asymptotic relative efficiency
of \(T_{n,1}\) with respect to \(T_{n,2}\) is
\begin{equation}
\label{eq:qlan-wss-are}
\operatorname{ARE}_{\mathrm P}(T_1,T_2)
=
\frac{V_{\mathrm{WSS}}}{V_{\mathrm{eff}}}.
\end{equation}
\end{enumerate}
\end{theorem}

The variance $V_{\mathrm{WSS}}$ is the limiting variance
of the unrestricted WSS plug-in estimator, whereas
$V_{\mathrm{eff}}$ is that of the scalar SLD estimator.
A measurement and a decision rule based on its outcome
are different statistical objects. We therefore
introduce $V_{\mathrm{sp}}$: the minimum variance of
a linear statistic in the Gaussian WSS limit whose
mean equals $c_0^Th$ for every local parameter $h$.

Put
\[
Z_n:=\sqrt n(\widehat{\boldsymbol{\alpha}}_n
                 -\boldsymbol{\alpha}_0),
\qquad
\mathsf H:=\{z\in\R^r:\mathbf1^Tz=0\}.
\]
Under \(\theta_n=\theta_0+h/\sqrt n\), Lemma \ref{lem:wss-local-delta-method} gives
\[
Z_n\weak Z\sim N_r(A_0h,C_0).
\]
The first-order fluctuation of the unrestricted plug-in is
\[
\sqrt n\{\varphi(\widehat{\boldsymbol{\alpha}}_n)
             -\varphi(\boldsymbol{\alpha}_0)\}
=b_0^TZ_n+o_{\mathbb P}(1).
\]
Its limiting mean is \(c_0^Th\), because \(c_0=A_0^Tb_0\), and
its limiting variance is \(b_0^TC_0b_0=V_{\mathrm{WSS}}\).
We now ask whether a different linear statistic computed from the same
WSS outcome can have the same limiting mean \(c_0^Th\) for every \(h\),
but smaller limiting variance. Since the limit
\(Z\sim N_r(A_0h,C_0)\), for a linear statistic \(w^TZ\) we have,
\[
\mathbb E_h[w^TZ]=w^TA_0h,
\qquad
\operatorname{Var}_h(w^TZ)=w^TC_0w.
\]
Thus \(w^TZ\) has mean \(c_0^Th\) for every \(h\) precisely when
\begin{equation}
A_0^Tw=c_0.
\label{eq:wssp-influence-constraint}
\end{equation}
Among all such \(w\), we therefore seek one minimizing the limiting
variance \(w^TC_0w\). Adding a multiple of \(\mathbf1\) to \(w\)
changes neither the statistic on \(\mathsf H\) nor its variance,
so we choose the representative satisfying
\(\boldsymbol{\alpha}_0^Tw=0\).
For this representative,
\[
w^TC_0w
=
w^T\operatorname{diag}(\boldsymbol{\alpha}_0)w,
\]
giving the optimization problem
\begin{equation}
\min_{w\in\R^r}
\left\{
w^T\operatorname{diag}(\boldsymbol{\alpha}_0)w:
A_0^Tw=c_0,\ \boldsymbol{\alpha}_0^Tw=0
\right\}.
\label{eq:wssp-optimization}
\end{equation}
Put
$
I_{\mathrm{sp},0}
:=
A_0^T\operatorname{diag}(\boldsymbol{\alpha}_0)^{-1}A_0$ and $I_{\mathrm{sp},0}^{+}$ denote its Moore--Penrose inverse.
\begin{remark}
    We record the observation that $v^TI_{\mathrm{sp},0}v = (A_0v)^T\operatorname{diag}(\boldsymbol{\alpha}_0)^{-1}(A_0v)$ for all $v\in\R^d$, and hence we have $\operatorname{ker}(I_{\mathrm{sp},0}) = \operatorname{ker}(A_0)$. Also, as $I_{\mathrm{sp},0}$ is symmetric, \(
\operatorname{Ran}(I_{\mathrm{sp},0})
=
\ker(I_{\mathrm{sp},0})^\perp
=
\ker(A_0)^\perp
=
\operatorname{Ran}(A_0^T).
\) In particular, $c_0=A_0^Tb_0$ belongs to
$\operatorname{Ran}(I_{\mathrm{sp},0})$, so \(
I_{\mathrm{sp},0}I_{\mathrm{sp},0}^{+}c_0=c_0.
\)
\end{remark}

\noindent The unique minimizer of \eqref{eq:wssp-optimization} is
\[
w_0
:=
\operatorname{diag}(\boldsymbol{\alpha}_0)^{-1}
A_0I_{\mathrm{sp},0}^{+}c_0,
\]
and its minimum value is
\[
V_{\mathrm{sp}}
:=
c_0^TI_{\mathrm{sp},0}^{+}c_0.
\]
When \(A_0\) has full column rank, the Moore--Penrose inverse
is the ordinary inverse.

The corresponding linearly post-processed WSS test
rejects when
$
w_0^TZ_n>z_\gamma\sqrt{V_{\mathrm{sp}}}.
$
By Lemma~\ref{lem:wss-local-delta-method},
\[
w_0^TZ_n
\weak
N(c_0^Th,V_{\mathrm{sp}})
\]
under $\theta_n=\theta_0+h/\sqrt n$.
Its limiting power is therefore
\[
1-\Phi\left(
z_\gamma-\frac{c_0^Th}{\sqrt{V_{\mathrm{sp}}}}
\right).
\]
For $c_0^Th>0$, this power is maximal among
the corresponding upper-tail tests based on
linear statistics $w^TZ_n$ satisfying
$A_0^Tw=c_0$.

\begin{corollary}[Comparison of the three asymptotic variances]
\label{cor:wss-efficient}
Let $\Theta\subset\R^d$ be open, and suppose that
$\theta\mapsto\rho_\theta$ is a $C^1$ family of strictly
positive density matrices on $\C^r$ in a neighborhood of
$\theta_0$. Assume that $\rho_{\theta_0}$ has simple
spectrum, that $\varphi$ is symmetric and $C^1$ in a
neighborhood of $\boldsymbol{\alpha}_0$, and that
\[
c_0=A_0^Tb_0\ne0,
\qquad
\re J_{\theta_0}\succ0.
\]
For $V_{\mathrm{eff}}$ as defined in
\eqref{eq:qlan-spectral-efficient-variance}, take the
canonical SLD covariance
\[
\Sigma_{\theta_0}
=
(\re J_{\theta_0})^{-1}
J_{\theta_0}
(\re J_{\theta_0})^{-1}.
\]
Then the three variances satisfy
\begin{equation}
\label{eq:wss-three-variances}
0<V_{\mathrm{eff}}
\leq V_{\mathrm{sp}}
\leq V_{\mathrm{WSS}}.
\end{equation}
Moreover, the following equality criteria hold:
\begin{align}
V_{\mathrm{WSS}}=V_{\mathrm{sp}}
&\quad\Longleftrightarrow\quad
C_0b_0\in\operatorname{Ran}(A_0),
\label{eq:wss-sp-equality}
\\
V_{\mathrm{sp}}=V_{\mathrm{eff}}
&\quad\Longleftrightarrow\quad
[Y_{\theta_0},\rho_{\theta_0}]=0,
\label{eq:sp-eff-equality}
\end{align}
where $Y_{\theta_0}$ is the scalar SLD observable
defined in Section~\ref{sec:wald}.
When $A_0$ has full column rank, the second criterion
is equivalently
\[
\left(\re J_{\theta_0}-I_{\mathrm{sp},0}\right)
I_{\mathrm{sp},0}^{-1}c_0=0.
\]
If the SLD space at $\theta_0$ is $\mathcal D_{\rho_0}$-invariant, then
$
V_{\mathrm{sp}}=V_{\mathrm{eff}},
$
but equality with $V_{\mathrm{WSS}}$ additionally
requires $C_0b_0\in\operatorname{Ran}(A_0)$.
In particular, if the SLD space is $\mathcal D_{\rho_0}$-invariant
and $\operatorname{rank}(A_0)=r-1$, then
\[
V_{\mathrm{WSS}}=V_{\mathrm{sp}}=V_{\mathrm{eff}}.
\]
The Pitman asymptotic
relative efficiency of the corresponding scalar SLD
and WSS plug-in tests is
\[
\operatorname{ARE}_{\mathrm P}(T_1,T_2)
=
\frac{V_{\mathrm{WSS}}}{V_{\mathrm{eff}}}
=
\frac{V_{\mathrm{WSS}}}{V_{\mathrm{sp}}}
\frac{V_{\mathrm{sp}}}{V_{\mathrm{eff}}}
\geq 1,
\]
with equality if and only if both equality
conditions hold.
\end{corollary}

\begin{remark}
    The optimization defining $V_{\mathrm{sp}}$ concerns
linear statistics of the Gaussian WSS limit having
mean $c_0^Th$ for every $h$. It also applies to estimators with
an expansion
\[
\sqrt n\{\widehat\psi_n-\psi(\theta_0)\}
=w^TZ_n+o_{\mathbb P}(1),\qquad A_0^Tw=c_0,
\]
under every fixed local alternative. Such an
estimator need not be linear in the finite-sample
WSS outcome, as the next example illustrates.
\end{remark}

\subsubsection{Two explicit examples}
\label{subsubsec:wssp-examples}

\begin{example}[A removable post-processing loss]
\label{ex:spec-constrained-purity}

For $|\theta|<1/20$, let
\[
\rho_\theta
=
\operatorname{diag}\left(
\frac12+\theta,\frac3{10}-\theta,\frac15
\right),
\qquad
\psi(\theta)
=
\operatorname{Tr}(\rho_\theta^2)-\frac{19}{50}
=
\frac25\theta+2\theta^2.
\]
At $\theta_0=0$, we have
\[
\begin{gathered}
\boldsymbol\alpha_0
=
\begin{pmatrix}1/2\\3/10\\1/5\end{pmatrix},
\qquad
A_0
=
\left.\frac{d\boldsymbol\alpha_\theta}{d\theta}
\right|_{\theta=0}
=
\begin{pmatrix}1\\-1\\0\end{pmatrix},\\
b_0
=
2\boldsymbol\alpha_0
=
\begin{pmatrix}1\\3/5\\2/5\end{pmatrix},
\qquad
c_0=A_0^Tb_0=\frac25.
\end{gathered}
\]
Writing $D_0=\operatorname{diag}(\boldsymbol\alpha_0)$,
the spectral information is
$I_{\mathrm{sp},0}=A_0^TD_0^{-1}A_0
=2+10/3=16/3$.
Since the model is diagonal, its SLD at zero is
$L_{0,1}=\operatorname{diag}(2,-10/3,0)$, and hence
\[
G_0
=\operatorname{Tr}(\rho_0L_{0,1}^2)
=\frac12\cdot 2^2
+\frac3{10}\left(-\frac{10}{3}\right)^2
=\frac{16}{3}
=I_{\mathrm{sp},0}.
\]
Consequently,
$V_{\mathrm{eff}}=c_0^2/G_0=3/100$
and $V_{\mathrm{sp}}=c_0^2/I_{\mathrm{sp},0}=3/100$.
The efficient spectral weight is
\[
w_0
=D_0^{-1}A_0I_{\mathrm{sp},0}^{-1}c_0
=
\begin{pmatrix}2\\-10/3\\0\end{pmatrix}
\frac3{16}\frac25
=
\begin{pmatrix}3/20\\-1/4\\0\end{pmatrix}.
\]

For the unrestricted plug-in,
$\boldsymbol\alpha_0^Tb_0=19/25$, so
\[
\bar b_0
=b_0-(\boldsymbol\alpha_0^Tb_0)\mathbf1
=
\begin{pmatrix}6/25\\-4/25\\-9/25\end{pmatrix}.
\]
Its asymptotic variance is
\[
\begin{aligned}
V_{\mathrm{WSS}}
&=b_0^TC_0b_0
=b_0^TD_0b_0-(\boldsymbol\alpha_0^Tb_0)^2\\
&=\frac12+\frac3{10}\left(\frac35\right)^2
+\frac15\left(\frac25\right)^2
-\left(\frac{19}{25}\right)^2
=\frac{39}{625}.
\end{aligned}
\]
Thus
\begin{equation}
V_{\mathrm{eff}}=V_{\mathrm{sp}}=\frac3{100},
\qquad
V_{\mathrm{WSS}}=\frac{39}{625},
\qquad
V_{\mathrm{WSS}}-V_{\mathrm{sp}}=\frac{81}{2500}.
\label{eq:wssp-commuting-example}
\end{equation}

We now identify the additional fluctuation of the unrestricted
plug-in. Let $h_n\to h$, put $\theta_n=h_n/\sqrt n$, and define
$Z_n:=\sqrt n(\widehat{\boldsymbol\alpha}_n-
\boldsymbol\alpha_0)
=(Z_{n,1},Z_{n,2},Z_{n,3})^T$.
By Lemma~\ref{lem:wss-local-delta-method}, under
$\rho_{\theta_n}^{\otimes n}$,
\[
Z_n\weak Z=(Z_1,Z_2,Z_3)^T
\sim N_3(A_0h,C_0).
\]
In particular, $\mathbf1^TZ=0$ almost surely.
Now observe that
\[
\bar b_0-w_0
=
\begin{pmatrix}9/100\\9/100\\-9/25\end{pmatrix}
=
\frac9{100}\mathbf1
-\frac9{20}\begin{pmatrix}0\\0\\1\end{pmatrix}.
\]
Since $\mathbf1^TZ=0$, we obtain
\begin{equation}
b_0^TZ=\bar b_0^TZ=w_0^TZ-\frac9{20}Z_3.
\label{eq:wssp-example-noise}
\end{equation}
Moreover,
$\operatorname{Cov}(w_0^TZ,Z_3)
=\alpha_{0,3}(w_0)_3
-(\boldsymbol\alpha_0^Tw_0)\alpha_{0,3}=0$,
and
$\operatorname{Var}(Z_3)
=\alpha_{0,3}(1-\alpha_{0,3})=4/25$.
Since $Z$ is Gaussian, $Z_3\sim N(0,4/25)$ is
independent of $w_0^TZ$.
The model fixes the third eigenvalue, and the unrestricted
plug-in adds the independent noise term in
\eqref{eq:wssp-example-noise}, whose variance is
$(9/20)^2(4/25)=81/2500
=V_{\mathrm{WSS}}-V_{\mathrm{sp}}$.

An explicit nonlinear post-processing with the efficient
limit is
\[
\widehat\theta_n^{\mathrm{sp}}
=
\frac45
\frac{\widehat\alpha_{n,1}}
{\widehat\alpha_{n,1}+\widehat\alpha_{n,2}}
-\frac12,
\qquad
\widehat\psi_n^{\mathrm{sp}}
=
\frac25\widehat\theta_n^{\mathrm{sp}}
+2(\widehat\theta_n^{\mathrm{sp}})^2.
\]
The denominator never vanishes for $n\ge1$, since
$\widehat{\boldsymbol\alpha}_n=\lambda^{(n)}/n$
is a decreasingly ordered probability vector and
$\widehat\alpha_{n,1}+\widehat\alpha_{n,2}
\geq\widehat\alpha_{n,1}\geq1/3$.

To compute the local limit, the gradient of
$(x_1,x_2)\mapsto (4/5)x_1/(x_1+x_2)-1/2$
at $(1/2,3/10)$ is $(3/8,-5/8)^T$.
Since $\psi'(0)=2/5$, the chain rule gives
\[
\begin{aligned}
\sqrt n\,\widehat\psi_n^{\mathrm{sp}}
&=\frac25\left(
\frac38Z_{n,1}-\frac58Z_{n,2}
\right)+o_{\mathbb P}(1)\\
&=\frac3{20}Z_{n,1}-\frac14Z_{n,2}
+o_{\mathbb P}(1)
=w_0^TZ_n+o_{\mathbb P}(1).
\end{aligned}
\]
Consequently, under $\rho_{\theta_n}^{\otimes n}$,
\[
\sqrt n\,\widehat\psi_n^{\mathrm{sp}}
\weak N(c_0h,V_{\mathrm{sp}})
=N\left(\frac25h,\frac3{100}\right).
\]
Thus efficient post-processing of the same Schur outcome
achieves $V_{\mathrm{sp}}=V_{\mathrm{eff}}$, whereas
the unrestricted plug-in has variance $V_{\mathrm{WSS}}$.
The corresponding variance ratio is
$V_{\mathrm{WSS}}/V_{\mathrm{sp}}=52/25$.

\end{example}

\begin{example}[A gap beyond optimal linear WSS post-processing]
\label{ex:wssp-rotating}

For $|\theta|<1/20$, let
\[
D_\theta
=
\operatorname{diag}\left(
\frac12+\theta,\frac3{10}-\theta,\frac15
\right),
\qquad
H=
\begin{pmatrix}
0&1&0\\1&0&0\\0&0&0
\end{pmatrix},
\]
and put $U_\theta=\exp(-i\omega\theta H)$ and
$\rho_\theta^{(\omega)}=U_\theta D_\theta U_\theta^*$,
where $\omega$ is a known fixed real constant.
We consider the same functional
$\psi(\theta)
=\operatorname{Tr}((\rho_\theta^{(\omega)})^2)-19/50
=2\theta/5+2\theta^2$.

Since unitary conjugation preserves the spectrum,
\[
\boldsymbol\alpha_\theta
=
\begin{pmatrix}
1/2+\theta\\3/10-\theta\\1/5
\end{pmatrix},
\qquad
A_0
=
\left.\frac{d\boldsymbol\alpha_\theta}{d\theta}
\right|_{\theta=0}
=
\begin{pmatrix}1\\-1\\0\end{pmatrix}.
\]
As in Example~\ref{ex:spec-constrained-purity},
$\boldsymbol\alpha_0=(1/2,3/10,1/5)^T$,
$b_0=2\boldsymbol\alpha_0=(1,3/5,2/5)^T$
and $c_0=A_0^Tb_0=2/5$.
Writing $D_0=\operatorname{diag}(\boldsymbol\alpha_0)$,
we obtain
\[
\begin{aligned}
I_{\mathrm{sp},0}
&=A_0^TD_0^{-1}A_0=\frac{16}{3},\\
w_0
&=D_0^{-1}A_0I_{\mathrm{sp},0}^{-1}c_0
=\begin{pmatrix}3/20\\-1/4\\0\end{pmatrix}.
\end{aligned}
\]
Hence $V_{\mathrm{sp}}=c_0^2/I_{\mathrm{sp},0}
=3/100$ and
$V_{\mathrm{WSS}}=b_0^TC_0b_0=39/625$.
These quantities and the WSS law are unchanged for every $\omega$.

We now compute the SLD information.
Since $U_0=I_3$, we have $\rho_0=D_0$, and
\[
\left.\partial_\theta\rho_\theta^{(\omega)}
\right|_{\theta=0}
=
\operatorname{diag}(1,-1,0)-i\omega[H,D_0]
=
\begin{pmatrix}
1&i\omega/5&0\\
-i\omega/5&-1&0\\
0&0&0
\end{pmatrix}.
\]
The SLD at zero, denoted by $L_0^{(\omega)}$, satisfies
$2\partial_\theta\rho_\theta^{(\omega)}|_{\theta=0}
=\rho_0L_0^{(\omega)}+L_0^{(\omega)}\rho_0$.
Thus
\[
L_0^{(\omega)}
=
\begin{pmatrix}
2&i\omega/2&0\\
-i\omega/2&-10/3&0\\
0&0&0
\end{pmatrix}.
\]
The diagonal contribution to the SLD information is
$\frac12\cdot2^2+
\frac3{10}(-10/3)^2=16/3=I_{\mathrm{sp},0}$.
The off-diagonal contribution, denoted by $K_0$, is
\[
K_0
=
(\alpha_{0,1}+\alpha_{0,2})
|(L_0^{(\omega)})_{12}|^2
=
\left(\frac12+\frac3{10}\right)\frac{\omega^2}{4}
=
\frac{\omega^2}{5}.
\]
Consequently,
$G_0=\operatorname{Tr}(\rho_0(L_0^{(\omega)})^2)
=I_{\mathrm{sp},0}+K_0
=16/3+\omega^2/5$.
It follows that
\begin{equation}
V_{\mathrm{eff}}
=
\frac{c_0^2}{G_0}
=
\frac{(2/5)^2}{16/3+\omega^2/5}
=
\frac{12}{400+15\omega^2}.
\label{eq:wssp-rotating-variance}
\end{equation}
For example, $\omega=4$ gives
\[
V_{\mathrm{eff}}=\frac3{160}
<
V_{\mathrm{sp}}=\frac3{100}
<
V_{\mathrm{WSS}}=\frac{39}{625}.
\]
Thus the additional information is supplied by the
parameter-dependent rotation, which WSS does not record.

Finally, for $\omega\ne0$,
$\mathcal D_{\rho_0}L_0^{(\omega)}$
is a nonzero matrix with zero diagonal, whereas
every nonzero real multiple of $L_0^{(\omega)}$
has nonzero diagonal.
Hence the one-dimensional SLD space is not
$\mathcal D_{\rho_0}$-invariant.

\end{example}

\subsubsection{Connection with quantum-to-classical sufficiency}
\label{subsubsec:wssp-sufficiency}

Our results are directly related to quantum-to-classical sufficiency, which separates selection of a measurement
from optimization of its classical decision rule. In
\cite[Theorem~5.1]{Lahiry2026QCS}, WSS is shown to be sufficient relative to
unitarily invariant POVMs: for every such measurement $N$, a
parameter-independent Markov kernel $K_N$ satisfies
\begin{equation}
 P_\rho^N=K_NP_\rho^{\mathrm{WSS}}
 \quad\text{for every state }\rho.
 \label{eq:wssp-blackwell}
\end{equation}
Together with appropriate symmetrization, this yields exact Bayes and
minimax reductions to WSS for the invariant decision problems treated
in \cite[Theorems~5.4 and~5.6]{Lahiry2026QCS}.

Equation~\eqref{eq:wssp-blackwell} continues to hold when the state model
is restricted to any subfamily. What may fail on a restricted family
is the separate justification for restricting all quantum procedures
to invariant ones. For a family
$\rho_\theta=U_\theta\operatorname{diag}(\boldsymbol\alpha_\theta)U_\theta^*$,
arbitrary conjugation can break the defining relation between spectrum
and eigenbasis. Symmetrization then need not preserve the decision
problem. Example~\ref{ex:wssp-rotating} exhibits precisely this effect.
In contrast, when all conjugates of every admissible spectrum belong
to the parameter space and the eigenbasis is treated as nuisance,
the relevant symmetry reduction is compatible with the model.

Finally, equality of the variances for one functional is much weaker
than experiment-level sufficiency. It concerns one scalar derivative
at one reference point. Quantum-to-classical sufficiency concerns all
decision rules induced by the prescribed measurement class and, in
\eqref{eq:wssp-blackwell}, holds at every sample size. Conversely,
first-order inefficiency of one plug-in does not contradict sufficiency
of the observation from which it is computed.

\section{Discussion}\label{sec:disc}


In this paper, we have developed a local asymptotic
maximin theory for quantum testing using the
asymptotic representation theorem of \cite{Fujiwara2023}.
Under the stated asymptotic $\mathcal D$-invariance
assumptions, the Gaussian bounds constrain the local power of the admissible tests,
and scalar SLD procedures attain these envelopes in
local asymptotic power. We also construct adaptive tests which
have the same local powers under their own smoothness
assumptions; identifying those powers with the general
testing bound additionally requires the hypotheses
of Section~\ref{sec:iid_opt}.
For spectral functionals, we compare the plug-in
rule, the optimal linear first-order WSS statistic,
and the scalar SLD statistic and characterize conditions under which the gaps
$V_{\mathrm{WSS}}-V_{\mathrm{sp}}$ and
$V_{\mathrm{sp}}-V_{\mathrm{eff}}$ vanish.

Several limitations of the present work suggest natural directions for future research. First, our analysis relies crucially on the smoothness of the functional being tested. Extending the theory to nonsmooth functionals will likely require substantially different techniques. Second, although our general testing results apply
to smooth functionals that may depend on both
eigenvalues and eigenvectors, our detailed
measurement comparison concerns spectral
functionals. It would be valuable to develop
corresponding comparisons for eigenvector-dependent
functionals, particularly in models outside the
asymptotic $\mathcal D$-invariance framework, and
to identify practically implementable procedures
with sharp local power guarantees.

A further direction concerns the assumption of a simple spectrum. Our theory of relative efficiency relies critically on this assumption because the limiting power is derived using the Gaussian approximation of \cite{Kahn&Guta}. When the spectrum contains repeated eigenvalues, this approximation no longer applies. Indeed, as shown in \cite{Meliot12}, the relevant limiting distribution is then non-Gaussian, and its analysis requires fundamentally different techniques. Developing an efficiency theory for quantum models with repeated eigenvalues lies beyond the scope of the present paper and remains an important problem for future work.

Finally, although the underlying quantum state may be high-dimensional, the functionals of interest are often low-dimensional, suggesting that suitably targeted inferential procedures may still be possible. However, as the literature on high-dimensional classical statistics demonstrates, conventional LAN theory generally cannot be applied directly in such regimes. Addressing high-dimensional quantum models will therefore require alternative notions of asymptotic optimality, together with inferential and computational methods that go beyond those developed here. We leave these questions for future investigation.

\section*{Acknowledgements}
K.Y. was supported by JSPS KAKENHI Grant Numbers JP23H01090 and JP22K03466 and JST ERATO Grant Number JPMJER2402, Japan. R.N.K. was supported in part by JSPS KAKENHI Grant Numbers JP22H05105 and JP23K25774. R.N.K. also acknowledges financial support from the MEXT Scholarship for Research Students. 

\section*{Disclosure of AI Use}

GPT 5 and 6 were used to assist with reviewing the manuscript, checking mathematical arguments, and copyediting.
The authors have independently verified the content and take full responsibility for its correctness.

\bibliographystyle{plain}
\bibliography{qtest}

\appendix

\section{Weak Schur sampling and the Schur-Weyl distribution} \label{sec:SW}

Let $\cH_d:=\C^d$.  We first define the  representations of the group $U(d)$ and $S(n)$ (the group of unitary matrices and permutations respectively) on $(\mathbb{C}^d)^{\otimes n}$. Let $f_{\mathbf{a}}=f_{{\mathbf{a}_1}}\otimes\ldots\otimes f_{{\mathbf{a}_n}}$, where $f_1,\ldots,f_d$ are basis elements of $\mathbb{C}^d$ and $\mathbf{a}_i\in\{1,\ldots,d\}$ (note that $f_{\mathbf{a}}$'s span $(\mathbb{C}^d)^{\otimes n}$) and consider the following actions:
\begin{align}
    \pi_n(T): f_{\mathbf{a}_1}\otimes\ldots\otimes f_{\mathbf{a}_n}&\rightarrow Tf_{\mathbf{a}_1}\otimes\ldots\otimes Tf_{\mathbf{a}_n}, \quad T\in U(d)\label{T_action}\\
    \tilde{\pi}_n(\sigma): f_{\mathbf{a}_1}\otimes\ldots\otimes f_{\mathbf{a}_n}&\rightarrow f_{\mathbf{a}_{\sigma^{-1}(1)}}\otimes \ldots\otimes f_{\mathbf{a}_{\sigma^{-1}(n)}}, \quad \sigma\in S(n)\label{sigma_action}.
\end{align}
It can be shown that the representation space can be decomposed into direct sum of subspaces each of which is a tensor product of irreducible representations of $U(d)$ and $S(n)$. This decomposition is called Schur-Weyl duality in the literature (for a detailed account see \cite{fulton,goodman}).

Define a partition of $n$ as a vector
$\lambda=(\lambda_1,\ldots,\lambda_d)$
with
\[
    \lambda_1\geq\cdots\geq\lambda_d\geq0,
    \qquad
    \sum_{i=1}^d\lambda_i=n.
\]

A Young diagram is defined by an ordered tuple of integers $\lambda=(\lambda_1,\ldots,\lambda_d)$ with $\lambda_1\geq\ldots\geq\lambda_d\geq0$ and can be represented graphically by a diagram with $d$ lines each containing $\lambda_i$ boxes. For example a typical Young diagram looks like:
\begin{center}
\begin{ytableau}
 $\text{ }$&& \cr
 & \cr
 \cr
\end{ytableau}
\end{center}
with $\lambda=(3,2,1)$.

Thus, a partition of $n$ can be identified with a Young diagram with
$n$ boxes. Recall that
\[
    \Young
    =
    \{\lambda\vdash n:\ell(\lambda)\leq d\}
\]
is the set of partitions of $n$ with at most $d$ nonzero parts,
equivalently Young diagrams with at most $d$ rows.

Schur-Weyl duality gives a unitary decomposition
\[
    \cH_d^{\ot n}
    \cong
    \bigoplus_{\lambda\in\Young}
    \cU_\lambda\ot\cV_\lambda,
\]
where 
$\cU_\lambda$ carries an irreducible representation of $U(d)$ and $\cV_\lambda$ carries an irreducible representation of $S_n$.




\subsection{Weak Schur sampling and Schur-Weyl distribution}
Let $\Pi_\lambda$ denote the orthogonal projection onto the
isotypic component $\mathcal U_\lambda\otimes\mathcal V_\lambda$.

\begin{definition}[Weak Schur sampling]
Weak Schur sampling is the projective measurement
\[
    \WSS_n:=\{\Pi_\lambda:\lambda\in\Young\}.
\]
Its outcome is the Young diagram $\lambda$.
\end{definition}

The induced classical experiment is
\[
    \cF_n
    :=
    \left\{
        P_\rho^{\WSS_n}:
        \rho\in\mathcal S(\mathbb C^d)
    \right\},
\]
where
\[
    P_\rho^{\WSS_n}(\lambda)
    =\Tr\bigl(\rho^{\ot n}\Pi_\lambda\bigr).
\]

It can be shown that (cf. \cite{Lahiry2026QCS}) the weak Schur-sampling distribution is constant on unitary
orbits and therefore depends on \(\rho\) only through its spectrum. In particular the following holds.
\begin{corollary}[Unitary invariance of the weak Schur distribution]
Let \(\rho,\sigma\in\mathcal S(\mathbb C^d)\) be unitarily
conjugate:
\[
\sigma=U\rho U^\dagger
\]
for some \(U\in U(d)\). Then
\[
\operatorname{Tr}
\left(
\sigma^{\otimes n}\Pi_\lambda
\right)
=
\operatorname{Tr}
\left(
\rho^{\otimes n}\Pi_\lambda
\right)
\]
for every \(\lambda\in\mathbb Y_{n,d}\).
\end{corollary}



We may therefore write
\[
    P_\theta^{(n)}(\lambda)
    :=P_\rho^{\WSS_n}(\lambda)
    \quad\text{whenever }\spec(\rho)=\theta.
\]
Next, we describe the Schur-Weyl distribution using standard and semistandard
Young tableaux.

For a cell \(u=(i,j)\) of \(\lambda\), let \(\lambda'_j\) be the
number of rows containing at least \(j\) cells. Its hook length is
\[
h_\lambda(i,j)=\lambda_i-j+\lambda'_j-i+1.
\]
The number of standard Young tableaux of shape \(\lambda\) is
\[
f^\lambda
=
\frac{n!}{\displaystyle\prod_{u\in\lambda}h_\lambda(u)}.
\]

Let \(\operatorname{SSYT}_d(\lambda)\) denote the semistandard Young
tableaux of shape \(\lambda\) with entries in \(\{1,\ldots,d\}\).
Entries increase weakly along rows and strictly down columns. If
\(m_i(T)\) counts the entries equal to \(i\), define the Schur polynomial by
\[
s_\lambda(\theta)
=
\sum_{T\in\operatorname{SSYT}_d(\lambda)}
\prod_{i=1}^{d}\theta_i^{m_i(T)}.
\]

The Schur-Weyl distribution is
\[ P_\theta^{(n)}(\lambda)
=
f^\lambda s_\lambda(\theta),
\qquad \lambda\in\Young.
\]

\section{Local Asymptotic Normality and Asymptotic Representation Theorem in the classical case}

\label{sec:class-lan-testing}

In this section, we recall the classical local asymptotic normality framework and the associated
Gaussian representation of limiting statistics. These results reduce
the analysis of local testing power to a Gaussian shift experiment. Our presentation follows
van der Vaart~\cite[Chapters~7 and~15]{vaart}.

\subsection{Local asymptotic normality}
\label{app:lan-definition}

Consider a sequence of statistical experiments
\[
 \mathcal E_n
   =\bigl(\mathcal X_n,\mathcal A_n,
                \{P_{n,\theta}:\theta\in\Theta\}\bigr),
\]
where \(\Theta\subset\R^p\) is open, and fix \(\theta_0\in\Theta\).
We work with dominated experiments and write \(p_{n,\theta}\) for
densities relative to a common measure \(\mu_n\). Let \(R_n\) be
nonsingular \(p\times p\) matrices satisfying
\(\|R_n^{-1}\|\to0\), and set
\[
 P_{n,h}:=P_{n,\theta_0+R_n^{-1}h},
 \qquad P_{n,0}:=P_{n,\theta_0}.
\]
For every bounded sequence of local parameters, these distributions
are defined for all sufficiently large \(n\).

\begin{definition}[Local asymptotic normality]
\label{def:lan}
The sequence \((\mathcal E_n)\) is locally asymptotically normal at
\(\theta_0\), with normalization \(R_n\) and information matrix \(I\),
if \(I\) is symmetric and nonnegative definite and there are
statistics \(\Delta_n\in\R^p\) such that, under \(P_{n,0}\),
\begin{equation}
 \Delta_n\weak\mathcal{N}_p(0,I),
 \label{eq:lan-central}
\end{equation}
and, for every deterministic sequence \(h_n\to h\in\R^p\),
\begin{equation}
 \log\frac{p_{n,\theta_0+R_n^{-1}h_n}}
               {p_{n,\theta_0}}
   =h^{\mathsf T}\Delta_n
      -\frac12h^{\mathsf T}Ih+o_{P_{n,0}}(1).
 \label{eq:lan-expansion}
\end{equation}
The likelihood ratio is evaluated at the observed data, and
\(\Delta_n\) is called a central sequence.
\end{definition}

This is the dominated-experiment formulation of
\cite[Definition~7.14, p.~104]{vaart}. Tightness of \(\Delta_n\)
allows \(h\) to be replaced by \(h_n\) on the right-hand side of
\eqref{eq:lan-expansion}. In the scalar normalization used below,
\(R_n=r_n\operatorname{Id}_p\), where \(r_n\to\infty\), so that local
alternatives have the form \(\theta_0+h/r_n\).

For independent observations from a model with density \(p_\theta\),
a standard sufficient condition is differentiability in quadratic
mean at \(\theta_0\): there exists
\(\dot\ell_{\theta_0}\in L^2(P_{\theta_0})^p\) such that
\begin{equation}
 \int\left(
   \sqrt{p_{\theta_0+t}}-\sqrt{p_{\theta_0}}
   -\frac12t^{\mathsf T}\dot\ell_{\theta_0}\sqrt{p_{\theta_0}}
       \right)^2\,d\mu=o(\|t\|^2).
 \label{eq:dqm}
\end{equation}
For \(P_{n,\theta}=P_\theta^{\otimes n}\), this gives LAN with
\[
 R_n=\sqrt n\,\operatorname{Id}_p,\qquad
 \Delta_n=\frac1{\sqrt n}\sum_{j=1}^n\dot\ell_{\theta_0}(X_j),
 \qquad
 I=\E_{\theta_0}
       [\dot\ell_{\theta_0}\dot\ell_{\theta_0}^{\mathsf T}],
\]
as in \cite[Theorem~7.2 and Example~7.15]{vaart}. The general
definition also applies to experiments without a product structure.

Throughout the remainder of the appendix, assume that \(I\) is
positive definite. The corresponding Gaussian shift experiment
observes
\begin{equation}
 Z\sim Q_h:=\mathcal{N}_p(h,I^{-1}),\qquad h\in\R^p.
 \label{eq:gaussian-shift}
\end{equation}
Its likelihood ratio satisfies
\begin{equation}
 \log\frac{dQ_h}{dQ_0}(z)
    =h^{\mathsf T}Iz-\frac12h^{\mathsf T}Ih.
 \label{eq:gaussian-likelihood}
\end{equation}
Comparing \eqref{eq:lan-expansion} and
\eqref{eq:gaussian-likelihood} identifies \(I^{-1}\Delta_n\) as
the limiting Gaussian observation. More precisely, LAN implies
mutual contiguity of \(P_{n,h_n}\) and \(P_{n,0}\) for every
\(h_n\to h\), and Le Cam's third lemma gives
\begin{equation}
 \Delta_n\weak\mathcal{N}_p(Ih,I),
 \qquad
 I^{-1}\Delta_n\weak\mathcal{N}_p(h,I^{-1})
 \quad\text{under }P_{n,h_n}.
 \label{eq:central-local-law}
\end{equation}
Contiguity also transfers \(o_{P_{n,0}}(1)\) remainders to these
local alternatives.

\subsection{Asymptotic representation}
\label{app:lan-representation}

The Gaussian experiment describes the possible limiting laws of
statistics under local perturbations. The following formulation
extends the independent-sample statement of
\cite[Theorem~7.10, p.~98]{vaart} to LAN experiments, as explained
in \cite[Section~7.6, p.~104]{vaart}.

\begin{theorem}[Asymptotic representation]
\label{thm:classical-representation}
Suppose \((\mathcal E_n)\) is LAN at \(\theta_0\), with nonsingular
information \(I\). Let \(T_n\) be \(\R^q\)-valued statistics such
that, for every fixed \(h\in\R^p\),
\[
 \Law_{P_{n,h}}(T_n)\weak L_h.
\]
There exists a Markov kernel \(K\) from \(\R^p\) to \(\R^q\),
independent of \(h\), satisfying
\begin{equation}
 L_h(B)=\int_{\R^p}K(z,B)\,Q_h(dz)
 \label{eq:representation}
\end{equation}
for every Borel set \(B\subset\R^q\) and every \(h\).
Equivalently, there is a measurable map \(t\), also independent of
\(h\), such that
\[
 L_h=\Law\bigl(t(Z,U)\bigr),
 \qquad Z\sim Q_h,\quad U\sim\operatorname{Unif}(0,1),
 \quad U\perp Z.
\]
\end{theorem}

\subsection{A local power bound}
\label{app:lan-power}

We now use the scalar normalization \(R_n=r_n\operatorname{Id}_p\).
Let \(\psi:\Theta\to\R\) be differentiable at \(\theta_0\), with
\[
 \psi(\theta_0)=0,\qquad c:=\nabla\psi(\theta_0)\ne0,
\]
where \(c\) is a column vector. Consider the hypotheses
\[
 H_0:\psi(\theta)\le0
 \qquad\text{and}\qquad
 H_1:\psi(\theta)>0.
\]
Because
\[
 r_n\psi(\theta_0+h/r_n)\rightarrow c^{\mathsf T}h,
\]
the local null and alternative are described by the half-spaces
\(c^{\mathsf T}h\le0\) and \(c^{\mathsf T}h>0\).
For a significance level \(0<\eta<1\), write
\[
 \sigma^2=c^{\mathsf T}I^{-1}c,
 \qquad z_\eta=\Phi^{-1}(1-\eta),
\]
where \(\Phi\) is the standard normal distribution function.

\begin{theorem}[Local asymptotic power bound]
\label{thm:power-bound}
Suppose the experiments are LAN at \(\theta_0\), with scalar
normalization \(r_n\to\infty\) and positive definite information
matrix \(I\). Under the preceding assumptions on \(\psi\), every
sequence of level-\(\eta\) tests, satisfying
\begin{equation}
 \sup_{\theta\in\Theta:\,\psi(\theta)\le0}
       \E_{n,\theta}\varphi_n\le\eta
 \quad\text{for each }n,
 \label{eq:test-size}
\end{equation}
obeys, for every \(h\) with \(c^{\mathsf T}h>0\),
\begin{equation}
 \limsup_{n\to\infty}\E_{n,\theta_0+h/r_n}\varphi_n
 \le
 \beta_\eta^\star(h)
 :=
 1-\Phi\left(z_\eta
            -\frac{c^{\mathsf T}h}{\sqrt{c^{\mathsf T}I^{-1}c}}
       \right).
 \label{eq:power-envelope}
\end{equation}
\end{theorem}

This is the power bound of
\cite[Theorem~15.4, p.~219]{vaart}, expressed in the present notation.

\subsection{Attainment of the bound}
\label{app:lan-attainment}

The following result gives the attainment statement corresponding
to \cite[Addendum~15.5, p.~219]{vaart}.

\begin{prop}[Attainment of the local power envelope]
\label{prop:attainment}
Under the assumptions of Theorem~\ref{thm:power-bound},
suppose that
\begin{equation}
S_n=\frac{c^{\mathsf T}I^{-1}\Delta_n}{\sigma}
+o_{P_{n,0}}(1).
\label{eq:efficient-test-statistic}
\end{equation}
For every $h_n\to h$,
\[
S_n\weak
\mathcal N\!\left(\frac{c^{\mathsf T}h}{\sigma},1\right)
\qquad\text{under }P_{n,h_n}.
\]
Consequently, the tests
$\varphi_n^\star=\ind\{S_n>z_\eta\}$ have limiting
power $\beta_\eta^\star(h)$ for every fixed $h$.
For every sequence $h_n\to h$ satisfying
$\psi(\theta_0+h_n/r_n)\le0$ eventually,
\[
\limsup_{n\to\infty}
\E_{n,\theta_0+h_n/r_n}\varphi_n^\star\le\eta.
\]
\end{prop}



An efficient estimator provides a standard implementation.
Suppose
\[
 r_n(\widehat\theta_n-\theta_0)
     =I^{-1}\Delta_n+o_{P_{n,0}}(1),
 \qquad
 \widehat\sigma_n\rightarrow\sigma
       \quad\text{in }P_{n,0}\text{-probability},
\]
where \(\widehat\sigma_n>0\). Differentiability of \(\psi\)
and Slutsky's lemma yield
\[
 W_n:=\frac{r_n\psi(\widehat\theta_n)}{\widehat\sigma_n}
      =\frac{c^{\mathsf T}I^{-1}\Delta_n}{\sigma}
           +o_{P_{n,0}}(1).
\]
The Wald test rejecting when \(W_n>z_{\eta}\) therefore attains
the local power envelope. In regular models,
\(\widehat\sigma_n^2\) can be obtained by consistently estimating
\(c^{\mathsf T}I^{-1}c\); see
\cite[Example~15.6, p.~220]{vaart}.

\begin{remark}[Dependence on the hypotheses]
For the simple Gaussian hypotheses \(h=0\) versus one specified
alternative \(h\), the Neyman-Pearson power is
\[
 1-\Phi\left(z_{\eta}-\sqrt{h^{\mathsf T}Ih}\right).
\]
The half-space envelope \eqref{eq:power-envelope} instead reflects
size control along the entire local null, including nuisance
directions tangent to its boundary.
\end{remark}

\section{Von Neumann algebra formulation of degenerate Gaussian models}
\label{sec:vn_alg}

In Subsection~\ref{subsec:gaussian-states} we described Gaussian
experiments using the usual Hilbert-space representation when the
commutation matrix is nonsingular, and as classical-quantum Gaussian
experiments when the commutation matrix is degenerate. In this appendix we give the corresponding von Neumann algebra
formulation needed for the Gaussian testing results and for later
quasi-characteristic-function calculations.

\paragraph{Weyl representations.}

Let $S\in\R^{d\times d}$ be a real skew-symmetric matrix with $d\geq1$. A
\emph{regular Weyl representation}  associated with $S$ is a strongly
continuous family of unitaries
$\{W_S(\xi):\xi\in\R^d\}
$
on a separable Hilbert space satisfying
\begin{equation}
    W_S(\xi)W_S(\eta)
    =
    e^{-\i\xi^TS\eta}W_S(\xi+\eta),
    \qquad
    W_S(\xi)^*=W_S(-\xi).
    \label{eq:weyl-CCR}
\end{equation}
We denote by
\[
    {\rm CCR}(S)
    :=
    \{W_S(\xi):\xi\in\R^d\}''
\]
the von Neumann algebra generated by the Weyl operators.
This notation depends on the representation; the concrete representation
used throughout the remainder of this appendix is fixed below.

When $S$ is nonsingular, $d=2m$ and the Stone--von Neumann theorem
implies that every regular irreducible representation is unitarily
equivalent to the Schr\"odinger representation on
$L^2(\R)^{\otimes m}$. If
\[
    W_S(\xi)=e^{\i\xi^TX},
    \qquad
    X=(X_1,\ldots,X_{2m})^T,
\]
then, on the common invariant Schwartz core, (cf. \cite{Folland1989}, Section 1.3) \eqref{eq:weyl-CCR} gives
\begin{equation}
    [X_i,X_j]=2\i S_{ij}I.
    \label{eq:CCR-Heisenberg}
\end{equation}
For the standard canonical observables
$(Q_1,P_1,\ldots,Q_m,P_m)$,
\[
    S
    =
    \frac12
    \bigoplus_{j=1}^m
    \Omega,
    \qquad
    \Omega=
    \begin{pmatrix}
        0&1\\
        -1&0
    \end{pmatrix}.
\]
In this representation
${\rm CCR}(S)=\mathcal B(L^2(\R)^{\otimes m})$.

\paragraph{The degenerate case.}

Suppose now that $S$ is possibly singular, with
\[
    \operatorname{rank}(S)=2k,
    \qquad
    r=d-2k.
\]
There exists an invertible real matrix $T$ such that
\begin{equation}
    T^TST
    =
    0_{r\times r}
    \oplus
    \frac12
    \bigoplus_{j=1}^k\Omega.
    \label{eq:degenerate-canonical-form}
\end{equation}

Let $X_j$ be the self-adjoint generator determined by
$W_S(te_j)=e^{\i tX_j}$, and define the transformed observables by
\[
    \widetilde X=T^TX,
    \qquad
    W_{\widetilde X}(\eta):=W_S(T\eta).
\]
Their commutation matrix is $\widetilde S=T^TST$. We realize them as
\[
    \widetilde X=
    (\widetilde Q_1,\ldots,\widetilde Q_r,
      Q_1,P_1,\ldots,Q_k,P_k)^T,
\]
where the first $r$ coordinates commute with all coordinates and
$(Q_j,P_j)$ are canonical pairs. The original coordinates are recovered
as $X=T^{-T}\widetilde X$; in particular, $X_1$ will always refer to the
first original coordinate.

To realize this system on an ordinary Hilbert space, introduce
auxiliary conjugate observables
$\widetilde P_1,\ldots,\widetilde P_r$ and consider
$\widetilde{\mathcal H}
    =
    L^2(\R)^{\otimes(r+k)}.$
The enlarged collection
\[
    (\widetilde Q_1,\widetilde P_1,\ldots,
      \widetilde Q_r,\widetilde P_r,
      Q_1,P_1,\ldots,Q_k,P_k)
\]
has a nondegenerate commutation matrix.
We use the regular irreducible Schr\"odinger representation of this
enlarged system on $\widetilde{\mathcal H}$.

The algebra corresponding to
the original degenerate system is the von Neumann subalgebra generated
by
\[
    \{e^{\i t\widetilde Q_j}:t\in\R,\ 1\le j\le r\}
\]
and by the Weyl operators of the remaining $k$ canonical pairs.
Consequently,
\begin{equation}
    {\rm CCR}(S)
    \simeq
    L^\infty(\R^r)
    \,\bar\otimes\,
    \mathcal B\!\left(L^2(\R)^{\otimes k}\right).
    \label{eq:degenerate-CCR-product}
\end{equation}

Here $L^\infty(\R^r)$ is taken with respect to Lebesgue measure.
From now on, $\mathcal M=\operatorname{CCR}(S)$ denotes this concrete
von Neumann algebra, with the original Weyl operators $W_S(\xi)$.
Zero-dimensional factors are omitted.

This is the von Neumann algebra counterpart of the
classical--quantum representation used in
Subsection~\ref{subsec:gaussian-states}. The auxiliary conjugate
variables $\widetilde P_j$ do not belong to the original
${\rm CCR}(S)$; they are introduced only to obtain a nondegenerate
ambient representation. We refer the reader to \cite{dammeier2023quantum} for further details on quantum-classical hybrid systems.

\paragraph{Gaussian states.}

As in the Gaussian testing model, assume that $V$ is real symmetric
and positive definite and that
\[
    \Sigma:=V+\i S\succeq0.
\]
For every $\mu\in\R^d$, there is a unique normal state
$\phi_{\mu,\Sigma}$ on $\mathcal M$ with the characteristic function
below. We call it the Gaussian state with mean $\mu$ and covariance
$\Sigma$, and write $\phi=\phi_{\mu,\Sigma}\sim N(\mu,\Sigma)$.

\begin{equation}
    \phi(W_S(\xi))
    =
    \exp\left(
        \i\xi^T\mu-\frac12\xi^TV\xi
    \right),
    \qquad
    \xi\in\R^d.
    \label{eq:degenerate-Gaussian-state}
\end{equation}

The condition $V\succ0$ ensures a nondegenerate classical marginal and
normality in the fixed Lebesgue representation. The matrix $\Sigma$
itself is allowed to be singular.
Uniqueness follows because the linear span of the Weyl operators is
ultraweakly dense in $\mathcal M$.

When $S$ is nonsingular this is the usual density-operator Gaussian
state introduced in Subsection~\ref{subsec:gaussian-states}. When
$S$ is singular, the representation
\eqref{eq:degenerate-CCR-product} identifies it with the corresponding
classical--quantum Gaussian experiment.

For clarity, the tensor-product decomposition of $\mathcal M$ does not
require the state to be a product in the chosen coordinates.
Under the change $\widetilde X=T^TX$, the mean and covariance become
\[
 \widetilde\mu=T^T\mu,\qquad
 \widetilde V=T^TVT,\qquad
 \widetilde\Sigma=T^T\Sigma T.
\]
Write $\widetilde X=(C,R)$, with classical coordinates
$C=(\widetilde Q_1,\ldots,\widetilde Q_r)^T$ and quantum coordinates
$R=(Q_1,P_1,\ldots,Q_k,P_k)^T$.
When both parts are present, partition
\[
 \widetilde\mu=
 \begin{pmatrix}\mu_c\\ \mu_q\end{pmatrix},
 \qquad
 \widetilde\Sigma=
 \begin{pmatrix}
    V_c & K\\
    K^T & V_q+\i S_q
 \end{pmatrix},
 \qquad
 S_q=\frac12\bigoplus_{j=1}^k\Omega.
\]
Then $V_c\succ0$ and the Schur complement gives
\[
 \Gamma+\i S_q\succeq0,
 \qquad
 \Gamma:=V_q-K^TV_c^{-1}K.
\]
Let $f$ be the density of $N(\mu_c,V_c)$ and let $\rho_c$ be the usual
quantum Gaussian density operator with mean
$\mu_q+K^TV_c^{-1}(c-\mu_c)$ and covariance $\Gamma+\i S_q$.
In the representation \eqref{eq:degenerate-CCR-product}, the state is
\[
 \phi(A)=\int_{\R^r} f(c)\,
       \operatorname{Tr}\!\bigl(\rho_c A(c)\bigr)\,dc,
 \qquad A\in\mathcal M.
\]
Here $A(c)$ is the operator-valued representative of $A$ in the
tensor-product algebra. This construction proves existence and
normality and allows classical--quantum correlations.
The further linear change
\[
    R'=R-K^TV_c^{-1}C
\]
preserves the commutation relations and gives $(C,R')$ the block-diagonal
covariance $V_c\oplus(\Gamma+\i S_q)$, hence a product Gaussian state.
This is the product form used in Subsection~\ref{subsec:gaussian-states}.
If $r=0$ or $k=0$, the construction reduces to the usual quantum or
classical Gaussian state, respectively.

For later use, the ordered products of Weyl operators satisfy
\begin{equation}
\label{eqn:q-GaussianCharFnc}
\begin{split}
    \phi\left(
        \prod_{t=1}^q W_S(\xi_t)
    \right)
    &=
    \exp\Bigg[
       \sum_{t=1}^q
       \left(
          \i\xi_t^T\mu
          -\frac12\xi_t^T\Sigma\xi_t
       \right)  \\
    &\hspace{2.7cm}
       -
       \sum_{1\le t<u\le q}
          \xi_t^T\Sigma\xi_u
    \Bigg].
\end{split}
\end{equation}
Indeed, \eqref{eq:weyl-CCR} gives
\[
    \prod_{t=1}^q W_S(\xi_t)
    =
    \exp\left(
        -\i\sum_{t<u}\xi_t^TS\xi_u
    \right)
    W_S\left(\sum_{t=1}^q\xi_t\right),
\]
and \eqref{eqn:q-GaussianCharFnc} follows from
\eqref{eq:degenerate-Gaussian-state}. Since
$\xi^TS\xi=0$, one may equivalently replace
$\xi_t^T\Sigma\xi_t$ in the first sum by $\xi_t^TV\xi_t$.

\paragraph{Translation of Gaussian mean}.

The enlarged representation also permits arbitrary translations of
the Gaussian mean. For every $a\in\R^d$ there exists a unitary
$U(a)$ on $\widetilde{\mathcal H}$ such that
\begin{equation}
    U(a)^*W_S(\xi)U(a)
    =
    e^{\i\xi^Ta}W_S(\xi),
    \qquad
    \xi\in\R^d.
    \label{eq:Gaussian-translation-action}
\end{equation}

The unitaries $U(a)$ can be chosen as a strongly continuous projective
representation. Indeed, let $(\widetilde E,\widetilde\omega)$ be the
ambient symplectic space and $\widetilde W$ its Weyl representation,
with convention
\[
 \widetilde W(v)\widetilde W(w)
 =e^{-\i\widetilde\omega(v,w)/2}\widetilde W(v+w).
\]
Identify the original coordinate space $E=\R^d$ with its embedded
subspace, so that $\widetilde W(\xi)=W_S(\xi)$ for $\xi\in E$.
Nondegeneracy permits a linear choice of $v(a)\in\widetilde E$ such that
\[
 \widetilde\omega(v(a),\xi)=\xi^Ta,\qquad \xi\in E.
\]
Taking $U(a)=\widetilde W(v(a))$ gives
\eqref{eq:Gaussian-translation-action}.
Since conjugation maps every Weyl generator to a scalar multiple of
itself, it defines a $*$-automorphism
\[
 \alpha_a(A):=U(a)^*AU(a),\qquad A\in\mathcal M.
\]
These maps are ultraweakly continuous and satisfy
\[
 \alpha_0=\operatorname{id},
 \qquad \alpha_a\circ\alpha_b=\alpha_{a+b}.
\]
The scalar phases in the projective representation disappear under
conjugation. Strong continuity of $U(a)$ also makes this action
point-ultraweakly continuous.
The implementing unitaries may lie outside $\mathcal M$.

Finally, we note that the characteristic function and normality give the covariance
relation
\begin{equation}
    \phi_{\mu,\Sigma}\circ\alpha_a=\phi_{\mu+a,\Sigma}.\label{gauss_translation}
\end{equation}
The following lemma identifies the invariant algebra used in the proof of Theorem~\ref{thm:qGausTest_HS1}.

\begin{lemma}[Invariant algebra under nuisance translations]
\label{lem:Gaussian-invariant-algebra}
For the CCR representation $\mathcal M$ and the displacement unitaries
$U(a)$ constructed above,
\[
\begin{aligned}
 &\bigl\{A\in\mathcal M:
 U(a)^*AU(a)=A\ \text{for all }a\in\R^d\text{ with }a_1=0\bigr\}\\
 &\hspace{3cm}=\{e^{\i tX_1}:t\in\R\}''.
\end{aligned}
\]
In particular, a test $T\in\mathcal M$, $0\leq T\leq I$, is invariant
under these translations if and only if $T=g(X_1)$ for some Borel
function $g:\R\to[0,1]$.
\end{lemma}

\begin{proof}
Denote the algebra on the left by $\mathcal N$.
In the ambient symplectic space introduced above, put
\[
 F:=\{v\in\widetilde E:\widetilde\omega(v,e_1)=0\}.
\]
For $v\in F$, define $a(v)_j:=\widetilde\omega(v,e_j)$.
Then $a(v)_1=0$, and the Weyl relations give
\[
 \widetilde W(v)^*W_S(\xi)\widetilde W(v)
 =e^{\i\xi^Ta(v)}W_S(\xi),\qquad \xi\in E.
\]
Thus conjugation by $\widetilde W(v)$ preserves $\mathcal M$ and agrees
there with conjugation by $U(a(v))$: both maps are ultraweakly continuous
and agree on the Weyl generators. Consequently,
\[
 \mathcal N\subseteq\{\widetilde W(v):v\in F\}'.
\]

Extend $q_1=e_1$ to a symplectic basis
$(q_1,p_1,\ldots,q_n,p_n)$ of $\widetilde E$.
This is possible because $\widetilde\omega$ is nondegenerate:
choose $p_1$ with $\widetilde\omega(e_1,p_1)=1$ and complete the
basis in the symplectic orthogonal complement of their span.
By the Stone--von Neumann theorem, we may work in the corresponding
Schr\"odinger representation on
$L^2(\R)\otimes\mathcal H_{\mathrm{rest}}$, where
$\mathcal H_{\mathrm{rest}}=L^2(\R^{n-1})$ and $X_1=Q_1\otimes I$.
Since
\[
 F=\operatorname{span}\{q_1,q_2,p_2,\ldots,q_n,p_n\},
\]
its Weyl operators generate
\[
 \{\widetilde W(v):v\in F\}''
 =\mathcal D\,\bar\otimes\,\mathcal B(\mathcal H_{\mathrm{rest}}),
\]
where $\mathcal D=L^\infty(\R)$ acts by multiplication on $L^2(\R)$.
The multiplication algebra is maximal abelian (cf. \cite{kadisonringrose1983fundamentals} Example 5.1.6,) so
\[
 \{\widetilde W(v):v\in F\}'
 =\mathcal D\,\bar\otimes\,\mathbb C I_{\mathrm{rest}}
 =\{e^{\i tX_1}:t\in\R\}''.
\]
Here $\mathcal H_{\mathrm{rest}}=\mathbb C$ when $n=1$.
This proves one inclusion.
Conversely, if $a_1=0$, conjugation by $U(a)$ fixes every $e^{\i tX_1}$
and hence their von Neumann algebra by ultraweak continuity.
The assertion for tests follows from the spectral theorem, choosing a
Borel representative with values in $[0,1]$.
\end{proof}

\section{The Quantum Hunt--Stein Theorem}
\label{sec:quantum-hunt-stein}

This section establishes the nuisance-symmetry reduction used in
Section~\ref{sec:gauss}. We treat Euclidean nuisance groups, which
include all translations needed for the Gaussian testing problems.
The formulation uses normal states on a von Neumann algebra and thus
covers classical, quantum, and classical--quantum experiments.
Kumagai and Hayashi~\cite[Theorem~3.2 and Remark~3.2]{KumagaiHayashi2013}
formulate a related noncompact reduction on $\mathcal B(\mathcal H)$.
Since a hybrid Gaussian experiment may instead have a nontrivial
classical center, we prove the required algebraic version directly.

\paragraph{Statistical model and tests.}

Let $\mathcal M$ be a von Neumann algebra with identity $I$ and
separable predual $\mathcal M_*$. We use the duality
$\mathcal M=(\mathcal M_*)^*$; its weak-star topology
$\sigma(\mathcal M,\mathcal M_*)$ is the ultraweak topology.
Thus $A_j\to A$ ultraweakly means that
$\omega(A_j)\to\omega(A)$ for every $\omega\in\mathcal M_*$.
A normal state is a positive ultraweakly continuous linear functional
$\phi$ with $\phi(I)=1$.

Let
\[
 \mathcal S
 =\{\phi_{\theta,\xi}:\theta\in\Theta,\ \xi\in\Xi\}
\]
be a family of normal states on $\mathcal M$, where $\theta$ is the
parameter of interest and $\xi$ is nuisance. Suppose that $\Xi$ is
nonempty and that $\Theta_0,\Theta_1$ are disjoint nonempty sets
with $\Theta=\Theta_0\cup\Theta_1$. Consider
\begin{equation}
 H_0:\theta\in\Theta_0,
 \qquad
 H_1:\theta\in\Theta_1.
 \label{eq:hypothesis}
\end{equation}
When $\mathcal M=\mathcal B(\mathcal H)$, this notation includes
$\phi_{\theta,\xi}(A)=\operatorname{Tr}(\rho_{\theta,\xi}A)$.
For a hybrid experiment, $\phi_{\theta,\xi}$ is evaluated on its
specified observable algebra.

A test is an effect $T\in\mathcal M$, meaning $0\le T\le I$.
Write $\pi_T(\theta,\xi):=\phi_{\theta,\xi}(T)$ for its
rejection probability. At an alternative $\theta\in\Theta_1$,
its type II error probability is
\[
 \beta_T(\phi_{\theta,\xi})
 :=1-\phi_{\theta,\xi}(T).
\]
Fix $\alpha\in(0,1)$. Define
\begin{align*}
 \mathcal T_\alpha
 &:=
 \{T\in\mathcal M:0\le T\le I,\
   \phi_{\theta,\xi}(T)\le\alpha
   \text{ for all }(\theta,\xi)\in\Theta_0\times\Xi\},
 \\
 \mathcal T_{\alpha,u}
 &:=
 \{T\in\mathcal T_\alpha:\
   \phi_{\theta,\xi}(T)\ge\alpha
   \text{ for all }(\theta,\xi)\in\Theta_1\times\Xi\}.
\end{align*}
Thus unbiasedness here means that power is at least the prescribed
level $\alpha$ throughout the alternative.
For $\theta\in\Theta_1$, the minimum nuisance power is
\[
 \underline\pi_T(\theta)
 :=\inf_{\xi\in\Xi}\phi_{\theta,\xi}(T).
\]
Maximizing this quantity is equivalent to minimizing
$\sup_{\xi\in\Xi}\beta_T(\phi_{\theta,\xi})
 =1-\underline\pi_T(\theta)$.

\paragraph{Covariant nuisance action.}

Let $G=(\mathbb R^q,+)$, where $q\ge0$, act on $\Xi$, with action
denoted by $t\circ\xi$. Suppose that
$\{\tau_t:t\in\mathbb R^q\}$ is an action on $\mathcal M$ by normal
$*$-automorphisms:
\[
 \tau_0=\operatorname{id},
 \qquad
 \tau_s\circ\tau_t=\tau_{s+t}.
\]
Assume that the action is point-ultraweakly continuous: for every
$A\in\mathcal M$ and $\omega\in\mathcal M_*$, the function
$t\mapsto\omega(\tau_t(A))$ is continuous.
We assume that the action is covariant, i.e.
\begin{equation}
 \phi_{\theta,\xi}(\tau_t(A))
 =
 \phi_{\theta,t\circ\xi}(A),
 \qquad
 A\in\mathcal M,\quad
 \theta\in\Theta,\quad
 \xi\in\Xi,\quad
 t\in\mathbb R^q.
 \label{eq:covariance}
\end{equation}
Only the nuisance parameter is changed by this action.
For example, if $\tau_t(A)=V_t^*AV_t$ for a strongly continuous
unitary or projective representation on a Hilbert space, then
\eqref{eq:covariance} is the usual relation
$V_t\rho_{\theta,\xi}V_t^*=\rho_{\theta,t\circ\xi}$.

Next we define the class of invariant elements with respect to this action.
\[
 \mathcal M^G
 :=
 \{A\in\mathcal M:\tau_t(A)=A
                   \text{ for every }t\in\mathbb R^q\},\]
This naturally gives us the definition of invariant tests, and invariant unbiased tests.
 \[\mathcal T_{\alpha,G}
 :=\mathcal T_\alpha\cap\mathcal M^G,
 \quad
 \mathcal T_{\alpha,u,G}
 :=\mathcal T_{\alpha,u}\cap\mathcal M^G.
\]

\begin{theorem}[Quantum Hunt--Stein reduction for Euclidean nuisance groups]
\label{thm:hunt-stein}
Under the preceding assumptions, for every effect $T\in\mathcal M$
there exists an invariant effect $\bar T\in\mathcal M^G$ such that
\begin{equation}
 \inf_{\eta\in\Xi}\phi_{\theta,\eta}(T)
 \le
 \phi_{\theta,\xi}(\bar T)
 \le
 \sup_{\eta\in\Xi}\phi_{\theta,\eta}(T),
 \qquad
 \theta\in\Theta,\quad \xi\in\Xi.
 \label{eq:hs-symmetrization-bounds}
\end{equation}
In particular, symmetrization preserves level and level-$\alpha$
unbiasedness, and
\[
 \underline\pi_{\bar T}(\theta)
 \ge
 \underline\pi_T(\theta),
 \qquad \theta\in\Theta_1.
\]
Consequently, for every $\theta\in\Theta_1$,
\begin{align}
 \inf_{T\in\mathcal T_{\alpha,G}}
     \sup_{\xi\in\Xi}\beta_T(\phi_{\theta,\xi})
 &=
 \inf_{T\in\mathcal T_\alpha}
     \sup_{\xi\in\Xi}\beta_T(\phi_{\theta,\xi}),
 \label{eq:hs-minimax-level}
 \\
 \inf_{T\in\mathcal T_{\alpha,u,G}}
     \sup_{\xi\in\Xi}\beta_T(\phi_{\theta,\xi})
 &=
 \inf_{T\in\mathcal T_{\alpha,u}}
     \sup_{\xi\in\Xi}\beta_T(\phi_{\theta,\xi}).
 \label{eq:hs-minimax-unbiased}
\end{align}
\end{theorem}

\begin{proof}
If $q=0$, the group is trivial and we take $\bar T=T$.
Assume henceforth that $q\ge1$.

\emph{Construction of the averages.}
For $R>0$, let $F_R=[-R,R]^q$, and write $|E|$ for the Lebesgue
measure of a measurable set $E\subset\mathbb R^q$. Define
\begin{equation}
 T_R
 :=
 \frac1{|F_R|}
 \int_{F_R}\tau_t(T)\,dt,
 \label{eq:hs-folner-average}
\end{equation}
where the integral is defined in the ultraweak sense. More explicitly, $T_R$ is determined
by
\begin{equation}
 \omega(T_R)
 =
 \frac1{|F_R|}
 \int_{F_R}\omega(\tau_t(T))\,dt,
 \qquad \omega\in\mathcal M_*.
 \label{eq:hs-average-pairing}
\end{equation}
The scalar integrand is continuous. Since $*$-automorphisms are
isometric, the right-hand side is a linear functional of $\omega$
with absolute value at most $\|\omega\|\,\|T\|$.
The identity $\mathcal M=(\mathcal M_*)^*$ therefore gives a unique
$T_R\in\mathcal M$. Moreover,
$0\le\tau_t(T)\le I$ for every $t$. Evaluating
\eqref{eq:hs-average-pairing} on positive normal functionals gives
$0\le\omega(T_R)\le\omega(I)$, and hence $0\le T_R\le I$.

Using the covariance property of the action \eqref{eq:covariance} yields
\begin{equation}
 \phi_{\theta,\xi}(T_R)
 =
 \frac1{|F_R|}
 \int_{F_R}\phi_{\theta,t\circ\xi}(T)\,dt.
 \label{eq:hs-averaged-power}
\end{equation}
Consequently,
\begin{equation}
 \inf_{\eta\in\Xi}\phi_{\theta,\eta}(T)
 \le
 \phi_{\theta,\xi}(T_R)
 \le
 \sup_{\eta\in\Xi}\phi_{\theta,\eta}(T)
 \label{eq:hs-finite-average-bounds}
\end{equation}
for every $R$, $\theta$, and $\xi$.

\emph{Approximate invariance.}
Fix $s\in\mathbb R^q$. Normality of $\tau_s$ means that
$\omega\circ\tau_s\in\mathcal M_*$ whenever $\omega\in\mathcal M_*$.
Thus \eqref{eq:hs-average-pairing}, and a change of
variables imply
\[
 \tau_s(T_R)
 =
 \frac1{|F_R|}
 \int_{F_R}\tau_{s+t}(T)\,dt
 =
 \frac1{|F_R|}
 \int_{F_R+s}\tau_t(T)\,dt.
\]
Subtracting \eqref{eq:hs-folner-average}, the integrals over
$F_R\cap(F_R+s)$ cancel. Therefore
\[
 |\omega(\tau_s(T_R)-T_R)|
 \le
 \|\omega\|\,\|T\|
 \frac{|(F_R+s)\triangle F_R|}{|F_R|},
 \qquad \omega\in\mathcal M_*,
\]
where $\triangle$ denotes symmetric difference. Taking the supremum
over $\|\omega\|\le1$ gives
\begin{equation}
 \|\tau_s(T_R)-T_R\|
 \le
 \|T\|
 \frac{|(F_R+s)\triangle F_R|}{|F_R|}.
 \label{eq:hs-operator-folner-bound}
\end{equation}
For $R>\frac12\max_j|s_j|$, the intersection of the two cubes has
volume $\prod_{j=1}^q(2R-|s_j|)$, so
\begin{equation}
 \frac{|(F_R+s)\triangle F_R|}{|F_R|}
 =
 2\left\{
 1-\prod_{j=1}^q
       \left(1-\frac{|s_j|}{2R}\right)
 \right\}
 \le
 \frac1R\sum_{j=1}^q|s_j|
 \rightarrow0.
 \label{eq:hs-cube-folner-property}
\end{equation}
This is the F\o lner property of the cubes. In particular,
$\|\tau_s(T_R)-T_R\|\to0$ for every fixed $s$.

\emph{Compactness and an invariant limit.}
The effect interval $[0,I]_{\mathcal M}$ is ultraweakly compact.
Indeed, the closed unit ball is weak-star compact by the
Banach--Alaoglu theorem, and $[0,I]_{\mathcal M}$ is weak-star closed:
its defining inequalities can be tested against positive normal
functionals. Since $\mathcal M_*$ is separable, this compact set is
metrizable and therefore sequentially compact; see also
\cite[Chapter~III]{TakesakiI} for the predual framework.
Taking $R=1,2,\ldots$, choose a subsequence $R_j\to\infty$ and an
effect $\bar T\in\mathcal M$ such that
\begin{equation}
 T_{R_j}\rightarrow\bar T
 \quad\text{in }\sigma(\mathcal M,\mathcal M_*).
 \label{eq:hs-ultraweak-limit}
\end{equation}

For fixed $s\in\mathbb R^q$ and $\omega\in\mathcal M_*$,
normality of $\tau_s$, \eqref{eq:hs-ultraweak-limit}, and
\eqref{eq:hs-operator-folner-bound} give
\[
 \omega(\tau_s(\bar T)-\bar T)
 =
 \lim_{j\to\infty}
 \omega(\tau_s(T_{R_j})-T_{R_j})
 =0.
\]
Because the predual separates points of $\mathcal M$,
$\tau_s(\bar T)=\bar T$. The element $s$ was arbitrary, so
$\bar T\in\mathcal M^G$.

\emph{Level, unbiasedness, and minimax values.}
Every $\phi_{\theta,\xi}$ is normal. Passing to the limit in
\eqref{eq:hs-finite-average-bounds} proves
\eqref{eq:hs-symmetrization-bounds} simultaneously for all
$\theta\in\Theta$ and $\xi\in\Xi$.
If $T\in\mathcal T_\alpha$, its upper bound is at most $\alpha$
on $\Theta_0$, so $\bar T\in\mathcal T_{\alpha,G}$.
If $T\in\mathcal T_{\alpha,u}$, its lower bound is at least
$\alpha$ on $\Theta_1$, giving
$\bar T\in\mathcal T_{\alpha,u,G}$ as well.
For every $\theta\in\Theta_1$,
\[
 \sup_{\xi\in\Xi}\beta_{\bar T}(\phi_{\theta,\xi})
 \le
 \sup_{\xi\in\Xi}\beta_T(\phi_{\theta,\xi}).
\]
Taking the infimum over either admissible class proves one direction
of \eqref{eq:hs-minimax-level}--\eqref{eq:hs-minimax-unbiased}; the
reverse direction follows because the invariant class is a subclass
of the original class.
\end{proof}

The proof requires neither completeness of the statistical family
nor a Gaussian assumption. It selects an invariant cluster point for
each test; it does not assert convergence of the full family of
averages or construct a normal conditional expectation onto
$\mathcal M^G$.

\begin{corollary}[Optimal invariant tests are nuisance-maximin]
\label{cor:hunt-stein}
Assume the hypotheses of Theorem~\ref{thm:hunt-stein}. Let
$\mathcal C$ be either $\mathcal T_\alpha$ or
$\mathcal T_{\alpha,u}$, and put
$\mathcal C_G=\mathcal C\cap\mathcal M^G$.
Suppose that $T_0\in\mathcal C_G$ is uniformly most powerful within
$\mathcal C_G$, in the sense that
\[
 \phi_{\theta,\xi}(T_0)
 \ge
 \phi_{\theta,\xi}(A),
 \qquad
 A\in\mathcal C_G,\quad
 \theta\in\Theta_1,\quad
 \xi\in\Xi.
\]
Then
\begin{equation}
 \underline\pi_{T_0}(\theta)
 =
 \sup_{T\in\mathcal C}\underline\pi_T(\theta),
 \qquad \theta\in\Theta_1.
 \label{eq:hs-invariant-optimality}
\end{equation}
Thus $T_0$ is uniformly maximin over the nuisance parameter in the
same admissible class, including the unbiased class when it is used.
\end{corollary}

\begin{proof}
For $T\in\mathcal C$, let $\bar T\in\mathcal C_G$ be supplied by
Theorem~\ref{thm:hunt-stein}. The assumed invariant optimality gives
\[
 \underline\pi_{T_0}(\theta)
 \ge
 \underline\pi_{\bar T}(\theta)
 \ge
 \underline\pi_T(\theta),
 \qquad \theta\in\Theta_1.
\]
Taking the supremum over $T\in\mathcal C$, and using
$T_0\in\mathcal C$, proves \eqref{eq:hs-invariant-optimality}.
\end{proof}

\paragraph{Application to the Gaussian nuisance translations.}

In Section~\ref{sec:gauss}, take
$\mathcal M=\operatorname{CCR}(S)$ and the normal Gaussian states
$\phi_{y,z}\sim N((y,z),\Sigma)$, where
$\Sigma=V+iS\succeq0$, $V\succ0$, $y\in\mathbb R$, and
$z\in\mathbb R^{d-1}$. No nonsingularity assumption on $S$ is
required. If $\operatorname{rank}(S)=2k$ and $r=d-2k$, the concrete
representation described above gives
\[
 \mathcal M
 \simeq
 L^\infty(\mathbb R^r)\,\bar\otimes\,
 \mathcal B\bigl(L^2(\mathbb R^k)\bigr).
\]
Here $L^2(\mathbb R^0)=\mathbb C$ when $k=0$.
This algebra acts on the separable Hilbert space
$\widetilde{\mathcal H}=L^2(\mathbb R^{r+k})$ and has separable
predual: the predual is a quotient of
$\mathcal T_1(\widetilde{\mathcal H})$ under restriction of normal
functionals, and the latter trace-class space is separable.

For $t\in\mathbb R^{d-1}$, set
$\tau_t:=\alpha_{(0,t)}$, with $\alpha_a$ as in
\eqref{eq:Gaussian-translation-action}. Thus
\[
 \tau_t(W_S(\xi))
 =
 e^{i\xi^T(0,t)}W_S(\xi).
\]
These are normal automorphisms. Their implementations in the enlarged
Schr\"odinger representation may be chosen as a strongly continuous
family of Weyl displacement unitaries $U_t$. In particular, for a
trace-class operator $K$ and $A\in\mathcal M$,
\[
 \left|
 \operatorname{Tr}\{K(U_t^*AU_t-U_s^*AU_s)\}
 \right|
 \le
 \|U_tKU_t^*-U_sKU_s^*\|_1\,\|A\|
 \rightarrow0
 \qquad(t\to s).
\]
The trace-norm convergence follows first for finite-rank $K$ from
strong continuity and then for arbitrary $K$ by trace-norm
approximation. Restricting these normal functionals to $\mathcal M$
verifies point-ultraweak continuity. The unitaries $U_t$ may lie
outside $\mathcal M$; only the induced automorphisms of $\mathcal M$
are used.

The Gaussian characteristic function gives
\[
 \phi_{y,z}(\tau_t(W_S(\xi)))
 =
 \exp\left\{
 i\xi^T(y,z+t)-\tfrac12\xi^TV\xi
 \right\}
 =
 \phi_{y,z+t}(W_S(\xi)).
\]
Both sides define normal functionals of the observable. Since the
linear span of the Weyl operators is ultraweakly dense in
$\mathcal M$, this proves
\[
 \phi_{y,z}(\tau_t(A))
 =
 \phi_{y,z+t}(A),
 \qquad A\in\mathcal M,
\]
which is \eqref{eq:covariance}. Hence all hypotheses of
Theorem~\ref{thm:hunt-stein} are satisfied with $q=d-1$ and
$\Xi=\mathbb R^{d-1}$. For $d=1$, the nuisance group is trivial.

For one-sided testing use $\Theta_0=(-\infty,0]$ and
$\Theta_1=(0,\infty)$; for two-sided testing use
$\Theta_0=\{0\}$ and $\Theta_1=\mathbb R\setminus\{0\}$,
together with the unbiased class. The nuisance translations are
transitive on $\Xi$, so invariant power is independent of $z$.
For a general nonzero contrast $c^Th$, use the group
$\ker(c^T)$, identified linearly with $\mathbb R^{d-1}$.

This establishes the symmetrization step required in the Gaussian
testing proofs. The form of the invariant effects is supplied
separately by Lemma~\ref{lem:Gaussian-invariant-algebra}; no
fixed-point identification is inferred from the averages themselves.
The conclusion concerns minimum power over nuisance parameters,
not pointwise domination of every competing test at each nuisance
value.

\section{Proofs of the main theorems}
\begin{proof}[Proof of Theorem~\ref{thm:common_sub}] We follow the proof of the asymptotic quantum representation theorem in \cite[Theorem~2.4]{Fujiwara2023}. The hypothesis that \(\cL_h^{(n)}\weak\cL_h\) is not used in Steps~1 and~2 of that proof, nor in Step~3 before the final paragraph beginning ``Fix \(h\in\R^d\) arbitrarily.'' Consequently, the same argument up to that point gives the following objects: \begin{enumerate}[label=\textup{(\roman*)}] \item a subsequence \(\{n_k\}_{k\geq1}\), chosen independently of \(h\); \item a coordinatewise nondecreasing family \(\{\widetilde M_\alpha:\alpha\in\Q^s\}\subset\CCR(\Imm\Sigma)_+\), with \(0\leq\widetilde M_\alpha\leq I\); \item a POVM \(\overline M\) on \(\overline{\R}^{\,s}\), determined by \begin{equation} \label{eq:bar-M-from-rational} \overline M_t := \overline M((-\infty,t]) = \inf_{\substack{\alpha\in\Q^s\\\alpha>t}} \widetilde M_\alpha, \qquad t\in\R^s, \end{equation} where the infimum is taken in the weak operator topology; and \item the POVM on \(\R^s\) defined by \begin{equation} \label{eq:limiting-povm-mass-transfer} M^{(\infty)}(B) := \overline M(B) + \delta_0(B)\, \overline M(\overline{\R}^{\,s}\setminus\R^s), \qquad B\in\B(\R^s). \end{equation} \end{enumerate} Moreover, the sandwiched Le Cam third lemma, exactly as used in Step~3 of the cited proof, gives \begin{equation} \label{eq:rational-half-box-limits} \lim_{k\to\infty} \cL_h^{(n_k)}((-\infty,\alpha]) = \varrho_h(\widetilde M_\alpha), \qquad h\in\R^d,\quad \alpha\in\Q^s. \end{equation} Fix \(h\in\R^d\), and define \[ \mu_h(B) := \varrho_h\bigl(M^{(\infty)}(B)\bigr), \qquad B\in\B(\R^s). \] We prove that \(\cL_h^{(n_k)}\weak\mu_h\). For \(\alpha\in\Q^s\), set \[ a_h(\alpha):=\varrho_h(\widetilde M_\alpha). \] We first show that \(\overline M\) places no mass at infinity under the state \(\varrho_h\). Let \(\varepsilon>0\). By tightness of \(\{\cL_h^{(n)}\}_{n\geq1}\), there exists \(R>0\) such that \[ \cL_h^{(n)}([-R,R]^s)\geq1-\varepsilon \qquad\text{for every }n\geq1. \] Choose \(r\in\Q\) with \(r>R\), and let \(\alpha_r=(r,\ldots,r)\in\Q^s\). Since \([-R,R]^s\subset(-\infty,\alpha_r]\), \eqref{eq:rational-half-box-limits} gives \[ a_h(\alpha_r) = \lim_{k\to\infty} \cL_h^{(n_k)}((-\infty,\alpha_r]) \geq1-\varepsilon. \] Because \(\alpha\mapsto\widetilde M_\alpha\) is coordinatewise nondecreasing, \eqref{eq:bar-M-from-rational} implies \[ \overline M_{\alpha_r}\geq\widetilde M_{\alpha_r}, \] and hence \[ \varrho_h(\overline M_{\alpha_r})\geq1-\varepsilon. \] Letting \(r\to\infty\) through the rationals and using normality of \(\varrho_h\), we obtain \[ \varrho_h(\overline M(\R^s))=1. \] Therefore, \begin{equation} \label{eq:no-mass-at-infinity} \varrho_h\bigl( \overline M(\overline{\R}^{\,s}\setminus\R^s) \bigr) =0. \end{equation} It follows from \eqref{eq:limiting-povm-mass-transfer} and \eqref{eq:no-mass-at-infinity} that \[ \mu_h(B)=\varrho_h(\overline M(B)), \qquad B\in\B(\R^s). \] Consequently, for every \(t\in\R^s\), normality of \(\varrho_h\) and \eqref{eq:bar-M-from-rational} yield \begin{align} \mu_h((-\infty,t]) &= \varrho_h(\overline M_t) \notag\\ &= \inf_{\substack{\alpha\in\Q^s\\\alpha>t}} \varrho_h(\widetilde M_\alpha) \notag\\ &= \inf_{\substack{\alpha\in\Q^s\\\alpha>t}} a_h(\alpha). \label{eq:mu-right-regularization} \end{align} Equations \eqref{eq:rational-half-box-limits} and \eqref{eq:mu-right-regularization} verify the hypotheses of Lemma~\ref{lem:rational-lower-orthants} with \[ P_k=\cL_h^{(n_k)}, \qquad P=\mu_h, \qquad a(\alpha)=a_h(\alpha). \] Hence \[ \cL_h^{(n_k)}\weak\mu_h. \] Since \(h\in\R^d\) was arbitrary and the subsequence \(\{n_k\}\) was chosen independently of \(h\), this proves the assertion for every fixed local alternative \(h\). \end{proof}


\begin{proof}[Proof of Theorem~\ref{thm:qGausTest_HS1}]
Write the mean parameter as $h=(y,z)$, where
$y\in\R$ is the parameter of interest and
$z\in\R^{d-1}$ is the nuisance parameter. Let
\[
\mathcal M=\operatorname{CCR}(S),
\qquad
\phi_{y,z}\sim N((y,z),\Sigma),
\qquad
\Sigma=V+\i S,
\]
where $V\succ0$. Throughout, a test is an effect
$T\in\mathcal M$.

Consider the nuisance-translation group
$G=(\R^{d-1},+)$, acting on $\mathcal M$ through
\[
\tau_t:=\alpha_{(0,t)},
\qquad t\in\R^{d-1},
\]
where $\alpha_a$ denotes the normal automorphism
constructed in Appendix~\ref{sec:vn_alg}.
By the Gaussian covariance relation
\eqref{gauss_translation},
\begin{equation}
\phi_{y,z}(\tau_t(A))
=
\phi_{y,z+t}(A),
\qquad A\in\mathcal M.
\label{eq:gaussian-nuisance-covariance}
\end{equation}
The action is point-ultraweakly continuous,
and hence satisfies the hypotheses of
Theorem~\ref{thm:hunt-stein}.

Let $T\in\mathcal M$ be any level-$\alpha$
test. By Theorem~\ref{thm:hunt-stein},
there exists an invariant level-$\alpha$
test $\bar T\in\mathcal M$ such that
\begin{equation}
\inf_{z\in\R^{d-1}}\phi_{y,z}(T)
\le
\inf_{z\in\R^{d-1}}\phi_{y,z}(\bar T),
\qquad y>0.
\label{eq:gaussian-hs-minimum-power}
\end{equation}
Since the nuisance translations act transitively
on $\R^{d-1}$, the power of $\bar T$ is
independent of $z$.

By Lemma~\ref{lem:Gaussian-invariant-algebra},
the invariant algebra is
\[
\mathcal M^G
=
\mathcal A[X_1]
=
\{e^{\i tX_1}:t\in\R\}''.
\]
Thus every invariant test has the form
$\bar T=f(X_1)$ for some Borel function
$f:\R\to[0,1]$.

The Gaussian characteristic function gives
\[
\phi_{y,z}(e^{\i tX_1})
=
\exp\left\{
\i ty-\frac12t^2V_{11}
\right\}.
\]
Consequently, the restriction of $\phi_{y,z}$
to $\mathcal A[X_1]$ is the classical normal
location model
\[
X_1\sim N(y,V_{11}),
\]
which does not depend on the nuisance parameter
$z$. In particular,
\begin{equation}
\phi_{y,z}(\bar T)
=
\int_{\R}f(x)\,N(y,V_{11})(dx).
\label{eq:gaussian-invariant-classical-power}
\end{equation}

Since $\bar T$ has level $\alpha$, its
restriction satisfies
\[
\int_{\R}f(x)\,N(0,V_{11})(dx)
\le\alpha.
\]
The normal location family has a monotone
likelihood ratio. Therefore, by the
Neyman--Pearson lemma, for every $y>0$,
\[
\int_{\R}f(x)\,N(y,V_{11})(dx)
\le
1-\Phi\left(
z_\alpha-\frac{y}{\sqrt{V_{11}}}
\right),
\]
where $z_\alpha=\Phi^{-1}(1-\alpha)$.

Combining this inequality with
\eqref{eq:gaussian-hs-minimum-power} gives
\[
\inf_{z\in\R^{d-1}}\phi_{y,z}(T)
\le
1-\Phi\left(
z_\alpha-\frac{y}{\sqrt{V_{11}}}
\right),
\qquad y>0.
\]

Finally, consider the test
\[
T_\alpha
=
\left\{
X_1>z_\alpha\sqrt{V_{11}}
\right\}.
\]
By the spectral theorem,
$T_\alpha\in\mathcal A[X_1]\subseteq\mathcal M$.
Its rejection probability is
\[
\phi_{y,z}(T_\alpha)
=
1-\Phi\left(
z_\alpha-\frac{y}{\sqrt{V_{11}}}
\right),
\]
independently of $z$. This probability is
increasing in $y$ and equals $\alpha$ at
$y=0$. Hence $T_\alpha$ has level $\alpha$
for the composite null $y\le0$ and attains
the preceding bound for every $y>0$.
\end{proof}

\begin{proof}[Proof of Theorem~\ref{thm:qGausTest_two}]
As in the proof of Theorem~\ref{thm:qGausTest_HS},
choose an invertible real matrix $L$ whose first row
is $c^T$. Introduce the transformed observables
$\widetilde X=LX$ and mean parameter
$\widetilde h=Lh=(y,z)$, where
\[
Y:=\widetilde X_1=c^TX,
\qquad y=c^Th,\qquad z\in\R^{d-1}.
\]
The transformed covariance is
$\widetilde\Sigma=L\Sigma L^T$, and
\[
\sigma^2:=c^TVc=\widetilde\Sigma_{11}>0.
\]
Since $L$ is invertible, the transformed observables
generate the same algebra $\mathcal M=\operatorname{CCR}(S)$.
The testing problem becomes
\[
H_0:y=0,\qquad H_1:y\ne0.
\]

To distinguish the transformed parametrization from
the original Gaussian shift family, write
\[
\widetilde\phi_{y,z}
:=
\phi_{L^{-1}(y,z),\Sigma},
\qquad (y,z)\in\R\times\R^{d-1}.
\]
Equivalently, $\widetilde\phi_{y,z}$ is the Gaussian
state with mean $(y,z)$ and covariance
$\widetilde\Sigma$ in the transformed coordinates.

Let $T\in\mathcal M$ be any level-$\alpha$ unbiased
test. The nuisance translations act through
$\tau_t=\alpha_{L^{-1}(0,t)}$, $t\in\R^{d-1}$,
and satisfy
\[
\widetilde\phi_{y,z}(\tau_t(A))
=
\widetilde\phi_{y,z+t}(A),
\qquad A\in\mathcal M.
\]
By Theorem~\ref{thm:hunt-stein}, there exists an
invariant level-$\alpha$ unbiased test $\bar T$
such that
\begin{equation}
\inf_{z\in\R^{d-1}}\widetilde\phi_{y,z}(T)
\le
\inf_{z\in\R^{d-1}}\widetilde\phi_{y,z}(\bar T),
\qquad y\ne0.
\label{eq:two-sided-hs-reduction}
\end{equation}
By Lemma~\ref{lem:Gaussian-invariant-algebra},
applied in the transformed coordinates, every
invariant test belongs to
\[
\mathcal M^G
=
\{e^{\i tY}:t\in\R\}''.
\]
Thus $\bar T=f(Y)$ for some Borel function
$f:\R\to[0,1]$. Since
\[
\widetilde\phi_{y,z}(e^{\i tY})
=
\exp\left(\i ty-\frac12t^2\sigma^2\right),
\]
the invariant testing problem reduces to the
classical normal location family
$
Y\sim N(y,\sigma^2),
$
independently of the nuisance parameter $z$.

By the classical UMPU theorem for the normal
location family, the level-$\alpha$ unbiased
test maximizing power against every $y\ne0$ is
\[
T_\alpha
=
\left\{|Y|>z_{\alpha/2}\sigma\right\},
\]
where $z_{\alpha/2}=\Phi^{-1}(1-\alpha/2)$.
Its power function is
\[
\begin{aligned}
\beta_\alpha(y)
&:=\widetilde\phi_{y,z}(T_\alpha)\\
&=
1-\Phi\left(z_{\alpha/2}-\frac{y}{\sigma}\right)
+\Phi\left(-z_{\alpha/2}-\frac{y}{\sigma}\right).
\end{aligned}
\]
Consequently, \eqref{eq:two-sided-hs-reduction}
implies
\[
\inf_{z\in\R^{d-1}}\widetilde\phi_{y,z}(T)
\le\beta_\alpha(y),
\qquad y\ne0.
\]

The function $\beta_\alpha$ is even, satisfies
$\beta_\alpha(0)=\alpha$, and is nondecreasing
in $|y|$. Hence, for every $s>0$,
\[
\begin{aligned}
\inf_{\substack{|y|\ge s\\z\in\R^{d-1}}}
\widetilde\phi_{y,z}(T)
&\le
\inf_{z\in\R^{d-1}}\widetilde\phi_{s,z}(T)\\
&\le\beta_\alpha(s).
\end{aligned}
\]
The proposed test $T_\alpha$ is invariant,
level-$\alpha$, and unbiased, and its minimum
power over $|y|\ge s$ equals $\beta_\alpha(s)$.
Thus the bound is attained.

Returning to the original coordinates,
where $y=c^Th$ and $\sigma^2=c^TVc$,
proves the theorem.
\end{proof}

\begin{proof}[Proof of Theorem~\ref{thm:qLAN-test}]
Set
\[
c:=\dot\psi_{\theta_0}^T,\qquad
\sigma^2:=c^T\Sigma c>0,\qquad
b_\alpha(t):=
1-\Phi\left(z_\alpha-\frac{t}{\sigma}\right),
\quad t\in\R.
\]
For $h\in\R^d$, write
$\beta_n(h):=
\Tr\rho_{\theta_0+h/\sqrt n}^{(n)}T^{(n)}$. We first prove the upper bound. Suppose,
to the contrary, that it fails for some $s>0$.
Then there exists $\delta>0$ such that
\begin{equation}
\label{eq:strict-maximin-gap}
\inf_{\substack{D\subset\R^d\\0<|D|<\infty}}
\limsup_{n\to\infty}
\inf_{h\in D\cap H_{1,s}^{(n)}}
\beta_n(h)
>
b_\alpha(s)+2\delta.
\end{equation}

By continuity of $b_\alpha$, choose
$\epsilon>0$ sufficiently small that
$b_\alpha(s+\epsilon)<b_\alpha(s)+\delta$.
Put
\[
H_{1,s,\epsilon}
:=
\{h\in\R^d:c^Th\ge s+\epsilon\}.
\]
Since $H_{1,s,\epsilon}$ is closed and
separable, choose a countable dense subset
$\{h_1,h_2,\ldots\}\subset H_{1,s,\epsilon}$
and define $D_m:=\{h_1,\ldots,h_m\}$.

For each fixed $m$, differentiability of
$\psi$ at $\theta_0$ gives
\[
\sqrt n\,
\psi\left(\theta_0+\frac{h}{\sqrt n}\right)
=c^Th+o(1)
\]
for every $h\in D_m$. Since $D_m$ is finite
and $c^Th\ge s+\epsilon$ on $D_m$, we have
$D_m\subseteq H_{1,s}^{(n)}$ for all
sufficiently large $n$. Consequently,
\eqref{eq:strict-maximin-gap} implies
\[
\limsup_{n\to\infty}
\min_{1\le j\le m}\beta_n(h_j)
>
b_\alpha(s)+2\delta
\qquad(m\ge1).
\]
We may therefore choose a strictly increasing
sequence $\{n_m\}_{m\ge1}$ such that
\begin{equation}
\label{eq:diagonal-power-lower-bound}
\beta_{n_m}(h_j)
\ge b_\alpha(s)+\delta,
\qquad j=1,\ldots,m.
\end{equation}

For each $m$, regard $T^{(n_m)}$ as the binary
POVM $\{T^{(n_m)},I-T^{(n_m)}\}$.
Its outcome distributions are supported on
$\{0,1\}$ and hence are tight under every
fixed local alternative. By
Theorem~\ref{thm:common_sub}, there exist
a further subsequence $\{n_{m_k}\}_{k\ge1}$
and a limiting Gaussian test
$0\le T^{(\infty)}\le I$ such that
\begin{equation}
\label{eq:limiting-test-power-convergence}
\beta_{n_{m_k}}(h)
\rightarrow
\varrho_h(T^{(\infty)})
\qquad\text{for every fixed }h\in\R^d.
\end{equation}

Fix $j\ge1$. Since $m_k\ge j$ eventually,
\eqref{eq:diagonal-power-lower-bound}
and \eqref{eq:limiting-test-power-convergence}
give
\begin{equation}
\label{eq:dense-power-lower-bound}
\varrho_{h_j}(T^{(\infty)})
\ge b_\alpha(s)+\delta,
\qquad j\ge1.
\end{equation}

\noindent We next verify that $T^{(\infty)}$
is a level-$\alpha$ test for the Gaussian
hypotheses $H_0:c^Th\le0$ versus $H_1:c^Th>0$.

Let $h\in\R^d$ satisfy $c^Th<0$.
Differentiability of $\psi$ and
$\psi(\theta_0)=0$ give
\[
\sqrt n\,
\psi\left(\theta_0+\frac{h}{\sqrt n}\right)
=c^Th+o(1)<0
\]
for all sufficiently large $n$.
Since $T^{(n)}$ has level $\alpha$ under
the original null, $\beta_n(h)\le\alpha$
eventually. Hence
\eqref{eq:limiting-test-power-convergence}
implies
\begin{equation}
\label{eq:strict-limiting-null}
\varrho_h(T^{(\infty)})\le\alpha,
\qquad c^Th<0.
\end{equation}

By Lemma~\ref{lem:gaussian-shift-continuity},
the map $h\mapsto\varrho_h(T^{(\infty)})$
is continuous. It follows that
\begin{equation}
\label{eq:limiting-test-level}
\varrho_h(T^{(\infty)})\le\alpha,
\qquad c^Th\le0,
\end{equation}
so $T^{(\infty)}$ is level $\alpha$ in
the limiting Gaussian experiment.

By the same continuity argument,
\eqref{eq:dense-power-lower-bound} extends
from the dense set $\{h_j\}_{j\ge1}$ to
all of $H_{1,s,\epsilon}$:
\begin{equation}
\label{eq:global-power-lower-bound}
\varrho_h(T^{(\infty)})
\ge b_\alpha(s)+\delta,
\qquad h\in H_{1,s,\epsilon}.
\end{equation}
In particular,
\[
\inf_{c^Th\ge s+\epsilon}
\varrho_h(T^{(\infty)})
\ge b_\alpha(s)+\delta.
\]
On the other hand, since $T^{(\infty)}$
is level $\alpha$,
Theorem~\ref{thm:qGausTest_HS} gives
\[
\inf_{c^Th\ge s+\epsilon}
\varrho_h(T^{(\infty)})
\le b_\alpha(s+\epsilon)
<b_\alpha(s)+\delta,
\]
contradicting the preceding inequality.
This proves the upper bound.

\medskip
\noindent We finally verify attainment of the
Gaussian local power envelope. Consider
\[
T_\alpha^{(n)}
=
\left\{
\dot\psi_{\theta_0}X^{(n)}
>z_\alpha\sigma
\right\}.
\]
Applying the sandwiched Le Cam third lemma
under $\mathcal D$-extendibility
\cite[Corollary~3.7]{Fujiwara2023},
we obtain, for every fixed $h\in\R^d$,
\[
\left(
\dot\psi_{\theta_0}X^{(n)},
\rho_{\theta_0+h/\sqrt n}^{(n)}
\right)
\weak N(c^Th,\sigma^2).
\]
Since the limiting normal distribution
assigns zero mass to $z_\alpha\sigma$,
\begin{equation}
\begin{aligned}
\beta_n^\star(h)
&:=
\Tr\left(
\rho_{\theta_0+h/\sqrt n}^{(n)}
T_\alpha^{(n)}
\right)\\
&\rightarrow
1-\Phi\left(
z_\alpha-\frac{c^Th}{\sigma}
\right)
=b_\alpha(c^Th).
\end{aligned}
\label{eq:asymptotic-wald-power}
\end{equation}
In particular, setting $h=0$ gives
$\lim_{n\to\infty}\beta_n^\star(0)=\alpha$.

\noindent We now show that these tests attain
the stated power bound for the original
$n$-dependent alternative sets.

Fix $s>0$, a nonempty finite set
$D\subset\R^d$, and $\eta\in(0,s)$.
Define
\[
D_{-\eta}:=
\{h\in D:c^Th\ge s-\eta\}.
\]
For every $h\in D\setminus D_{-\eta}$,
differentiability of $\psi$ gives
\[
\sqrt n\,
\psi\left(\theta_0+\frac{h}{\sqrt n}\right)
=c^Th+o(1)<s
\]
for all sufficiently large $n$.
Since $D$ is finite, it follows that
$D\cap H_{1,s}^{(n)}\subseteq D_{-\eta}$
eventually. Consequently,
\[
\inf_{h\in D\cap H_{1,s}^{(n)}}
\beta_n^\star(h)
\ge
\inf_{h\in D_{-\eta}}
\beta_n^\star(h)
\]
for all sufficiently large $n$.

If $D_{-\eta}\ne\varnothing$, its finiteness
and \eqref{eq:asymptotic-wald-power} imply

\begin{equation*}
\liminf_{n\to\infty}
\inf_{h\in D\cap H_{1,s}^{(n)}}
\beta_n^\star(h)\ge
\min_{h\in D_{-\eta}}b_\alpha(c^Th)
\ge b_\alpha(s-\eta).
\end{equation*}

If $D_{-\eta}=\varnothing$, then
$D\cap H_{1,s}^{(n)}$ is eventually empty
and the same lower bound holds by the
convention $\inf\varnothing=+\infty$.

Since $\eta>0$ is arbitrary and
$b_\alpha$ is continuous, letting
$\eta\downarrow0$ yields
\[
\liminf_{n\to\infty}
\inf_{h\in D\cap H_{1,s}^{(n)}}
\beta_n^\star(h)
\ge b_\alpha(s).
\]
Since $D$ was arbitrary,
\begin{equation}
\label{eq:attainment-lower}
\inf_{\substack{D\subset\R^d\\0<|D|<\infty}}
\limsup_{n\to\infty}
\inf_{h\in D\cap H_{1,s}^{(n)}}
\beta_n^\star(h)
\ge b_\alpha(s).
\end{equation}

\noindent Conversely, fix $\eta>0$ and put
$h_\eta:=(s+\eta)c/\|c\|^2$.
Then $c^Th_\eta=s+\eta$, and
differentiability gives
\[
\sqrt n\,
\psi\left(
\theta_0+\frac{h_\eta}{\sqrt n}
\right)
=s+\eta+o(1)>s
\]
for all sufficiently large $n$.
Thus $h_\eta\in H_{1,s}^{(n)}$ eventually.

Taking $D=\{h_\eta\}$ and using
\eqref{eq:asymptotic-wald-power}, we get

\begin{equation*}
\inf_{\substack{D\subset\R^d\\0<|D|<\infty}}
\limsup_{n\to\infty}
\inf_{h\in D\cap H_{1,s}^{(n)}}
\beta_n^\star(h)\le
\lim_{n\to\infty}\beta_n^\star(h_\eta)
=b_\alpha(s+\eta).
\end{equation*}
Letting $\eta\downarrow0$ and using
continuity of $b_\alpha$ gives
\begin{equation}
\label{eq:attainment-upper}
\inf_{\substack{D\subset\R^d\\0<|D|<\infty}}
\limsup_{n\to\infty}
\inf_{h\in D\cap H_{1,s}^{(n)}}
\beta_n^\star(h)
\le b_\alpha(s).
\end{equation}
Combining \eqref{eq:attainment-lower}
and \eqref{eq:attainment-upper}, we obtain
\[
\inf_{\substack{D\subset\R^d\\0<|D|<\infty}}
\limsup_{n\to\infty}
\inf_{h\in D\cap H_{1,s}^{(n)}}
\Tr\left(
\rho_{\theta_0+h/\sqrt n}^{(n)}
T_\alpha^{(n)}
\right)
=b_\alpha(s).
\]
Thus the proposed scalar tests attain
the Gaussian local power envelope in the
stated finite-submodel criterion for every $s>0$.
\end{proof}

\begin{proof}[Proof of Theorem~\ref{thm:qLAN-test2}]
Set
\[
c:=\dot\psi_{\theta_0}^{T},
\qquad
\sigma^2:=c^T\Sigma c>0,
\]
and define
\[
b_\alpha(t)
:=
1-\Phi\left(
-z_{\alpha/2}-\frac{t}{\sigma},
z_{\alpha/2}-\frac{t}{\sigma}
\right),\qquad t\in\R.
\]
For $h\in\R^d$, write
$\beta_n(h):=
\Tr\rho_{\theta_0+h/\sqrt n}^{(n)}T^{(n)}$.
We first prove the upper bound.
Suppose, to the contrary, that it fails
for some $s>0$. Then there exists
$\delta>0$ such that
\begin{equation}
\label{eq:strict-maximin-gap-two-sided}
\inf_{\substack{D\subset\R^d\\0<|D|<\infty}}
\limsup_{n\to\infty}
\inf_{h\in D\cap H_{1,s}^{(n)}}
\beta_n(h)
>
b_\alpha(s)+3\delta.
\end{equation}
By continuity of $b_\alpha$, choose
$\eta>0$ sufficiently small that
\begin{equation}
\label{eq:two-sided-eta-choice}
b_\alpha(s+\eta)<b_\alpha(s)+\delta.
\end{equation}

Define
\[
A_{s+\eta}
:=
\{h\in\R^d:|c^Th|\ge s+\eta\}.
\]
Since $A_{s+\eta}$ is closed and
separable, choose a countable dense
subset $\{h_1,h_2,\ldots\}\subseteq
A_{s+\eta}$ and put
$D_m:=\{h_1,\ldots,h_m\}$.

For each fixed $h_j\in A_{s+\eta}$,
differentiability of $\psi$ gives
\[
\sqrt n\,
\psi\left(\theta_0+\frac{h_j}{\sqrt n}\right)
=c^Th_j+o(1).
\]
Since $|c^Th_j|\ge s+\eta$,
we have $h_j\in H_{1,s}^{(n)}$
for all sufficiently large $n$.
Because $D_m$ is finite,
$D_m\subseteq H_{1,s}^{(n)}$
eventually for each fixed $m$.
Thus \eqref{eq:strict-maximin-gap-two-sided}
implies
\[
\limsup_{n\to\infty}
\min_{1\le j\le m}\beta_n(h_j)
>
b_\alpha(s)+3\delta,
\qquad m\ge1.
\]
Consequently, we may choose a strictly
increasing sequence $\{n_m\}_{m\ge1}$
such that
\begin{equation}
\label{eq:diagonal-power-lower-bound-two-sided}
\beta_{n_m}(h_j)
\ge b_\alpha(s)+2\delta,
\qquad j=1,\ldots,m.
\end{equation}

For each $m$, regard $T^{(n_m)}$
as the binary POVM
$\{T^{(n_m)},I-T^{(n_m)}\}$.
Its outcome distributions are
supported on $\{0,1\}$ and hence
are tight under every fixed local
alternative. By Theorem~\ref{thm:common_sub},
there exist a further subsequence
$\{n_{m_k}\}_{k\ge1}$ and a limiting
Gaussian test $0\le T^{(\infty)}\le I$
such that
\begin{equation}
\label{eq:limiting-test-power-convergence-two-sided}
\beta_{n_{m_k}}(h)
\rightarrow
\varrho_h(T^{(\infty)})
\qquad\text{for every fixed }h\in\R^d.
\end{equation}

Fix $j\ge1$. Since $m_k\ge j$
eventually, the preceding two displays
give
\begin{equation}
\label{eq:dense-power-lower-bound-two-sided}
\varrho_{h_j}(T^{(\infty)})
\ge b_\alpha(s)+2\delta,
\qquad j\ge1.
\end{equation}
By Lemma~\ref{lem:gaussian-shift-continuity},
the map
$h\mapsto\varrho_h(T^{(\infty)})$
is continuous. Hence this inequality
extends to all of $A_{s+\eta}$:
\begin{equation}
\label{eq:global-power-lower-bound-two-sided}
\varrho_h(T^{(\infty)})
\ge b_\alpha(s)+2\delta,
\qquad |c^Th|\ge s+\eta.
\end{equation}

\noindent We next verify that
$T^{(\infty)}$ is a level-$\alpha$
unbiased test for the Gaussian
hypotheses $H_0:c^Th=0$ versus
$H_1:c^Th\ne0$.

Let $h\in\R^d$ satisfy $c^Th\ne0$.
Differentiability of $\psi$ at
$\theta_0$, together with
$\psi(\theta_0)=0$, gives
\[
\sqrt n\,
\psi\left(\theta_0+\frac{h}{\sqrt n}\right)
=c^Th+o(1).
\]
Thus
$\psi(\theta_0+h/\sqrt n)\ne0$
for all sufficiently large $n$.
Since $T^{(n)}$ is unbiased for
the original two-sided testing problem,
$\beta_n(h)\ge\alpha$ eventually.
By \eqref{eq:limiting-test-power-convergence-two-sided},
\begin{equation}
\label{eq:limiting-test-unbiased-alternative}
\varrho_h(T^{(\infty)})
\ge\alpha,
\qquad c^Th\ne0.
\end{equation}

\noindent It remains to verify the
level condition on the limiting null
hyperplane. Fix $h\in\R^d$ such that
$c^Th=0$. By Lemma~\ref{seq_level_lemma},
there exists $h_n\to h$ satisfying
\[
\psi\left(
\theta_0+\frac{h_n}{\sqrt n}
\right)=0
\]
for all sufficiently large $n$.
The finite-sample level condition
therefore gives
\[
\Tr\left(
\rho_{\theta_0+h_n/\sqrt n}^{(n)}
T^{(n)}
\right)\le\alpha
\]
eventually.

By 
assumption \eqref{trace_conv_local},
\[
\left\|
\rho_{\theta_0+h_n/\sqrt n}^{(n)}
-
\rho_{\theta_0+h/\sqrt n}^{(n)}
\right\|_1\rightarrow0.
\]
Since $0\le T^{(n)}\le I$,
\[
\begin{aligned}
\beta_n(h)
&\le
\Tr\left(
\rho_{\theta_0+h_n/\sqrt n}^{(n)}
T^{(n)}
\right)+
\left\|
\rho_{\theta_0+h_n/\sqrt n}^{(n)}
-
\rho_{\theta_0+h/\sqrt n}^{(n)}
\right\|_1\\
&\le\alpha+o(1).
\end{aligned}
\]
Consequently,
$\limsup_{n\to\infty}\beta_n(h)\le\alpha$.
Passing to the extracted subsequence,
we obtain
\begin{equation}
\label{eq:limiting-test-level-null-two-sided}
\varrho_h(T^{(\infty)})
\le\alpha,
\qquad c^Th=0.
\end{equation}

Combining
\eqref{eq:limiting-test-unbiased-alternative}
and
\eqref{eq:limiting-test-level-null-two-sided},
we conclude that $T^{(\infty)}$
is level-$\alpha$ unbiased in the
limiting Gaussian experiment. Therefore,
Theorem~\ref{thm:qGausTest_two} gives
\[
\begin{aligned}
\inf_{|c^Th|\ge s+\eta}
\varrho_h(T^{(\infty)})
&\le
1-\Phi\left(
-z_{\alpha/2}-\frac{s+\eta}{\sigma},
z_{\alpha/2}-\frac{s+\eta}{\sigma}
\right)\\
&=b_\alpha(s+\eta)
<b_\alpha(s)+\delta,
\end{aligned}
\]
where the last inequality follows
from \eqref{eq:two-sided-eta-choice}.
This contradicts
\eqref{eq:global-power-lower-bound-two-sided}
and proves the asserted upper bound.

\medskip
\noindent We finally verify attainment
of the Gaussian local power envelope.
Consider
\[
T_\alpha^{(n)}
=
\left\{
\left|\dot\psi_{\theta_0}X^{(n)}\right|
>
z_{\alpha/2}\sigma
\right\}.
\]
The sandwiched Le Cam third lemma
under $\mathcal D$-extendibility
\cite[Corollary~3.7]{Fujiwara2023}
gives, for every fixed $h\in\R^d$,
\[
\left(
\dot\psi_{\theta_0}X^{(n)},
\rho_{\theta_0+h/\sqrt n}^{(n)}
\right)
\weak N(c^Th,\sigma^2).
\]
Since the limiting normal law
assigns zero mass to the rejection
thresholds, we obtain
\begin{equation}
\beta_n^\star(h)
:=
\Tr\left(
\rho_{\theta_0+h/\sqrt n}^{(n)}
T_\alpha^{(n)}
\right)\rightarrow
b_\alpha(c^Th).
\label{eq:asymptotic-wald-power-two-sided}
\end{equation}
The function $b_\alpha$ is even,
nondecreasing in $|t|$, and satisfies
$b_\alpha(0)=\alpha$. In particular,
$\lim_{n\to\infty}\beta_n^\star(0)=\alpha$.

\noindent We now show attainment for the
original $n$-dependent alternative sets
\[
H_{1,s}^{(n)}
=
\left\{
h\in\R^d:
\left|
\sqrt n\,
\psi\left(\theta_0+\frac{h}{\sqrt n}\right)
\right|\ge s
\right\}.
\]
Fix $s>0$, a nonempty finite set
$D\subset\R^d$, and
$\epsilon\in(0,s)$. Put
\[
D_{-\epsilon}
:=
\{h\in D:|c^Th|\ge s-\epsilon\}.
\]
For every $h\in D\setminus D_{-\epsilon}$,
differentiability gives
\[
\sqrt n\,
\psi\left(\theta_0+\frac{h}{\sqrt n}\right)
=c^Th+o(1),
\]
so
$\left|\sqrt n\,\psi(
\theta_0+h/\sqrt n)\right|<s$
for all sufficiently large $n$.
Since $D$ is finite,
\[
D\cap H_{1,s}^{(n)}
\subseteq D_{-\epsilon}
\]
eventually. Hence
\[
\inf_{h\in D\cap H_{1,s}^{(n)}}
\beta_n^\star(h)
\ge
\inf_{h\in D_{-\epsilon}}
\beta_n^\star(h)
\]
for all sufficiently large $n$.

If $D_{-\epsilon}\ne\varnothing$,
its finiteness and
\eqref{eq:asymptotic-wald-power-two-sided}
give

\begin{equation*}
\liminf_{n\to\infty}
\inf_{h\in D\cap H_{1,s}^{(n)}}
\beta_n^\star(h)\ge
\min_{h\in D_{-\epsilon}}
b_\alpha(c^Th)
\ge b_\alpha(s-\epsilon).
\end{equation*}

If $D_{-\epsilon}=\varnothing$,
then $D\cap H_{1,s}^{(n)}$
is eventually empty and the same
inequality holds under the convention
$\inf\varnothing=+\infty$. Letting $\epsilon\downarrow0$,
we obtain
\[
\liminf_{n\to\infty}
\inf_{h\in D\cap H_{1,s}^{(n)}}
\beta_n^\star(h)
\ge b_\alpha(s).
\]
Since $D$ was arbitrary,
\begin{equation}
\label{eq:attainment-lower-two-sided}
\inf_{\substack{D\subset\R^d\\0<|D|<\infty}}
\limsup_{n\to\infty}
\inf_{h\in D\cap H_{1,s}^{(n)}}
\beta_n^\star(h)
\ge b_\alpha(s).
\end{equation}

\noindent Conversely, fix $\epsilon>0$
and put $h_\epsilon:=(s+\epsilon)c/\|c\|^2$.
Then $c^Th_\epsilon=s+\epsilon$,
and differentiability yields
\[
\sqrt n\,
\psi\left(
\theta_0+\frac{h_\epsilon}{\sqrt n}
\right)
=s+\epsilon+o(1).
\]
Thus $h_\epsilon\in H_{1,s}^{(n)}$
for all sufficiently large $n$.
Taking $D=\{h_\epsilon\}$ and using
\eqref{eq:asymptotic-wald-power-two-sided},
we obtain

\begin{equation*}
\inf_{\substack{D\subset\R^d\\0<|D|<\infty}}
\limsup_{n\to\infty}
\inf_{h\in D\cap H_{1,s}^{(n)}}
\beta_n^\star(h)\le
\lim_{n\to\infty}
\beta_n^\star(h_\epsilon)
=b_\alpha(s+\epsilon).
\end{equation*}

Letting $\epsilon\downarrow0$
and using continuity of $b_\alpha$ gives
\begin{equation}
\label{eq:attainment-upper-two-sided}
\inf_{\substack{D\subset\R^d\\0<|D|<\infty}}
\limsup_{n\to\infty}
\inf_{h\in D\cap H_{1,s}^{(n)}}
\beta_n^\star(h)
\le b_\alpha(s).
\end{equation}

Combining
\eqref{eq:attainment-lower-two-sided}
and
\eqref{eq:attainment-upper-two-sided},
we conclude that

\begin{equation*}
\inf_{\substack{D\subset\R^d\\0<|D|<\infty}}
\limsup_{n\to\infty}
\inf_{h\in D\cap H_{1,s}^{(n)}}
\Tr\left(
\rho_{\theta_0+h/\sqrt n}^{(n)}
T_\alpha^{(n)}
\right)=b_\alpha(s).
\end{equation*}

Thus the proposed scalar tests attain
the Gaussian local power envelope
for the original $n$-dependent
alternative sets $H_{1,s}^{(n)}$.
\end{proof}

\begin{proof}[Proof of Theorem~\ref{thm:qLAN-test3_exact_level}]
The proof follows the asymptotic minimax argument
of \cite{Fujiwara2023}.

Set
\[
c:=\dot\psi_{\theta_0}^T,\qquad
\sigma^2:=c^T\Sigma c>0,
\]
and define
\[
b_\alpha(t):=
1-\Phi\left(
-z_{\alpha/2}-\frac{t}{\sigma},
z_{\alpha/2}-\frac{t}{\sigma}
\right),\qquad t\in\R.
\]
For $h\in\R^d$, write
$\beta_n(h):=
\Tr\rho_{\theta_0+h/\sqrt n}^{(n)}T^{(n)}$.

\noindent We first prove the upper bound.
Fix $s>0$ and choose a countable dense
subset $\{h_1,h_2,\ldots\}$ of $H_{1,s}$.
Let $D_k:=\{h_1,\ldots,h_k\}$ and put
\[
\gamma_k^{(n)}
:=\min_{h\in D_k}\beta_n(h),\qquad
\gamma_k:=\limsup_{n\to\infty}
\gamma_k^{(n)},\qquad
\gamma_\infty:=\lim_{k\to\infty}\gamma_k.
\]
The last limit exists because $\gamma_k$
is nonincreasing and belongs to $[0,1]$.
Since each $D_k$ is an admissible finite
subset of $H_{1,s}$, the left-hand side of
\eqref{eq:qlan-tangent-two-sided-bound}
is at most $\gamma_\infty$.
It therefore suffices to establish
$\gamma_\infty\le b_\alpha(s)$.

Choose a strictly increasing sequence
$(n_k)$ such that
$\gamma_k^{(n_k)}-\gamma_k\to0$.
Then $\gamma_k^{(n_k)}\to\gamma_\infty$.

The binary POVMs
$\{T^{(n_k)},I-T^{(n_k)}\}$
are locally tight. By
Theorem~\ref{thm:common_sub}, there exist
a further subsequence $(n_{k_t})$
and a Gaussian test
$0\le T^{(\infty)}\le I$ such that
\begin{equation}
\label{eq:tangent-limiting-power}
\beta_{n_{k_t}}(h)
\rightarrow
\varrho_h(T^{(\infty)})
\qquad\text{for every fixed }h\in\R^d,
\end{equation}
where $\varrho_h=\phi_{h,\Sigma}$.

By the finite-sample level and unbiasedness
conditions, for every fixed $h$,
$\beta_n(h)\le\alpha$ eventually when
$c^Th=0$, and $\beta_n(h)\ge\alpha$
eventually when $c^Th\ne0$.
Passing to the limit in
\eqref{eq:tangent-limiting-power} gives
\[
\varrho_h(T^{(\infty)})\le\alpha
\quad(c^Th=0),\qquad
\varrho_h(T^{(\infty)})\ge\alpha
\quad(c^Th\ne0).
\]
Thus $T^{(\infty)}$ is a level-$\alpha$
unbiased test in the limiting Gaussian
experiment. By
Lemma~\ref{lem:gaussian-shift-continuity},
its power function
$h\mapsto\varrho_h(T^{(\infty)})$
is continuous.

For every fixed $h\in\bigcup_kD_k$,
we have $h\in D_{k_t}$ eventually.
Consequently,
\[
\begin{aligned}
\gamma_\infty
&=\lim_{t\to\infty}
\gamma_{k_t}^{(n_{k_t})}\\
&=\lim_{t\to\infty}
\min_{u\in D_{k_t}}\beta_{n_{k_t}}(u)\\
&\le\lim_{t\to\infty}
\beta_{n_{k_t}}(h)
=\varrho_h(T^{(\infty)}).
\end{aligned}
\]
Density and continuity extend this
inequality to every $h\in H_{1,s}$.
Therefore, by
Theorem~\ref{thm:qGausTest_two},
\[
\gamma_\infty
\le
\inf_{h\in H_{1,s}}
\varrho_h(T^{(\infty)})
\le b_\alpha(s).
\]
This proves the upper bound in
\eqref{eq:qlan-tangent-two-sided-bound}.

\medskip
\noindent We next verify attainment of
the Gaussian local power envelope.
Consider
\[
T_\alpha^{(n)}
=
\left\{
\left|\dot\psi_{\theta_0}X^{(n)}\right|
>z_{\alpha/2}\sigma
\right\}.
\]
By the sandwiched Le Cam third lemma under
$\mathcal D$-extendibility
\cite[Corollary~3.7]{Fujiwara2023},
for every fixed $h\in\R^d$,
\[
\left(
\dot\psi_{\theta_0}X^{(n)},
\rho_{\theta_0+h/\sqrt n}^{(n)}
\right)
\weak N(c^Th,\sigma^2).
\]
Since the limiting normal distribution
assigns zero mass to the rejection
thresholds, we obtain
\begin{equation}
\begin{aligned}
\beta_n^\star(h)
&:=
\Tr\left(
\rho_{\theta_0+h/\sqrt n}^{(n)}
T_\alpha^{(n)}
\right)\\
&\rightarrow
b_\alpha(c^Th).
\end{aligned}
\label{eq:tangent-asymptotic-power}
\end{equation}

The function $b_\alpha$ is even,
nondecreasing in $|t|$, and satisfies
$b_\alpha(0)=\alpha$.
Consequently, the proposed tests have
limiting rejection probability $\alpha$
for every fixed $h$ with $c^Th=0$,
and at least $\alpha$ for every fixed
$h$ with $c^Th\ne0$.

For every nonempty finite
$D\subset H_{1,s}$, finiteness
and \eqref{eq:tangent-asymptotic-power}
give
\[
\begin{aligned}
\lim_{n\to\infty}
\min_{h\in D}\beta_n^\star(h)
&=\min_{h\in D}b_\alpha(c^Th)\\
&\ge b_\alpha(s).
\end{aligned}
\]
Conversely, the vector
$h_s:=sc/\|c\|^2$ belongs to $H_{1,s}$
and satisfies $c^Th_s=s$.
Taking $D=\{h_s\}$ gives
\[
\lim_{n\to\infty}\beta_n^\star(h_s)
=b_\alpha(s).
\]
Hence
\[
\inf_{\substack{\varnothing\ne D\subset H_{1,s}\\
                 |D|<\infty}}
\limsup_{n\to\infty}
\inf_{h\in D}\beta_n^\star(h)
=b_\alpha(s),
\]
which proves attainment of the
Gaussian local power envelope for
every $s>0$.
\end{proof}

\section{Proofs of the technical lemmas}

\begin{lemma}[Continuity of Gaussian-shift power functions]
\label{lem:gaussian-shift-continuity}
Let $S\in\mathbb R^{m\times m}$ be real skew-symmetric
and $\Sigma=V+iS\succeq0$, with $V\succ0$.
Let $\mathcal M=\operatorname{CCR}(S)$ be the von Neumann
algebra constructed in Appendix~\ref{sec:vn_alg}, and
let $B\in\mathbb R^{m\times d}$. For $h\in\R^d$, let $\phi_{Bh,\Sigma}$ denote
the normal Gaussian state with mean $Bh$
and covariance $\Sigma$. Then
\begin{equation}
\lim_{g\to h}
\sup_{\substack{A\in\mathcal M\\\|A\|\le1}}
\left|
\phi_{Bg,\Sigma}(A)
-
\phi_{Bh,\Sigma}(A)
\right|=0.\label{eq:gaussian-shift-norm-continuity}
\end{equation}
In particular, $h\mapsto\phi_{Bh,\Sigma}(T)$ is continuous
for every bounded $T\in\mathcal M$, including
every test $0\le T\le I$.
\end{lemma}

\begin{proof}
Recall from Appendix~\ref{sec:vn_alg} that
$\mathcal M$ acts on a separable Hilbert space
$\widetilde{\mathcal H}$ and that translations of
the Gaussian mean are implemented by strongly
continuous displacement unitaries $U(a)$, with
\[
\alpha_a(A)=U(a)^*AU(a),
\qquad
\phi_{Bh,\Sigma}=\phi_{0,\Sigma}\circ\alpha_{Bh}.
\]

Since $\phi_{0,\Sigma}$ is a normal state on the concrete
von Neumann algebra $\mathcal M$, the normal-state
representation theorem
\cite[Theorem~7.1.12]{KadisonRingroseII}
provides a positive trace-class operator $K_0$
on $\widetilde{\mathcal H}$ such that
\[
\phi_{0,\Sigma}(A)=\Tr(K_0A),
\qquad A\in\mathcal M.
\]
Consequently,
\[
\phi_{Bh,\Sigma}(A)
=
\Tr\!\left(
U(Bh)K_0U(Bh)^*A
\right).
\]

For every $A\in\mathcal M$ with $\|A\|\le1$,
the trace inequality
$|\Tr(KA)|\le\|K\|_1\|A\|$ therefore gives
\begin{equation}
\begin{aligned}
|\phi_{Bg,\Sigma}(A)-\phi_{Bh,\Sigma}(A)|
&\le
\bigl\|
U(Bg)K_0U(Bg)^*
-
U(Bh)K_0U(Bh)^*
\bigr\|_1.
\end{aligned}
\label{eq:gaussian-continuity-trace-bound}
\end{equation}

We show that the right-hand side tends to zero.
For a rank-one operator $K=|u\rangle\langle v|$,
the identity
$\||u\rangle\langle v|\|_1=\|u\|\|v\|$
and the triangle inequality yield
\begin{align*}
&\|U(Bg)KU(Bg)^*-U(Bh)KU(Bh)^*\|_1\\
&\quad\le
\|(U(Bg)-U(Bh))u\|\,\|v\|\\
&\qquad+
\|u\|\,\|(U(Bg)-U(Bh))v\|.
\end{align*}
This tends to zero as $g\to h$ by strong
continuity of $U$. By linearity, the same holds
for every finite-rank operator.

Finally, let $K$ be finite-rank and approximate
$K_0$ in trace norm. Unitary invariance gives
\begin{align*}
&\|U(Bg)K_0U(Bg)^*
    -U(Bh)K_0U(Bh)^*\|_1\\
&\quad\le
2\|K_0-K\|_1\\
&\qquad+
\|U(Bg)KU(Bg)^*
    -U(Bh)KU(Bh)^*\|_1.
\end{align*}
Since finite-rank operators are dense in the
trace class, the right-hand side tends to zero
as $g\to h$. Taking the supremum over $\|A\|\le1$ in
\eqref{eq:gaussian-continuity-trace-bound}
proves \eqref{eq:gaussian-shift-norm-continuity}.
The asserted continuity of $\phi_{Bh,\Sigma}(T)$ follows
immediately.
\end{proof}

\begin{lemma} \label{lem:rational-lower-orthants} Let \(P_n\), \(n\geq1\), and \(P\) be probability measures on \(\R^m\). For \(t=(t_1,\ldots,t_m)\in\R^m\), write \[ (-\infty,t] := \prod_{j=1}^m(-\infty,t_j]. \] Suppose that, for every \(q\in\Q^m\), the limit \[ a(q) := \lim_{n\to\infty}P_n((-\infty,q]) \] exists. Assume moreover that \begin{equation} \label{eq:right-regularized-rational-limits} P((-\infty,t]) = \inf_{\substack{q\in\Q^m\\q>t}}a(q), \qquad t\in\R^m, \end{equation} where \(q>t\) means that \(q_j>t_j\) for every \(j=1,\ldots,m\). Then \[ P_n\weak P. \] \end{lemma}

\begin{proof}[Proof of Lemma \ref{lem:rational-lower-orthants}]
 Set \[ F_n(t):=P_n((-\infty,t]) \qquad\text{and}\qquad F(t):=P((-\infty,t]). \] It suffices to prove that \(F_n(t)\to F(t)\) at every continuity point \(t\) of \(F\). Fix such a point \(t\). For every \(q\in\Q^m\) satisfying \(q>t\), monotonicity gives \(F_n(t)\leq F_n(q)\). Therefore, \[ \limsup_{n\to\infty}F_n(t) \leq a(q). \] Taking the infimum over all rational \(q>t\) and using \eqref{eq:right-regularized-rational-limits}, we obtain \begin{equation} \label{eq:rational-orthant-upper-bound} \limsup_{n\to\infty}F_n(t) \leq F(t). \end{equation} For \(\ell\geq1\), let \[ y_\ell:=t-\frac{1}{\ell}(1,\ldots,1), \] and choose \(p_\ell\in\Q^m\) such that \[ y_\ell<p_\ell<t \] coordinatewise. Since \(p_\ell<t\), \[ \liminf_{n\to\infty}F_n(t) \geq \lim_{n\to\infty}F_n(p_\ell) =a(p_\ell). \] Since \(p_\ell>y_\ell\), \eqref{eq:right-regularized-rational-limits} also gives \[ F(y_\ell) = \inf_{\substack{q\in\Q^m\\q>y_\ell}}a(q) \leq a(p_\ell). \] Consequently, \[ \liminf_{n\to\infty}F_n(t) \geq F(y_\ell) \qquad(\ell\geq1). \] Because \(y_\ell\to t\) and \(t\) is a continuity point of \(F\), letting \(\ell\to\infty\) yields \begin{equation} \label{eq:rational-orthant-lower-bound} \liminf_{n\to\infty}F_n(t) \geq F(t). \end{equation} Combining \eqref{eq:rational-orthant-upper-bound} and \eqref{eq:rational-orthant-lower-bound} proves that \(F_n(t)\to F(t)\) at every continuity point of \(F\), and hence \(P_n\weak P\). \end{proof}

\begin{lemma}[Exact-null approximation]
\label{seq_level_lemma}
Let $\psi:\Theta\to\R$ be continuous
in a neighborhood of $\theta_0$ and
differentiable at $\theta_0$, with
$\psi(\theta_0)=0$ and
$c:=\dot\psi_{\theta_0}^T\ne0$.
For every $h\in\R^d$ satisfying
$c^Th=0$, there exists a sequence
$h_n\to h$ such that
\[
\psi\left(
\theta_0+\frac{h_n}{\sqrt n}
\right)=0
\]
for all sufficiently large $n$.
\end{lemma}

\begin{proof}[Proof of Lemma \ref{seq_level_lemma}]
Fix \(h\in\R^d\) such that
\(
c^Th=0
\),
and choose \(v\in\R^d\) such that
\(
c^Tv=1.
\)
For \(|t|\leq1\), define
\[
g_n(t)
:=
\sqrt n\,
\psi\left(
\theta_0+\frac{h+tv}{\sqrt n}
\right).
\]
By differentiability of \(\psi\) at \(\theta_0\),
\[
g_n(t)=t+r_n(t),
\]
where
\begin{equation}
\label{eq:two-sided-null-remainder}
\sup_{|t|\leq1}|r_n(t)|\rightarrow0.
\end{equation}

\noindent Indeed, writing
\[
\psi(\theta_0+u)
=
c^Tu+R(u),
\qquad
\frac{|R(u)|}{\|u\|}\rightarrow0
\quad(u\to0),
\]
we have, for \(|t|\leq1\),
\[
r_n(t)
=
\sqrt n\,
R\left(\frac{h+tv}{\sqrt n}\right).
\]
Since \(\|h+tv\|\leq C\) for \(|t|\leq1\), it follows that
\[
\sup_{|t|\leq1}|r_n(t)|
\leq
C
\sup_{\|u\|\leq C/\sqrt n}
\frac{|R(u)|}{\|u\|}
\rightarrow0.
\]

\noindent Let $
\varepsilon_n:=\sup_{|t|\leq1}|r_n(t)|.$
For all sufficiently large \(n\), \(2\varepsilon_n\leq1\). Moreover,
\[
g_n(-2\varepsilon_n)
\leq
-\varepsilon_n
\leq0, \qquad
g_n(2\varepsilon_n)
\geq
\varepsilon_n
\geq0.
\]
By continuity, there exists
$
t_n\in[-2\varepsilon_n,2\varepsilon_n]
$
such that
$
g_n(t_n)=0.
$
Set
$
h_n:=h+t_nv.$
Then \(h_n\to h\), and
\begin{equation}
\label{eq:exact-null-two-sided}
\psi\left(\theta_0+\frac{h_n}{\sqrt n}\right)=0.
\end{equation}

\end{proof}

\section{Trace norm bound under local perturbations}

\begin{lemma}[Trace-norm bound for local product states]
\label{lem:local-trace-norm-bound}
Let $\Theta\subset\mathbb R^d$ be open and let
$\{\rho_\vartheta:\vartheta\in\Theta\}$ be a quantum statistical model
on a finite-dimensional Hilbert space $\mathcal H$. Assume that
$\vartheta\mapsto\rho_\vartheta$ is continuously differentiable in
trace norm. Let \(L_j(\vartheta):=L_{\vartheta,j}\) denote the SLDs defined in
\eqref{eq:general-sld}, and let
\[
G_\vartheta
=
\left(
G_{\vartheta,ij}
\right)_{i,j=1}^d
\]
be the SLD Fisher information matrix defined in
\eqref{eq:general-sld-information}. Assume that
\(\vartheta\mapsto G_\vartheta\) is continuous.
For \(h\in\mathbb R^d\), write
\[
L_h(\vartheta)
:=
\sum_{j=1}^d h_jL_j(\vartheta),
\qquad
\rho_\vartheta^{(n)}
:=
\rho_\vartheta^{\otimes n}.
\]

For every $\theta\in\Theta$, $h\in\mathbb R^d$, and integer $n\ge1$
such that $\theta+t h/\sqrt n\in\Theta$ for all $t\in[0,1]$,
\begin{equation}
 \begin{split}
 \left\|
   \rho_{\theta+h/\sqrt n}^{(n)}
   -\rho_\theta^{(n)}
 \right\|_1
 &\le
 \int_0^1
 \sqrt{
   h^\top G_{\theta+t h/\sqrt n}h
 }\,dt\\
 &\le
 \sup_{0\le t\le1}
 \sqrt{
   h^\top G_{\theta+t h/\sqrt n}h
 }.
 \end{split}
 \label{eq:local-trace-norm-bound}
\end{equation}
In particular, for fixed $\theta$ and $h$,
\begin{equation}
 \limsup_{n\to\infty}
 \left\|
   \rho_{\theta+h/\sqrt n}^{(n)}
   -\rho_\theta^{(n)}
 \right\|_1
 \le
 \sqrt{h^\top G_\theta h}.
 \label{eq:local-trace-norm-limsup}
\end{equation}
Since $h$ is real, $h^\top J_\vartheta h=h^\top G_\vartheta h$,
so $G_\vartheta$ may equivalently be replaced by $J_\vartheta$
inside these quadratic forms.
\end{lemma}

\begin{proof}
Taking the trace in the SLD equation and differentiating
$\operatorname{Tr}(\rho_\vartheta)=1$ gives
\[
 \operatorname{Tr}\!\left(
   \rho_\vartheta L_j(\vartheta)
 \right)
 =
 \operatorname{Tr}(\partial_j\rho_\vartheta)
 =0.
\]
Consequently,
\begin{equation}
 \operatorname{Tr}\!\left(
   \rho_\vartheta L_h(\vartheta)
 \right)=0,
 \qquad
 \operatorname{Tr}\!\left(
   \rho_\vartheta L_h(\vartheta)^2
 \right)
 =
 h^\top G_\vartheta h
 =
 h^\top J_\vartheta h.
 \label{eq:directional-SLD-moments}
\end{equation}
The second identity follows by expanding $L_h(\vartheta)^2$.
In particular, $G_\vartheta$ is real symmetric and positive
semidefinite.

We first give a trace-norm estimate. For any density matrix
$\sigma$ and self-adjoint matrix $B$, the triangle inequality,
invariance of the trace norm under taking adjoints, and
$\|CD\|_1\le\|C\|_2\|D\|_2$ give
\begin{equation}
 \begin{split}
 \frac12\|\sigma B+B\sigma\|_1
 &\le \|\sigma B\|_1\\
 &=\|\sigma^{1/2}(\sigma^{1/2}B)\|_1\\
 &\le
 \|\sigma^{1/2}\|_2\|\sigma^{1/2}B\|_2
 =
 \sqrt{\operatorname{Tr}(\sigma B^2)},
 \end{split}
 \label{eq:SLD-trace-norm-estimate}
\end{equation}
where $\|C\|_2:=\{\operatorname{Tr}(C^*C)\}^{1/2}$
is the Hilbert--Schmidt norm.

Fix $\theta,h,n$ as in the statement, and put
$\theta_t=\theta+t h/\sqrt n$ and
$\sigma_t=\rho_{\theta_t}^{\otimes n}$.
Continuous differentiability of the model implies that
$t\mapsto\sigma_t$ is continuously differentiable in trace norm.
The fundamental theorem of calculus therefore gives
\begin{equation}
 \rho_{\theta+h/\sqrt n}^{(n)}-\rho_\theta^{(n)}
 =
 \int_0^1\dot\sigma_t\,dt.
 \label{eq:product-state-path}
\end{equation}

Define the collective directional SLD by
\[
 L_h^{(n)}(\vartheta)
 :=
 \sum_{a=1}^n L_h^{[a]}(\vartheta),
 \qquad
 L_h^{[a]}(\vartheta)
 :=
 I^{\otimes(a-1)}
 \otimes L_h(\vartheta)
 \otimes I^{\otimes(n-a)}.
\]
The chain rule gives
\[
 \frac{d}{dt}\rho_{\theta_t}
 =
 \frac1{2\sqrt n}
 \bigl(
   \rho_{\theta_t}L_h(\theta_t)
   +L_h(\theta_t)\rho_{\theta_t}
 \bigr).
\]
Applying the product rule to the $n$ tensor factors consequently yields
\begin{equation}
 \dot\sigma_t
 =
 \frac1{2\sqrt n}
 \bigl(
   \sigma_tL_h^{(n)}(\theta_t)
   +L_h^{(n)}(\theta_t)\sigma_t
 \bigr).
 \label{eq:collective-SLD-derivative}
\end{equation}
Combining \eqref{eq:product-state-path},
\eqref{eq:collective-SLD-derivative}, and
\eqref{eq:SLD-trace-norm-estimate}, we obtain
\begin{equation}
 \left\|
   \rho_{\theta+h/\sqrt n}^{(n)}
   -\rho_\theta^{(n)}
 \right\|_1
 \le
 \frac1{\sqrt n}
 \int_0^1
 \sqrt{
   \operatorname{Tr}\!\left[
     \sigma_t
     \bigl(L_h^{(n)}(\theta_t)\bigr)^2
   \right]
 }\,dt.
 \label{eq:collective-SLD-intermediate-bound}
\end{equation}

To evaluate the integrand, expand the square of the collective SLD.
Each diagonal term has expectation
$\operatorname{Tr}(\rho_{\theta_t}L_h(\theta_t)^2)$.
For $a\ne b$, the tensor-product structure and
\eqref{eq:directional-SLD-moments} give
\[
 \operatorname{Tr}\!\left[
   \sigma_tL_h^{[a]}(\theta_t)L_h^{[b]}(\theta_t)
 \right]
 =
 \left[
   \operatorname{Tr}\!\left(
     \rho_{\theta_t}L_h(\theta_t)
   \right)
 \right]^2
 =0.
\]
Hence
\[
 \operatorname{Tr}\!\left[
   \sigma_t\bigl(L_h^{(n)}(\theta_t)\bigr)^2
 \right]
 =
 n\,\operatorname{Tr}\!\left(
   \rho_{\theta_t}L_h(\theta_t)^2
 \right)
 =
 n\,h^\top G_{\theta_t}h.
\]
Substituting into \eqref{eq:collective-SLD-intermediate-bound}
proves the integral estimate in \eqref{eq:local-trace-norm-bound}.
The supremum estimate follows because the interval of integration
has length one.

Finally, for fixed $\theta,h$, the path lies in $\Theta$ for all
sufficiently large $n$, and
$\sup_{0\le t\le1}\|\theta_t-\theta\|\le\|h\|/\sqrt n\to0$.
Continuity of $G_\vartheta$ at $\theta$ therefore implies
\[
 \sup_{0\le t\le1}
 \|G_{\theta_t}-G_\theta\|_{\mathrm{op}}
 \rightarrow0,
\]
where $\|\cdot\|_{\mathrm{op}}$ is the matrix operator norm.
It follows that
$\sup_{0\le t\le1}\sqrt{h^\top G_{\theta_t}h}
\to\sqrt{h^\top G_\theta h}$, proving
\eqref{eq:local-trace-norm-limsup}.
\end{proof}

\begin{corollary}[Comparison of moving local parameters]
\label{cor:moving-local-trace-norm}
Under the assumptions of Lemma~\ref{lem:local-trace-norm-bound},
fix $\theta_0\in\Theta$ and let $(h_n)$ and $(k_n)$ be bounded
sequences in $\mathbb R^d$. Put
\[
 \delta_n=h_n-k_n,
 \qquad
 \vartheta_{n,t}
 =
 \theta_0+\frac{k_n+t\delta_n}{\sqrt n},
 \qquad 0\le t\le1.
\]
For all sufficiently large $n$,
\begin{equation}
 \begin{split}
 \left\|
   \rho_{\theta_0+h_n/\sqrt n}^{(n)}
   -\rho_{\theta_0+k_n/\sqrt n}^{(n)}
 \right\|_1
 &\le
 \int_0^1
 \sqrt{
   \delta_n^\top G_{\vartheta_{n,t}}\delta_n
 }\,dt\\
 &\le
 \|\delta_n\|
 \sup_{0\le t\le1}
 \|G_{\vartheta_{n,t}}\|_{\mathrm{op}}^{1/2}.
 \end{split}
 \label{eq:moving-local-trace-norm}
\end{equation}
Consequently, if $h_n-k_n\to0$, then
\begin{equation}
 \left\|
   \rho_{\theta_0+h_n/\sqrt n}^{(n)}
   -\rho_{\theta_0+k_n/\sqrt n}^{(n)}
 \right\|_1
 \rightarrow0.
 \label{eq:moving-local-trace-norm-convergence}
\end{equation}
The convergence of rejection probabilities is uniform over all tests:
\[
 \sup_{0\le T\le I_{\mathcal H^{\otimes n}}}
 \left|
   \operatorname{Tr}\!\left(
     \rho_{\theta_0+h_n/\sqrt n}^{(n)}T
   \right)
   -
   \operatorname{Tr}\!\left(
     \rho_{\theta_0+k_n/\sqrt n}^{(n)}T
   \right)
 \right|
 \rightarrow0.
\]
\end{corollary}

\begin{proof}
Boundedness of $h_n$ and $k_n$ gives
\[
 \sup_{0\le t\le1}\|\vartheta_{n,t}-\theta_0\|
 \le
 \frac{\max\{\|h_n\|,\|k_n\|\}}{\sqrt n}
 \rightarrow0.
\]
Since $\Theta$ is open, the entire segment
$\{\vartheta_{n,t}:0\le t\le1\}$ lies in $\Theta$ for all sufficiently
large $n$. Apply Lemma~\ref{lem:local-trace-norm-bound} with basepoint
$\theta=\theta_0+k_n/\sqrt n$ and local displacement $h=\delta_n$.
This gives the first inequality in
\eqref{eq:moving-local-trace-norm}; the second follows from
\[
 \delta_n^\top G_{\vartheta_{n,t}}\delta_n
 \le
 \|G_{\vartheta_{n,t}}\|_{\mathrm{op}}\|\delta_n\|^2.
\]
Continuity of $G$ at $\theta_0$ implies
\[
 \sup_{0\le t\le1}
 \|G_{\vartheta_{n,t}}-G_{\theta_0}\|_{\mathrm{op}}
 \rightarrow0.
\]
The supremum on the right-hand side of
\eqref{eq:moving-local-trace-norm} is therefore bounded.
If $\delta_n\to0$, this proves
\eqref{eq:moving-local-trace-norm-convergence}.
Finally, for any effect $0\le T\le I$,
\[
 \left|\operatorname{Tr}\bigl((\rho-\sigma)T\bigr)\right|
 \le \|\rho-\sigma\|_1\|T\|_{\mathrm{op}}
 \le \|\rho-\sigma\|_1.
\]
Applying this inequality to the two product states proves the last
assertion.
\end{proof}


\section{Proofs for the adaptive scalar Wald tests}
\label{app:adaptive-wald}

\subsection{Regularity and the one-step identities}
\label{app:wald-identities}

Let $Q_\eta$ be the orthogonal projection onto $\ker\rho_\eta$.
Constant rank and the $C^2$ assumption imply that $Q_\eta$ is $C^2$.
Differentiating $\rho_\eta Q_\eta=0$ gives
$Q_\eta(\partial_j\rho_\eta)Q_\eta=0$.
The $C^2$ family of linear maps on Hermitian matrices
\[
 \mathcal T_\eta(B)
 :=\frac{\rho_\eta B+B\rho_\eta}{2}+Q_\eta BQ_\eta
\]
is invertible: in an eigenbasis of $\rho_\eta$, it multiplies each
matrix entry by $(\lambda_a+\lambda_b)/2$ if at least one eigenvalue
is positive, and by $1$ otherwise. Thus
$L_{\eta,j}^{(0)}:=\mathcal T_\eta^{-1}(\partial_j\rho_\eta)$
is a $C^1$ choice of SLDs. Any other SLD differs from this one only
on the kernel block, which does not affect $J_\eta$ or $G_\eta$.
For any continuous SLD choice as in Subsection~\ref{subsec:wald-model},
$G_\eta$ and $Y_\eta$ are therefore continuous, and
$C_Y:=\sup_{\eta\in K}\|Y_\eta\|_{\mathrm{op}}<\infty$.

All derivatives and bounds used below are therefore ordinary
finite-dimensional ones.

Taking traces in \eqref{eq:general-sld} and differentiating
$\operatorname{Tr}\rho_\eta=1$ gives
$\operatorname{Tr}(\rho_\eta L_{\eta,j})=0$ and hence
$\mu(\eta,\eta)=0$. Holding $\eta$ fixed and differentiating only the
true-state argument, cyclicity of the trace yields
\begin{align*}
 \left.\partial_{\vartheta_j}\mu(\vartheta,\eta)
       \right|_{\vartheta=\eta}
 &=\operatorname{Tr}\{(\partial_j\rho_\eta)Y_\eta\}\\
 &=\frac12\operatorname{Tr}\{\rho_\eta
       (L_{\eta,j}Y_\eta+Y_\eta L_{\eta,j})\}
   =(G_\eta v_\eta)_j=(c_\eta)_j.
\end{align*}
Likewise,
\[
 \operatorname{Tr}(\rho_\eta Y_\eta^2)
 =v_\eta^TG_\eta v_\eta=c_\eta^TG_\eta^{-1}c_\eta.
\]
This proves \eqref{eq:wald-influence-identities}.

To establish \eqref{eq:wald-quadratic-cancellation}, set
$B(\vartheta,\eta)=\psi(\eta)+\mu(\vartheta,\eta)-\psi(\vartheta)$.
The identities just proved give $B(\eta,\eta)=0$ and
$\nabla_\vartheta B(\eta,\eta)=0$. On a closed ball contained in
$\Theta$, both second derivatives of $\rho_\vartheta$ and of
$\psi(\vartheta)$ are bounded. Moreover,
\[
 D_\vartheta^2\mu(\vartheta,\eta)[u,u]
 =\operatorname{Tr}\{D^2\rho_\vartheta[u,u]Y_\eta\},
 \qquad
 |D_\vartheta^2\mu(\vartheta,\eta)[u,u]|
 \le \|D^2\rho_\vartheta[u,u]\|_1\|Y_\eta\|_{\mathrm{op}}.
\]
Hence $|D_\vartheta^2B(\vartheta,\eta)[u,u]|\le C_1\|u\|^2$
uniformly for both parameters in a sufficiently small closed ball.
Writing $\delta=\vartheta-\eta$, Taylor's formula along the segment gives
\[
 B(\vartheta,\eta)
 =\int_0^1(1-t)
     D_\vartheta^2B(\eta+t\delta,\eta)[\delta,\delta]dt,
 \qquad
 |B(\vartheta,\eta)|\le\frac{C_1}{2}\|\delta\|^2.
\]
No derivative of $Y_\eta$ or of its spectral projections is taken in
this expansion.

\subsection{Proof of Proposition~\ref{prop:wald-ic-pilot}}

\begin{proof}[Proof of Proposition~\ref{prop:wald-ic-pilot}]
We verify that the POVM in \eqref{eq:wald-ic-povm} induces a regular
identifiable finite-dimensional classical model, after which the result
follows from standard maximum-likelihood theory.

First, since $\|H_j\|_{\mathrm{op}}\leq1$ and
$\varepsilon\in(0,1)$,
\[
    I\pm\varepsilon H_j
    \succeq
    (1-\varepsilon)I.
\]
Hence every $M_{j,\pm}^{\mathrm{pil}}$ is strictly positive. Moreover,
\[
    M_{j,+}^{\mathrm{pil}}
    +
    M_{j,-}^{\mathrm{pil}}
    =
    \frac{I}{s},
\]
and therefore
\[
    \sum_{j=1}^s
    \left(
        M_{j,+}^{\mathrm{pil}}
        +
        M_{j,-}^{\mathrm{pil}}
    \right)
    =I.
\]
Thus $\{M_{j,\pm}^{\mathrm{pil}}\}$ is a POVM. Furthermore,
\[
    M_{j,+}^{\mathrm{pil}}
    -
    M_{j,-}^{\mathrm{pil}}
    =
    \frac{\varepsilon}{s}H_j.
\]
Since $H_1,\ldots,H_s$ form a real basis of the traceless Hermitian
matrices, the real span of the POVM effects contains $I$ and every
$H_j$, and therefore equals the space of all Hermitian matrices.
Hence the POVM is informationally complete.

For $x\in\mathcal X$, write
\[
    p_\eta(x)
    =
    \operatorname{Tr}
    \left(
        \rho_\eta M_x^{\mathrm{pil}}
    \right).
\]
Strict positivity of the effects gives the uniform lower bound
\[
    p_\eta(x)
    \geq
    \frac{1-\varepsilon}{2s},
    \qquad
    \eta\in K,\quad x\in\mathcal X.
\]
Since $\eta\mapsto\rho_\eta$ is $C^2$ in trace norm, each
$\eta\mapsto p_\eta(x)$ is $C^2$. In particular,
$\log p_\eta(x)$ and its first two derivatives are continuous and
bounded on compact subsets of $\operatorname{int}K$.

The induced classical model is identifiable on $K$. Indeed, suppose
that
$p_\eta(x)=p_\vartheta(x)$
    for every $x\in\mathcal X.$
Then
\[
    \operatorname{Tr}
    \left[
        (\rho_\eta-\rho_\vartheta)M_x^{\mathrm{pil}}
    \right]
    =0
    \qquad
    \text{for every }x.
\]
Since the POVM effects span the Hermitian matrices, this implies
$
    \rho_\eta=\rho_\vartheta.
$
By the assumed injectivity of the parametrization on $K$,
$\eta=\vartheta$.

We next verify nonsingularity of the classical Fisher information.
For $a\in\mathbb R^d$,
\[
\begin{aligned}
    a^TI_\theta^{\mathrm{pil}}a
    &=
    \sum_{x\in\mathcal X}
    \frac{
        \{a^T\nabla p_\theta(x)\}^2
    }{
        p_\theta(x)
    }.
\end{aligned}
\]
If this quantity is zero, strict positivity of all
$p_\theta(x)$ implies
\[
    a^T\nabla p_\theta(x)
    =
    \operatorname{Tr}
    \left[
        (\partial_a\rho_\theta)
        M_x^{\mathrm{pil}}
    \right]
    =0
    \qquad
    \text{for every }x,
\]
where
$\partial_a\rho_\theta
    :=
    \sum_{j=1}^d a_j\partial_j\rho_\theta.
$
Informational completeness then gives
$ \partial_a\rho_\theta=0.$
If
\[
    L_{\theta,a}
    :=
    \sum_{j=1}^d a_jL_{\theta,j},
\]
the SLD equations imply
\[
    \partial_a\rho_\theta
    =
    \frac12
    \left(
        \rho_\theta L_{\theta,a}
        +
        L_{\theta,a}\rho_\theta
    \right).
\]

Multiplying by $L_{\theta,a}$ and taking the trace gives
\[
    a^TG_\theta a
    =
    \operatorname{Tr}
    \left(
        \rho_\theta L_{\theta,a}^2
    \right)
    =\operatorname{Tr}\{(\partial_a\rho_\theta)L_{\theta,a}\}
    =0.
\]

Because $G_\theta\succ0$, we conclude that $a=0$. Thus
$I_\theta^{\mathrm{pil}}\succ0.$
It remains only to invoke standard likelihood theory. For completeness,
let
\[
    M_m(\eta)
    :=
    \frac1m
    \sum_{a=1}^m
    \log p_\eta(W_a),
    \qquad
    M(\eta)
    :=
    \sum_{x\in\mathcal X}
    p_\theta(x)\log p_\eta(x).
\]
Because $\mathcal X$ is finite and
$\inf_{\eta\in K,x}p_\eta(x)>0$, the strong law for the empirical
frequencies gives
\[
    \sup_{\eta\in K}
    |M_m(\eta)-M(\eta)|
    \rightarrow0
    \qquad\text{almost surely}.
\]
Moreover,
\[
    M(\theta)-M(\eta)
    =
    \sum_{x\in\mathcal X}
    p_\theta(x)
    \log\frac{p_\theta(x)}{p_\eta(x)}
    \geq0,
\]
with equality only when $p_\eta=p_\theta$, and hence, by
identifiability, only when $\eta=\theta$. Thus the population
criterion has the unique maximizer $\theta$. Since $K$ is compact,
the consistency theorem for $M$-estimators
\cite[Theorem~5.7]{vaart} yields
\[
    \widetilde\theta_n
    \rightarrow
    \theta
    \qquad\text{in probability}.
\]

Finally, choose a closed ball
$U\subset\operatorname{int}K$ with $\theta$ in its interior.
Consistency implies that $\widetilde\theta_n\in U$ with probability
tending to one. On $U$, the induced finite-outcome model
\[
\mathcal P^{\mathrm{pil}}
=
\{p_\eta:\eta\in U\},
\qquad
p_\eta(x)
=
\operatorname{Tr}(\rho_\eta M_x^{\mathrm{pil}}),
\]
is a regular dominated parametric model with common support $\mathcal X$.

Write $s_\theta(x):=\nabla_\theta\log p_\theta(x)$. We first verify differentiability in quadratic mean at $\theta$. By the preceding positivity bound, there exists $\delta>0$ such that
\[
\inf_{\eta\in U}\inf_{x\in\mathcal X}p_\eta(x)\geq\delta.
\]
Moreover, since $\eta\mapsto\rho_\eta$ is $C^2$ in trace norm, each map
$
\eta\longmapsto p_\eta(x)
$
is $C^2$ on $U$, and consequently,
$
\eta\longmapsto \sqrt{p_\eta(x)}
$
is also $C^2$, and
 Taylor expansion at $\theta$ gives, for $h\to0$,
\[
\sqrt{p_{\theta+h}(x)}
=
\sqrt{p_\theta(x)}
+
\frac12 h^Ts_\theta(x)\sqrt{p_\theta(x)}
+
R_x(h),
\]
where, uniformly over $x\in\mathcal X$,
\[
|R_x(h)|
\leq C\|h\|^2
\]
for all sufficiently small $h$. Since $\mathcal X$ is finite,
\[
\sum_{x\in\mathcal X}R_x(h)^2
=
O(\|h\|^4)
=
o(\|h\|^2).
\]
Hence
\[
\sum_{x\in\mathcal X}
\left[
\sqrt{p_{\theta+h}(x)}
-
\sqrt{p_\theta(x)}
-
\frac12 h^Ts_\theta(x)\sqrt{p_\theta(x)}
\right]^2
=
o(\|h\|^2),
\]
which is precisely differentiability in quadratic mean at $\theta$, with score $s_\theta$.

The score has mean zero,
and its covariance matrix is the classical Fisher information
$
I_\theta^{\mathrm{pil}}\succ0.
$

It remains to verify the local Lipschitz condition on the log likelihood. Since $p_\eta(x)$ is $C^2$ and uniformly bounded away from zero on the compact set $U$,
\[
(\eta,x)\longmapsto
\nabla_\eta\log p_\eta(x)
=
\frac{\nabla_\eta p_\eta(x)}{p_\eta(x)}
\]
is continuous and uniformly bounded. Thus there is a finite constant $C_U$ such that
\[
\sup_{\eta\in U}\sup_{x\in\mathcal X}
\|\nabla_\eta\log p_\eta(x)\|
\leq C_U.
\]
By the mean-value theorem, for every $\eta,\vartheta\in U$ and
$x\in\mathcal X$,
\[
\left|
\log p_\eta(x)-\log p_\vartheta(x)
\right|
\leq
C_U\|\eta-\vartheta\|.
\]
Equivalently, the required Lipschitz envelope may be taken to be the constant function
$
F(x)\equiv C_U
$
and the finiteness of the outcome space gives
\[
\mathbb E_\theta[F(X)^2]
=
C_U^2<\infty.
\]
Thus the local Lipschitz and square-integrable-envelope conditions required in \cite[Theorem~5.39]{vaart} are satisfied. Hence the same theorem gives
\[
    \sqrt{m_n}(\widetilde\theta_n-\theta)
    \weak
    N_d
    \left(
        0,
        (I_\theta^{\mathrm{pil}})^{-1}
    \right).
\]
This proves
\eqref{eq:wald-pilot-mle-clt}. In particular,
\[
    \sqrt{m_n}
    (\widetilde\theta_n-\theta)
    =
    O_{\mathbb P_\theta}(1),
\]
which is \eqref{eq:wald-pilot-rate}.
\end{proof}

\subsection{Well-definedness of the adaptive measurement}
\label{app:wald-measurement}

Since $\eta\mapsto Y_\eta$ is continuous and bounded on $K$, for every
bounded continuous real function $f$ the map $\eta\mapsto f(Y_\eta)$ is
continuous. This follows by uniform polynomial approximation on an
interval containing all the spectra. Taking bounded monotone limits
and using finite-dimensional spectral calculus then shows that
$\eta\mapsto M_\eta(B)$ is Borel measurable for every Borel set $B$.
This justifies choosing the measurement at the random pilot estimate,
even when eigenvalues of $Y_\eta$ coincide.

More explicitly, let $P_n^{\mathrm{pil}}(dw)$ be the pilot POVM,
including any classical processing, and write
$\widetilde\theta_n=\tau_n(w)$ for its measurable estimate map.
The joint law of the pilot and second-stage outcomes is
\begin{equation}
 \operatorname{Tr}\{\rho_\vartheta^{\otimes m_n}P_n^{\mathrm{pil}}(dw)\}
 \prod_{a=1}^{N_n}
       \operatorname{Tr}\{\rho_\vartheta M_{\tau_n(w)}(dz_a)\}.
 \label{eq:wald-joint-law}
\end{equation}
This is a parameter-independent measurement rule on the product
experiment: only the outcome law depends on $\vartheta$.
For a measurable rejection region $D$ in the recorded outcomes, its
quantum effect is obtained by integrating
$P_n^{\mathrm{pil}}(dw)\otimes M_{\tau_n(w)}^{\otimes N_n}(D_w)$,
where $D_w$ is the section of $D$ at $w$. Positivity, normalization, and
countable additivity follow from those of the conditional POVMs; in
finite dimension the integrals may be evaluated entrywise. This
establishes $0\le T_{n,\alpha}^\bullet\le I$. Formula~\eqref{eq:wald-joint-law} also proves
conditional independence of the second-stage outcomes.

\subsection{Proof of Theorem~\ref{thm:wald-local-clt}}

\begin{proof}
Fix $h_n\to h$ and put $\theta_n=\theta_0+h_n/\sqrt n$.
All limits in probability in this proof are under
$\mathbb P_{\theta_n}^{(n)}$ unless another law is specified.

By Lemma \ref{lem:local-trace-norm-bound}, we have
\begin{equation}
 \left\|\rho_{\theta+b/\sqrt m}^{\otimes m}
                  -\rho_\theta^{\otimes m}\right\|_1
 \le\int_0^1\sqrt{b^TG_{\theta+t b/\sqrt m}b}\,dt,
 \label{eq:wald-product-trace-bound}
\end{equation}
valid whenever the segment is contained in $\Theta$.
Apply this bound with $m=m_n$ and
$b=b_n:=\sqrt{m_n/n}\,h_n$. The path is
$\theta_0+t h_n/\sqrt n$, so boundedness of $G$ near $\theta_0$ gives
\begin{equation}
 d_n:=\left\|\rho_{\theta_n}^{\otimes m_n}
                       -\rho_{\theta_0}^{\otimes m_n}\right\|_1
 \le C\sqrt{m_n/n}\,\|h_n\|\rightarrow0.
 \label{eq:wald-pilot-trace-transfer}
\end{equation}
Every event determined by the pilot has an associated effect on its
$m_n$ copies. Therefore, for any such event $B_n$,
\[
 |\mathbb P_{\theta_n}^{(n)}(B_n)
             -\mathbb P_{\theta_0}^{(n)}(B_n)|\le d_n.
\]
Applying the event bound to the tail events of
$\sqrt{m_n}\|\widetilde\theta_n-\theta_0\|$ transfers
\eqref{eq:wald-pilot-rate} from $\theta_0$ to $\theta_n$.
Since $\sqrt{m_n}(\theta_n-\theta_0)=b_n\to0$, we obtain
\begin{equation}
 \sqrt{m_n}(\widetilde\theta_n-\theta_n)
       =O_{\mathbb P_{\theta_n}^{(n)}}(1),
 \qquad \widetilde\theta_n\rightarrow\theta_0.
 \label{eq:wald-local-pilot}
\end{equation}

Set $\mu_n=\mu(\theta_n,\widetilde\theta_n)$ and
$B_n=\psi(\widetilde\theta_n)+\mu_n-\psi(\theta_n)$.
With probability tending to one, both parameters lie in the
neighborhood where \eqref{eq:wald-quadratic-cancellation} holds.
Consequently,
\begin{equation}
 \sqrt n|B_n|
 \le C\sqrt n\|\widetilde\theta_n-\theta_n\|^2
 =O_{\mathbb P_{\theta_n}^{(n)}}(\sqrt n/m_n)=o_{\mathbb P_{\theta_n}^{(n)}}(1).
 \label{eq:wald-negligible-bias}
\end{equation}
The displayed bound is used on that event; its complement has
probability tending to zero, which suffices for the probability
conclusion. No moment bound on the pilot's exceptional event is needed.

Let $\sigma_n^2=\nu(\theta_n,\widetilde\theta_n)$.
Continuity of \eqref{eq:wald-two-arguments} and
\eqref{eq:wald-local-pilot} imply
\begin{equation}
 \mu_n\rightarrow0,
 \qquad \sigma_n^2\rightarrow
       \nu(\theta_0,\theta_0)=V_0.
 \label{eq:wald-conditional-variance}
\end{equation}
Conditional on $\mathcal F_n^{\mathrm{pil}}$, the centered variables
$\xi_{n,a}=Z_{n,a}-\mu_n$ are independent and identically distributed,
with mean zero, variance $\sigma_n^2$, and absolute value at most
$2C_Y$. For fixed $t\in\mathbb R$, Taylor's formula for the scalar
exponential gives
\[
 \mathbb E_{\theta_n}^{(n)}
   [e^{it\xi_{n,1}/\sqrt{N_n}}\mid\mathcal F_n^{\mathrm{pil}}]
 =1-\frac{t^2\sigma_n^2}{2N_n}+R_{n,t},
 \qquad
 |R_{n,t}|\le\frac{|t|^3(2C_Y)^3}{6N_n^{3/2}}.
\]
Thus conditional independence and the expansion
$\log(1+z)=z+O(|z|^2)$ near zero imply
\begin{equation}
 \begin{split}
 &\mathbb E_{\theta_n}^{(n)}
 \left[\left.
 e^{it\sqrt{N_n}(\overline Z_n-\mu_n)}
 \right|\mathcal F_n^{\mathrm{pil}}\right]\\
 &\hspace{12mm}
 =\left(1-\frac{t^2\sigma_n^2}{2N_n}+R_{n,t}\right)^{N_n}
 \rightarrow e^{-t^2V_0/2}\quad\text{in probability}.
 \end{split}
 \label{eq:wald-conditional-cf}
\end{equation}
The remainder in the logarithmic expansion is uniform over pilot
outcomes, since $\sigma_n^2\le4C_Y^2$.
Conditional characteristic functions have modulus at most one;
convergence in probability to the constant in
\eqref{eq:wald-conditional-cf} therefore also holds in mean absolute
value. Taking expectations and applying the characteristic-function
continuity theorem yields
\begin{equation}
 \sqrt{N_n}(\overline Z_n-\mu_n)\weak N(0,V_0).
 \label{eq:wald-stage2-clt}
\end{equation}
Since $N_n/n\to1$, the same limit holds with $\sqrt n$.

Finally the exact identity
\[
 \widehat\psi_n-\psi(\theta_n)
 =B_n+(\overline Z_n-\mu_n)
\]
and \eqref{eq:wald-negligible-bias}--\eqref{eq:wald-stage2-clt}
prove the centered normal limit in
\eqref{eq:wald-centered-local-clt}. Continuity of
$V_{\mathrm{eff}}$ and \eqref{eq:wald-local-pilot} imply
$\widehat V_n\to V_0>0$; the lower bound $n^{-1}$ is irrelevant with
probability tending to one. Finally, if $\psi(\theta_0)=0$,
differentiability gives
\[
 \sqrt n\psi(\theta_n)=c_0^Th_n+o(1)\rightarrow c_0^Th.
\]
Adding this shift and applying Slutsky's theorem proves
\eqref{eq:wald-shifted-local-clt}.
\end{proof}

\subsection{Proof of Corollary~\ref{cor:wald-local-power}}

\begin{proof}
For every $h_n\to h$, Theorem~\ref{thm:wald-local-clt} gives the normal
limit of $W_n$ under $\theta_0+h_n/\sqrt n$.
Therefore the one- and two-sided rejection probabilities converge to
$\beta_\alpha^+(c_0^Th;V_0)$ and
$\beta_\alpha^{\pm}(c_0^Th;V_0)$, respectively.

To prove compact uniformity, suppose that the supremum in
\eqref{eq:wald-uniform-local-power} does not converge to zero for a
compact $H$. Then along a subsequence there are $h_n\in H$ for which
the error is bounded below by a fixed positive constant. Passing to a
further subsequence gives $h_n\to h\in H$. The sequential local limit
just proved and continuity of the functions in
\eqref{eq:wald-power-functions} force that error to zero, a
contradiction.

Differentiability at $\theta_0$ also gives, uniformly for $h\in H$,
$\sqrt n\psi(\theta_0+h/\sqrt n)=c_0^Th+o(1)$.
At an actual one-sided local null, this implies $c_0^Th\le o(1)$;
since $\beta_\alpha^+(t;V_0)$ is increasing in $t$ and equals $\alpha$
at zero, \eqref{eq:wald-uniform-local-power} gives the first level bound.
At an actual equality null, $c_0^Th=o(1)$, and continuity of
$\beta_\alpha^{\pm}$ at zero gives the second.
At $h=0$, both limiting powers equal $\alpha$ directly.
Finally, $\beta_\alpha^{\pm}(t;V_0)$ is even and nondecreasing
in $|t|$, with value $\alpha$ at $t=0$, which also gives the stated
local asymptotic unbiasedness for the tangent hypotheses.
\end{proof}
\section{Proofs for weak Schur sampling}
\label{app:wssp}

\begin{proof}[Proof of Lemma~\ref{lem:wss-local-delta-method}]
We first establish the discrete Schur--multinomial approximation used
in the separated-spectrum LAN construction of
\cite[Section~7.5]{Kahn&Guta}.
Let $K$ be a compact neighborhood of $\boldsymbol\alpha_0$ in the
ordered simplex on which every eigenvalue and adjacent gap is at
least $\delta>0$. Write
$P_{n,\alpha}=\operatorname{SW}^n(\alpha)$,
extended by zero to unordered count vectors, and
$Q_{n,\alpha}=\operatorname{Mult}(n,\alpha)$ on
\[
 \mathcal K_n
 :=
 \left\{
 k\in\mathbb Z_{\ge0}^r:
 \sum_i k_i=n
 \right\}.
\]
We claim that
\begin{equation}
 \sup_{\alpha\in K}
 \|P_{n,\alpha}-Q_{n,\alpha}\|_{\mathrm{TV}}
 \rightarrow0.
 \label{eq:spec-schur-multinomial}
\end{equation}
For a partition $\lambda$, the dimension and Schur polynomial
identities give
\[
 P_{n,\alpha}(\lambda)
 =
 \frac{
 n!\prod_{i<j}(\lambda_i-\lambda_j+j-i)
 }{
 \prod_i(\lambda_i+r-i)!
 }
 \frac{
 \det[\alpha_i^{\lambda_j+r-j}]_{i,j=1}^r
 }{
 \prod_{i<j}(\alpha_i-\alpha_j)
 }.
\]
Consider the typical set
\[
 \mathcal A_{n,\alpha}
 :=
 \left\{
 k\in\mathcal K_n:
 \max_i|k_i-n\alpha_i|\le n^{3/4}
 \right\}.
\]
Uniformly over $\alpha\in K$, its elements are positive, strictly
decreasing partitions for all sufficiently large $n$. Moreover,
\[
 Q_{n,\alpha}(\mathcal A_{n,\alpha}^c)
 \le
 2r\exp(-2\sqrt n).
\]
On this set, the identity-permutation term in the determinant is
$\prod_i\alpha_i^{\lambda_i+r-i}$. Every other permutation term is
smaller in absolute value by a factor at most $e^{-cn}$, uniformly
in $\alpha\in K$: the adjacent gaps of $\lambda_i+r-i$ are at least
$n\delta/2$, those of $\log\alpha_i$ have a fixed positive lower
bound, and sorting any nonidentity permutation increases the
corresponding exponent by at least a constant times $n$.
Consequently,
\begin{align*}
 \frac{P_{n,\alpha}(\lambda)}{Q_{n,\alpha}(\lambda)}
 &=
 \frac{
 \prod_{i<j}(\lambda_i-\lambda_j+j-i)
 }{
 \prod_{i<j}(\alpha_i-\alpha_j)
 }
 \prod_i
 \frac{
 \lambda_i!\,\alpha_i^{r-i}
 }{
 (\lambda_i+r-i)!
 }
 \{1+O(e^{-cn})\}
 \\
 &=
 1+O(n^{-1/4}),
 \qquad \lambda\in\mathcal A_{n,\alpha},
\end{align*}
uniformly over $K$. Indeed, each factor in the first product is
$n\{1+O(n^{-1/4})\}$, and the $i$th factorial factor is
$n^{-(r-i)}\{1+O(n^{-1/4})\}$; the powers of $n$ cancel because
$\sum_i(r-i)=r(r-1)/2$.

If the uniform error in this likelihood ratio is $\eta_n\to0$
and $q_n=Q_{n,\alpha}(\mathcal A_{n,\alpha}^c)$, normalization gives
\[
 \sum_{\mathcal A_{n,\alpha}}
 |P_{n,\alpha}-Q_{n,\alpha}|
 \le\eta_n,
 \qquad
 P_{n,\alpha}(\mathcal A_{n,\alpha}^c)
 \le q_n+\eta_n.
\]
Thus
$\|P_{n,\alpha}-Q_{n,\alpha}\|_{\mathrm{TV}}\le\eta_n+q_n$,
proving \eqref{eq:spec-schur-multinomial}.

Eventually $\boldsymbol\alpha_n\in K$. Represent a multinomial
vector as $N_n=\sum_{k=1}^n e_{I_{n,k}}$, where the $I_{n,k}$ are
independent within each row and
$\mathbb P(I_{n,k}=i)=\alpha_{n,i}$.
The centered summands are uniformly bounded and their covariance
matrices converge to $C_0$. The triangular-array central limit theorem
and Cram\'er--Wold therefore imply
\[
 (N_n-n\boldsymbol\alpha_n)/\sqrt n
 \weak N_r(0,C_0).
\]
For scalar directions with zero limiting variance, convergence to zero
follows from Chebyshev's inequality. Total variation contracts under
measurable maps, so \eqref{eq:spec-schur-multinomial} transfers this
limit to WSS and proves \eqref{eq:spec-wss-vector-clt}.
The deterministic centering shift gives
\eqref{eq:spec-wss-shift-clt}.

In particular,
$\widehat{\boldsymbol\alpha}_n\to\boldsymbol\alpha_0$ in probability
and
$\sqrt n(\widehat{\boldsymbol\alpha}_n-\boldsymbol\alpha_n)$
is tight. On a convex neighborhood on which $\varphi$ is $C^1$,
\[
 \varphi(\widehat{\boldsymbol\alpha}_n)
 -
 \varphi(\boldsymbol\alpha_n)
 =
 b_0^T(\widehat{\boldsymbol\alpha}_n-\boldsymbol\alpha_n)+R_n,
\]
where
\[
 |R_n|
 \le
 \|\widehat{\boldsymbol\alpha}_n-\boldsymbol\alpha_n\|
 \sup_{0\le t\le1}
 \left\|
 \nabla\varphi\bigl(
 \boldsymbol\alpha_n
 +
 t(\widehat{\boldsymbol\alpha}_n-\boldsymbol\alpha_n)
 \bigr)
 -
 b_0
 \right\|.
\]
The supremum is $o_{\mathbb P}(1)$, and the event on which this
expansion applies has probability tending to one. Hence
$\sqrt nR_n=o_{\mathbb P}(1)$, and
\eqref{eq:spec-wss-delta-clt} follows from
\eqref{eq:spec-wss-vector-clt} and Slutsky's theorem.
\end{proof}

\begin{proof}[Proof of Theorem~\ref{thm:qlan-wss-comparison}]
\begin{enumerate}[
  label=\textbf{\ensuremath{\langle\arabic*\rangle}},
  leftmargin=3.5em,
  itemsep=0.8em
]

\item
By differentiability of \(\psi\) at \(\theta_0\),
\begin{equation}
\label{eq:spectral-functional-local-expansion}
\sqrt n
\left[
\psi\left(\theta_0+\frac{h}{\sqrt n}\right)
-
\psi(\theta_0)
\right]
=
c_0^Th+o(1).
\end{equation}

\item
By the q-LAN convergence under local alternatives,
\[
\left(
X^{(n)},
\rho_{\theta_0+h/\sqrt n}^{\otimes n}
\right)
\weak
N(h,\Sigma_{\theta_0}).
\]
Therefore,
\[
c_0^TX^{(n)}
\weak
N(c_0^Th,V_{\mathrm{eff}}).
\]
Combining this convergence with
\eqref{eq:spectral-functional-local-expansion} proves
\eqref{eq:qlan-functional-local-clt}.

\item

By Lemma~\ref{lem:wss-local-delta-method},
applied with
$\boldsymbol\alpha_n
=\boldsymbol\alpha_{\theta_0+h/\sqrt n}$,

\[
\sqrt n
\left[
\varphi\left(\frac{\lambda^{(n)}}n\right)
-
\varphi(\boldsymbol{\alpha}_{\theta_0+h/\sqrt n})
\right]
\weak
N(0,V_{\mathrm{WSS}}).
\]
This is exactly \eqref{eq:wss-functional-local-clt}.

\item
The rejection region in \eqref{eq:qlan-spectral-test}, together with
the convergence obtained in Step~\(\langle2\rangle\), gives
\[
\lim_{n\to\infty}
\Tr\left[
\rho_{\theta_0+h/\sqrt n}^{\otimes n}T_{n,1}
\right]
=
1-\Phi\left(
z_\gamma-\frac{\kappa(h)}{\sqrt{V_{\mathrm{eff}}}}
\right).
\]
The boundary of the limiting rejection region has zero Gaussian
probability.  The same argument applied to
\eqref{eq:wss-functional-local-clt} proves
\eqref{eq:wss-local-power-spectral}.

\item
Fix a direction \(h\) with \(\kappa(h)>0\), a significance level
\(\gamma\in(0,1)\), and a desired power
\[
\pi\in(\gamma,1).
\]
Write
\[
z_\pi:=\Phi^{-1}(1-\pi).
\]
At the alternative \(\theta_0+th\), a test with asymptotic variance
\(V_j\), \(j\in\{1,2\}\), has asymptotic power
\[
1-\Phi\left(
z_\gamma-\frac{\sqrt n\,t\kappa(h)}{\sqrt{V_j}}
\right).
\]
Equating this expression to \(\pi\) gives
\[
\frac{\sqrt n\,t\kappa(h)}{\sqrt{V_j}}
=
z_\gamma-z_\pi.
\]
Consequently, the sample sizes $n_j(t)$ prescribed
by the limiting Gaussian power satisfy
\begin{equation}
\label{eq:sample-size-from-local-power}
n_j(t)
\sim
\frac{
V_j(z_\gamma-z_\pi)^2
}{
t^2\kappa(h)^2
},
\qquad
t\downarrow0.
\end{equation}
Taking the quotient in the orientation used by the Pitman efficiency
of the first test relative to the second gives
\[
\operatorname{ARE}_{\mathrm P}(T_1,T_2)
=
\lim_{t\downarrow0}
\frac{n_2(t)}{n_1(t)}
=
\frac{V_{\mathrm{WSS}}}{V_{\mathrm{eff}}}.
\]
This proves \eqref{eq:qlan-wss-are}.

\end{enumerate}
\end{proof}

\begin{proof}[Proof of Corollary~\ref{cor:wss-efficient}]

\vspace{0.2in}

 \textbf{The difference $V_{\mathrm{WSS}}-V_{\mathrm{sp}}$}

\vspace{0.2in}

Let
\(
D_0:=\operatorname{diag}(\boldsymbol\alpha_0).
\)
Then it is easy to check that
\(
w_0:=D_0^{-1}A_0I_{\mathrm{sp},0}^{+}c_0
\)
is a feasible solution to the optimization problem
\eqref{eq:wssp-optimization}, and that
\(
w_0^TD_0w_0=V_{\mathrm{sp}}.
\)
For any other feasible $w$, write $w=w_0+r$. Since $w$ and $w_0$
satisfy the first constraint in \eqref{eq:wssp-optimization}, we get
$A_0^Tr=0$, and hence
\(
w_0^TD_0r=c_0^TI_{\mathrm{sp},0}^{+}A_0^Tr=0.
\)
Expanding the quadratic form $w^TD_0w$ gives
\[
w^TD_0w
=w_0^TD_0w_0+r^TD_0r
=V_{\mathrm{sp}}+r^TD_0r
\geq V_{\mathrm{sp}}.
\]
Since $D_0\succ0$, equality with $V_{\mathrm{sp}}$ occurs precisely
when $r=0$, or $w=w_0$. This establishes the uniqueness of the
minimizer.
Now, coming to $V_{\mathrm{WSS}}$, define
\[
\bar b_0=b_0-(\boldsymbol\alpha_0^Tb_0)\mathbf1
\]
It is easy to check that $\bar b_0$ is a feasible solution to the
optimization problem \eqref{eq:wssp-optimization}. Observe that
\[
\bar b_0^TD_0\bar b_0=b_0^TC_0b_0=V_{\mathrm{WSS}}
\]
This implies that $V_{\mathrm{WSS}}\geq V_{\mathrm{sp}}$, with
equality if and only if $\bar b_0=w_0$.
Now note that
\[
\bar b_0=w_0
\implies
\bar b_0\in\operatorname{Ran}(D_0^{-1}A_0)
\]
We show the converse. Suppose
$\bar b_0=D_0^{-1}A_0v$. The constraint
$A_0^T\bar b_0=c_0$ gives
$I_{\mathrm{sp},0}v=c_0$.
This implies
$$v - I_{\mathrm{sp},0}^{+}c_0\in \operatorname{ker}(I_{\mathrm{sp},0})=\operatorname{ker}(A_0)$$
and hence $D_0^{-1}A_0(v - I_{\mathrm{sp},0}^{+}c_0)=0 \implies \bar b_0=w_0$. Therefore,
\begin{equation}
\bar b_0=w_0
\iff
\bar b_0\in\operatorname{Ran}(D_0^{-1}A_0)
\label{eq:equality_cond}
\end{equation}
Now observe that
$D_0\bar b_0
=D_0b_0-\boldsymbol\alpha_0(\boldsymbol\alpha_0^Tb_0)
=C_0b_0$.
Hence the condition \eqref{eq:equality_cond} becomes
\[
C_0b_0\in\operatorname{Ran}(A_0)
\]
Thus equality holds if and only if the covariance-weighted derivative $C_0b_0$ lies in the space of first-order spectral changes $\operatorname{Ran}(A_0)$. This proves \eqref{eq:wss-sp-equality}.\\

\noindent\textbf{The difference $V_{\mathrm{sp}}-V_{\mathrm{eff}}$.}\\

Consider the Hilbert space
\begin{equation}
\mathcal H_0
:=
\{H=H^*:\operatorname{Tr}(\rho_0H)=0\},
\label{eq:space}
\end{equation}
with the inner product
\[
\langle H,K\rangle_{\rho_0}
=
\frac12\operatorname{Tr}\{\rho_0(HK+KH)\}.
\]
Its associated squared norm is
$\|H\|_{\rho_0}^2=\operatorname{Tr}(\rho_0H^2)$.
Then note that
$\mathcal T_0=\operatorname{span}_{\mathbb R}
\{L_{0,1},\ldots,L_{0,d}\}$
is a subspace of this space. Recall the SLD information matrix
$(G_0)_{jk}=\langle L_{0,j},L_{0,k}\rangle_{\rho_0}$.

\medskip
\noindent We now show that $V_{\mathrm{eff}}$ is the squared
norm of an observable in \eqref{eq:space}.
Define $Y_0=\sum_{j=1}^d(v_0)_jL_{0,j}$,
where $v_0=G_0^{-1}c_0$. Consequently,
\[\langle L_{0,j},Y_0\rangle_{\rho_0}=(c_0)_j , \qquad 
\|Y_0\|_{\rho_0}^2
=v_0^TG_0v_0
=c_0^TG_0^{-1}c_0
=V_{\mathrm{eff}}.\]
Similarly, $V_{\mathrm{sp}}$ is the squared
norm of an observable in \eqref{eq:space}.
Define $H_{\mathrm{sp},0}
=U_0\operatorname{diag}(w_0)U_0^*$.
Then observe that
$\operatorname{Tr}(\rho_0H_{\mathrm{sp},0})
=\boldsymbol\alpha_0^Tw_0=0$,
so $H_{\mathrm{sp},0}\in\mathcal H_0$.
Furthermore,
\[
\begin{aligned}
\|H_{\mathrm{sp},0}\|_{\rho_0}^2
&=\operatorname{Tr}(\rho_0H_{\mathrm{sp},0}^2)\\
&=\operatorname{Tr}\bigl(
U_0D_0U_0^*
U_0\operatorname{diag}(w_0)^2U_0^*
\bigr)\\
&=\operatorname{Tr}\bigl(
U_0D_0\operatorname{diag}(w_0)^2U_0^*
\bigr)\\
&=\operatorname{Tr}(D_0w_0w_0^T)\\
&=w_0^TD_0w_0\\
&=V_{\mathrm{sp}}.
\end{aligned}
\]
We now show that $H_{\mathrm{sp},0}-Y_0$
and $Y_0$ are orthogonal.
Indeed, observe that
\[
(U_0^*L_{0,j}U_0)_{ii}
=\frac{(A_0)_{ij}}{\alpha_{0,i}}.
\]
Therefore,
\[
\begin{aligned}
\langle L_{0,j},H_{\mathrm{sp},0}\rangle_{\rho_0}
&=\sum_{i=1}^r
\alpha_{0,i}
(U_0^*L_{0,j}U_0)_{ii}(w_0)_i\\
&=\sum_{i=1}^r(A_0)_{ij}(w_0)_i\\
&=(A_0^Tw_0)_j\\
&=(c_0)_j\\
&=\langle L_{0,j},Y_0\rangle_{\rho_0}.
\end{aligned}
\]
Hence
\begin{equation}
\langle L_{0,j},
H_{\mathrm{sp},0}-Y_0\rangle_{\rho_0}=0,
\qquad j=1,\ldots,d.
\label{eq:projection}
\end{equation}
Since $Y_0=\sum_{j=1}^d(v_0)_jL_{0,j}$,
\[
\begin{aligned}
\langle H_{\mathrm{sp},0}-Y_0,Y_0\rangle_{\rho_0}
&=\sum_{j=1}^d(v_0)_j
\langle H_{\mathrm{sp},0}-Y_0,
L_{0,j}\rangle_{\rho_0}\\
&=0.
\end{aligned}
\]

\medskip
\noindent\textbf{Then
$V_{\mathrm{sp}}-V_{\mathrm{eff}}
=\|H_{\mathrm{sp},0}-Y_0\|_{\rho_0}^2\geq0$.}
Indeed,
\begin{equation}
\begin{aligned}
\|H_{\mathrm{sp},0}-Y_0\|_{\rho_0}^2
&=\|H_{\mathrm{sp},0}\|_{\rho_0}^2
-\|Y_0\|_{\rho_0}^2\\
&=V_{\mathrm{sp}}-V_{\mathrm{eff}}
\geq0.
\end{aligned}
\label{eq:sp-eff-gap}
\end{equation}
Equality holds if and only if $H_{\mathrm{sp},0}=Y_0$.
Since $H_{\mathrm{sp},0}$ commutes with $\rho_0$, we have
\[
H_{\mathrm{sp},0}=Y_0
\implies [Y_0,\rho_0]=0.
\]
We prove the converse. If $Y_0$ commutes with $\rho_0$,
then since $\rho_0$ has simple spectrum, $Y_0$ is diagonal
in the eigenbasis of $\rho_0$, so that
\[
Y_0=U_0\operatorname{diag}(y)U_0^*
\]
for some $y\in\mathbb R^r$.
We show that $y$ is feasible under the optimization problem
\eqref{eq:wssp-optimization}. First,
\[
\begin{aligned}
\boldsymbol\alpha_0^Ty
&=\operatorname{Tr}(\rho_0Y_0)\\
&=0.
\end{aligned}
\]
Also,
\[
\begin{aligned}
(A_0^Ty)_j
&=\sum_{i=1}^r(A_0)_{ij}y_i\\
&=\sum_{i=1}^r
\alpha_{0,i}
(U_0^*L_{0,j}U_0)_{ii}y_i\\
&=\langle L_{0,j},Y_0\rangle_{\rho_0}\\
&=(c_0)_j.
\end{aligned}
\]
Thus $A_0^Ty=c_0$, and $y$ is feasible.
But the value of the objective function at $y$ is
\[
\begin{aligned}
y^TD_0y
&=\operatorname{Tr}(\rho_0Y_0^2)\\
&=\|Y_0\|_{\rho_0}^2\\
&=V_{\mathrm{eff}}.
\end{aligned}
\]
Hence
\[
V_{\mathrm{eff}}\geq V_{\mathrm{sp}}.
\]
We also had earlier, from \eqref{eq:sp-eff-gap}, that
$V_{\mathrm{sp}}\geq V_{\mathrm{eff}}$.
Hence $V_{\mathrm{sp}}=V_{\mathrm{eff}}$.
By the uniqueness of the minimizer, $y=w_0$.
Conjugating with $U_0$ implies $Y_0=H_{\mathrm{sp},0}$.
Hence we have proved
\[
H_{\mathrm{sp},0}=Y_0
\iff [Y_0,\rho_0]=0.
\]
This proves that
$V_{\mathrm{sp}}=V_{\mathrm{eff}}
\iff [Y_0,\rho_0]=0$,
which is \eqref{eq:sp-eff-equality}.
On the finite-dimensional real Hilbert space of centered Hermitian
matrices, the commutation operator satisfies
\begin{equation}
 \langle\mathcal D_{\rho_0}H,K\rangle_{\rho_0}
 =-\langle H,\mathcal D_{\rho_0}K\rangle_{\rho_0}.
 \label{eq:wssp-D-skew}
\end{equation}
This follows directly in the reference eigenbasis: on each real
plane formed by the real and imaginary parts of $H_{ij}$, the map
is a scalar multiple of a ninety-degree rotation, while the SLD
inner product is a positive multiple of the Euclidean inner product.
It vanishes on diagonal matrices. Since the spectrum is simple,
\[
 \ker\mathcal D_{\rho_0}
 =\{H=H^*:[H,\rho_0]=0\}.
\]
If $\mathcal T_0$ is invariant under $\mathcal D_{\rho_0}$,
\eqref{eq:wssp-D-skew} implies that $\mathcal T_0^\perp$ is invariant
as well. Hence its orthogonal projection commutes with the operator:
\[
 \mathcal D_{\rho_0}\Pi_{\mathcal T_0}
 =\Pi_{\mathcal T_0}\mathcal D_{\rho_0}.
\]
Now note that 
\begin{itemize}
    \item The observable $H_{\mathrm{sp},0}$ commutes with $\rho_0$.
    \item $\Pi_{\mathcal T_0}H_{\mathrm{sp},0} = Y_0$ by \eqref{eq:projection}
\end{itemize}    
Therefore
\[
 \mathcal D_{\rho_0}Y_0
 =\mathcal D_{\rho_0}\Pi_{\mathcal T_0}H_{\mathrm{sp},0}
 =\Pi_{\mathcal T_0}\mathcal D_{\rho_0}H_{\mathrm{sp},0}=0.
\]
It follows that $[Y_0,\rho_0]=0$. The equality
$V_{\mathrm{sp}}=V_{\mathrm{eff}}$ now follows.
Finally, if $\operatorname{rank}(A_0)=r-1$, then
since $\mathbf1^TA_0=0$,
\[
\operatorname{Ran}(A_0)
=
\{x\in\R^r:\mathbf1^Tx=0\}.
\]
Since $\mathbf1^TC_0=0$, we have
$C_0b_0\in\operatorname{Ran}(A_0)$.
Thus $V_{\mathrm{WSS}}=V_{\mathrm{sp}}$
by \eqref{eq:wss-sp-equality}.
If the SLD space is also $\mathcal D_{\rho_0}$-invariant,
then $V_{\mathrm{sp}}=V_{\mathrm{eff}}$,
and hence
\[
V_{\mathrm{WSS}}=V_{\mathrm{sp}}=V_{\mathrm{eff}}.
\]

\end{proof}

\end{document}